\documentclass[a4paper]{amsart}
\usepackage{graphicx} % Required for inserting images
\usepackage{a4wide}
\usepackage{amsmath, comment}
\usepackage[figuresright]{rotating}
\usepackage{amssymb,latexsym,graphics,subfigure,caption}
\usepackage{amsthm}
\newtheorem{theorem}{Theorem}[section]
\newtheorem{remark}{Remark}[section]
\newtheorem{proposition}[theorem]{Proposition} 
\newtheorem{lemma}[theorem]{Lemma}
\newtheorem{assumption}[theorem]{Assumption}
\usepackage{bbm}
\usepackage{cite}[compress,verbatim]
\allowdisplaybreaks
\title[Resonance-based stochastic numerical integrators]{Resonance-based integrators for stochastic Schrödinger equations. Convergence and long-time error bounds.}
\author{Stefano Di Giovacchino}
\address{Dipartimento di Ingegneria e Scienze dell'Informazione e Matematica, Università degli Studi dell'Aquila, Via Vetoio 67100, L'Aquila, Italy.}
\email{stefano.digiovacchino@univaq.it}
\thanks{The author is a member of the INdAM Research group GNCS and his work was supported by PRIN-MUR 2022 project 20229P2HEA “Stochastic numerical modelling for sustainable innovation” (CUP: E53C24002280006), granted by MUR within the scrolling of the final rankings of the PRIN 2022 call. The author thanks Prof. K. Schratz (Sorbonne, Paris) for hosting him in 2024 and for the valuable discussions. Also, the author is grateful to Dr. J. Cui (Hong Kong Polytechnic University) for all the suggestions. Finally, the author thanks Dr. G. Maierhofer (Univ. of Oxford) for sharing his material and for the interaction.}
\date{}
\begin{document}
\begin{abstract}
In this work, we develop and analyze resonance-based low-regularity numerical integrators for stochastic Schrödinger equations driven by additive $Q$-Wiener noise, covering both the linear equation with rough potential and the cubic nonlinear case.
For the linear equation, we prove strong and almost sure convergence under low regularity assumptions, achieving first-order convergence in $H^\sigma$ for solutions in $H^{\sigma+1}$, thereby improving on the classical $H^{\sigma+2}$ requirement. In a regime of $\mathcal{O}(\varepsilon^2)$ potentials and $\mathcal{O}(\varepsilon)$ noise, $\varepsilon \ll1 $, we establish uniform moment bounds for the exact solution up to times of length $\mathcal{O}(\varepsilon^{-2})$, and construct a non-resonant low-regularity scheme exhibiting uniform long-time error bounds of size $\mathcal{O}(\varepsilon^2 \tau)$, being $\tau>0$ the time step of the numerical integration, up to times of length $\mathcal{O}(\varepsilon^{-2}).$
For the cubic equation, we introduce both resonant and non-resonant integrators and derive analogous pathwise convergence results at low-regularity. In the weakly nonlinear stochastic regime, i.e., $\mathcal{O}(\varepsilon^2)$ nonlinearity and $\mathcal{O}(\varepsilon)$ noise, we prove long-time pathwise error bounds of size $\mathcal{O}(\varepsilon^2 \tau^\delta)$, for any $\delta<1$, up to times of length $\mathcal{O}(\varepsilon^{-2})$.
For both the equations, the long-term error analysis relies on a novel nontrivial extension of the {\em regularity-compensation oscillation technique} (RCO) to stochastic dispersive equations, allowing to get improved error estimates overcoming the loss of temporal regularity induced by stochastic convolutions. In particular, the long-term errors reduce by a factor of $\varepsilon^2$, compared to other classical and resonant low-regularity integrators.
To the best of our knowledge, this is the first work establishing long-time error bounds for low-regularity integrators applied to stochastic dispersive equations. The theoretical results will be finally confirmed by selected numerical experiments.
\end{abstract}
\maketitle
 
\section{Introduction}
In this work, we consider the stochastic Schr$\ddot{\rm o}$dinger equation (SSE) with additive noise, assuming the following formulation
\begin{equation}
    \label{SSE}
    \dot{\xi}(t)={\bf i}\Delta \xi(t)-{\bf i}F(\xi(t))+Q^{\frac{1}{2}}\dot{W}(t), \quad t\in[0,T],
\end{equation}
equipped with $\xi(0)=\xi_0$, where $F$ is in general a nonlinear mapping, $Q$ is the Covariance operator of a cylindrical complex-valued Wiener process $W(t)$ respect to a normal filtration $\{\mathcal{F}_t\}_{t\ge 0}$ on a filtered probability space 
$\left(\Omega, \mathcal{F}, \mathbb{P}, \{\mathcal{F}_t\}_{t\ge 0}\right)$. Further regularity on the nonlinear mapping $F$, on the operator $Q$ and of the Wiener process $W$ will be provided in the sequel.  

Nonlinear Schrödinger equations, and in particular, the cubic nonlinear equation, are widely spread in several fields, such as, plasma physics and hydrodynamics, since they represent effective mathematical models to describe many important phenomena. A possibility to incorporate stochastic effects in  physical models, described by Schrödinger equations, is given by model \eqref{SSE}. Indeed, in the scientific literature, several authors have addressed their studies to \eqref{SSE}, both in It$\hat{\rm o}$ and Stratonovich formulation; see, for instance, \cite{deb3,deb2,deb4, deb0,deb5,cui2,mora,hong,hong_wang} and reference therein.

The temporal discretization of (SSE) has captured much attention among numerical analysts within the last years; prior works are given by \cite{deb7,deb0}, in which the authors provided a semi-discrete temporal discretization to (SSE), in the case of $F=0$, $F$ globally Lipschitz and in the case of cubic nonlinearity. They showed convergence in various topology and also in the weak sense, both with additive and Stratonovich noise. In \cite{an_co}, the authors provided an explicit exponential integrator for \eqref{SSE}, in the linear case and also when $F(\xi)=V(x)\xi$, with $V$ smooth potential, and carried on a strong convergence analysis and showed some conservative properties of the provided integrators. Also the multiplicative linear case was addressed. Convergence and geometric properties of splitting type schemes for (SSE), both under additive and multiplicative noise, have been addressed, e.g., in \cite{bre-co,cui,jia,liu,liu2,hong,hong_wang}, while, for the theory regarding the numerical preservation of invariant measures for (SSE), we refer, e.g., to the monograph \cite{hong_wang} and the references therein.

However, in all the aforementioned paper, the provided theoretical analysis is limited to short time scales.

Recently, the long-time behaviour of solutions to semilinear evolutive equations, e.g. Schrödinger equations, has captured great interest and several works have been provided to study the effects of small nonlinearities and/or small initial data to the long-time evolution of the exact solutions \cite{bej,carl,fan,fuji,ik,oh,sun}.

The aim of this work is to study convergence and long-term behaviour of resonance-based low-regularity numerical integrators (i.e. integrators constructed exploiting Fourier analysis) for \eqref{SSE}, in the cases $F(\xi)=V(x)\xi$ and $F(\xi)=|\xi|^2 \xi$, when the exact solutions show a boundedness of moments up to large time scales.

In particular, we will design numerical integrators requiring $H^{\sigma+1}$ regularity for showing order one of convergence in $H^\sigma$, in contraposition with $H^{\sigma+2}$ regularity required by standard exponential and splitting type integrators \cite{bre-co}, and exhibiting excellent long-time behaviour in terms of boundedness of the errors.
%In this paper, we address our attention both on the case $F(\xi)=V(x)\xi$ and $F(\xi)=|\xi|^2 \xi$, with the aim of designing resonance-based low regularity numerical integrators in these two cases.  Beyond proving convergence for such integrators, the main focus of the presented work is to show long-term error bounds which make such integrators suitable for the time integration of (SSE) up to long time scales, in the scenario in which the exact solutions to (SSE) satisfies a boundedness of moments up to large time scales.

Although in the deterministic setting convergence and long-term error bounds of low-regularity numerical integrators applied to semilinear dispersive equations have been well-established (see, e.g., \cite{bao,fe_ka,os_ka,ma_ka,fe_ma,ala,bru_ka,kat1,li,kat0} and reference therein, for a non-exhaustive list of contributions on this direction), the application of these techniques to the setting of stochastic dispersive PDEs in nowadays a challenging open research field; recent first works in this area are given by \cite{arm,cao,cui1}; in particular, in \cite{arm}, the authors constructed symplectic low-regularity integrators for the nonlinear (SSE) and showed geometric properties and provided a local error analysis for the proposed schemes. In \cite{cui1}, a low-regularity scheme has been introduced and studied for the nonlinear Schrödinger equation with white noise dispersion. Finally, in \cite{cao}, strong convergence rates of an unfiltered low-regularity method have been investigated for the stochastic wave equation with a globally Lipschitz nonlinearity.

Also in these works, the results are restricted to short time scales.

At the best of our knowledge, there are no works establishing long-term error bounds for numerical integrators, and, in particular, for resonance-based low regularity numerical integrators, for (SSE), and, more in general, for stochastic dispersive equations. 

Following the idea in \cite{os_ka}, firstly, we will design low-regularity integrators by exploiting Fourier analysis, available for periodic square integrable complex-valued functions defined on the one-dimensional torus $\mathbb{T}$; in particular, such integrators will be obtained by integrating exactly the integral arising from the dominant term appearing in the Fourier expansion for the numerical approximation. Also, by exploiting the ideas developed in, e.g., \cite{bao}, non-resonant low-regularity integrators will be designed by integrating exactly the integral for the {\em zero mode term}. Then, in this situation, by introducing a novel extension of the (RCO) technique to the stochastic setting, we will be able to obtain improved long-term error bounds for such numerical solutions.

In the deterministic setting, the (RCO) technique has been exploited and applied to quadratic and cubic Schrödinger equations in \cite{fe_ka}, where the authors provided improved long-term errors bounds for non-resonant low-regularity integrators for the aforementioned equations under the assumption of an initial datum of size $\mathcal{O}(\varepsilon)$, $\varepsilon\ll 1$. In particular, for the cubic equations, they showed an error of size $\mathcal{O}(\varepsilon^2\tau)$ in $H^1$, for smooth solutions, up to time of size $\mathcal{O}(\varepsilon^{-2})$. To achieve these result, they also made use of the following well-known time regularity result for the solution to deterministic Schrödinger equations
\begin{equation}
    \label{intro_0}
    \|\xi(t+\tau)-  \xi(t)\|_{H^\sigma}\lesssim\tau.
\end{equation}

It is now worth mentioning that the present work is the first contribution in which (RCO) technique  is extended to stochastic dispersive equations and its application cannot descend directly from the deterministic setting. Indeed, the direct use of a time regularity result in the stochastic setting does not allow to get the desired first order in $\tau$, since, due to the presence of the stochastic convolution term, it translates into the well-know estimate
\begin{equation}
    \label{intro_1}
    \|\xi(t+\tau)- \xi(t)\|_{L^{2p}(\Omega; H^\sigma)}\lesssim\sqrt{\tau},
\end{equation}
i.e., one looses one-half order in $\tau$, if straightforwardly applies the estimate in \eqref{intro_1}.
Then, this introduces further challenges in the extension of the (RCO) technique in the stochastic scenario that will be treat through this work.

The major contributions of the present paper are now listed.
\begin{itemize}
\item For the case of linear equation with potential, we construct a resonant low-regularity integrator and establish strong and almost sure convergence with order $\gamma \in (0,1]$ in $L^{2p}(\Omega;H^\sigma), \ \sigma\in\mathbb{N}$, for solutions in $H^{\sigma+\gamma}$. In particular, we obtain first-order convergence in $H^\sigma$ under only $H^{\sigma+1}$ regularity, thereby improving the classical requirement of $H^{\sigma+2}$. In the regime of $\mathcal{O}(\varepsilon^2)$ potential and $\mathcal{O}(\varepsilon)$ noise, we prove uniform moment bounds for the exact solution up to times of length $\mathcal{O}(\varepsilon^{-2})$. We then design a non-resonant low-regularity integrator and establish uniform strong error bound in $L^2(\Omega; L^2)$ over the same long-time interval. Notably, the proposed method achieves an error of size $\mathcal{O}(\varepsilon^2 \tau)$, improving existing long-time error estimates by a factor of $\varepsilon^2$.
\item For the cubic equation, we introduce both resonant and non-resonant low-regularity schemes and derive pathwise convergence results at low-regularity. In particular, order $\gamma\in(0,1)$ is obtained in $H^\sigma, \ \sigma\in\mathbb{N}^+$, for solutions in $H^{\sigma+\gamma}$, improving classical regularity assumptions required by standard methods. In the weakly nonlinear and weak noise regime, i.e., for $\mathcal{O}(\varepsilon^{2})$ nonlinearity and $\mathcal{O}(\varepsilon)$ noise, we prove uniform long-time pathwise error bounds in $H^1$ of size $\mathcal{O}(\varepsilon^2 \tau^\delta)$, for any $\delta<1$, up to times of length $\mathcal{O}(\varepsilon^{-2})$. Also here, a  reduction  by a factor of $\varepsilon^2$ appears. Due to the lack of exponential integrability of the proposed resonance-based low-regularity numerical solutions, that is needed to get strong convergence of numerical approximation to nonlinear SPDEs \cite{cui}, we obtain the pathwise error rates and pathwise long-term error bounds by combining low-regularity techniques with a truncation-based strategy \cite{deb0}. 
\item
In both the scenarios, we introduce a novel and nontrivial extension of the (RCO) technique to stochastic dispersive equations,  developing a refined analysis that allows to get improved long-time error bounds previously unavailable for such equations. The analytical techniques introduced in this work are expected to be extended beyond the present setting and to be applicable to other classes of stochastic dispersive equations.
\end{itemize}

%We conclude by highlighting that in both the aforementioned situations, an improvement of a factor of $\varepsilon^2$ wil be obtained, in terms of long-term error bounds, compared with results existing for classical exponential and splitting schemes.

The paper is organized as follows. In Section 2, we recall all the material needed for understanding this work. Section 3 will be fully devoted to the case of $F(\xi)=V(x)\xi$, while Section 4 will contain the theory for the cubic case $F(\xi)=|\xi|^2\xi.$ Finally, several numerical examples will be provided in Section 5, confirming the effectiveness of the theoretical analysis.

\section{Setting and notation}
In this section, we provide all the material and the notation useful for the sequel of the present work. 

Let $\mathbb{T}=[-\pi,\pi]$ and $T>0$. We consider the stochastic Schrödinger equation \eqref{SSE}, here re-written by the matter of clarity,
\begin{equation}
    \label{SSE_2}
d\xi(t)={\bf i}\Delta \xi(t)dt-{\bf i}F(\xi(t))dt+Q^{\frac{1}{2}}dW(t),\qquad  t\in[0,T],
\end{equation}
where $\xi(t)=\xi(t,x)$, $\xi(0)=\xi_0(x), \ x\in\mathbb{T}$, and ${\bf i}$ denotes the imaginary unit. Periodic boundary conditions will be assumed. Through this work, $\xi(t)\in L^2(\mathbb{T})$, being $L^2(\mathbb{T})$ the space of the periodic square-integrable complex-valued functions defined on $\mathbb{T}$, i.e., $\xi(t):\mathbb{T}\to \mathbb{C}$, for $t\in[0,T]$. 

it is well-known that an orthonormal basis of $L^2(\mathbb{T})$ is given by $\{\frac{1}{\sqrt{2\pi}}{\rm e}^{{\bf i}\ell x}\}_{\ell\in\mathbb{Z}}$; see, e.g., \cite{lord}. Moreover, we assume the following form for the $L^2(\mathbb{T})$-valued Wiener process $W(t)$
$$
W(t,x)=\frac{1}{\sqrt{2\pi}}\sum_{\ell\in\mathbb{Z}}\beta_\ell(t){\rm e}^{{\bf i}\ell x},
$$
where $\{\beta_\ell(t)\}_{\ell\in\mathbb{Z}}$ is a sequence of i.i.d. standard Brownian motions. The symmetric positive-definite Covariance operator $Q$ satisfies the following assumption.
\begin{assumption}
\label{ass_q}
    The operator $Q^{\frac{1}{2}}$ is an Hilbert-Schmidt operator on $L^2(\mathbb{T})$ and satisfies the following: there exists a real sequence $\{q_\ell\}_{\ell\in\mathbb{Z}}$ such that it holds
    $$
    Q {\rm e}^{{\rm e}\ell x}=q_\ell{\rm e}^{{\rm e}\ell x}, \quad q_{\ell}\ge 0, \qquad \ell\in\mathbb{Z}.
    $$
\end{assumption}
Based on Assumption \ref{ass_q}, one the gets
$$
Q^{\frac{1}{2}}W(t,x)=\frac{1}{\sqrt{2\pi}} \sum_{\ell\in\mathbb{Z}} \sqrt{q_\ell} \beta_{\ell}(t){\rm e}^{{\bf i}\ell x}.
$$
In this work, we are interested in the time approximation to mild solutions to \eqref{SSE_2}, i.e., the exact solution to \eqref{SSE_2} assumes the form
$$
\xi(t)={\rm e}^{{\bf i}t \Delta}\xi_0-{\bf i}\int_0^t {\rm e}^{{\bf i}(t-s)\Delta}F(\xi(s))\,ds+\int_0^t {\rm e}^{{\bf i}(t-s)\Delta}Q^{\frac{1}{2}}\,dW(s).
$$
For classical results on mild solutions to SPDEs, on the theory of the semigroup operator ${\rm e}^{{\bf i}t\Delta}, \ t\ge 0$, as well as, for standard martingale inequalities, such as, the well-known Burkholder-Davis-Gundy inequality (in the sequel, denoted as BDG inequality) and the stochastic Fubini Theorem, we refer to the monograph \cite{dap}.

On $L^2(\mathbb{T})$, we define the scalar product
$$
(\xi_1,\xi_2)=\Re\int_{\mathbb{T}} \xi_1(x)\overline{\xi_2}(x)\,dx, \qquad \xi_1,\xi_2\in\mathbb{L}^2(\mathbb{T}),
$$
where $\overline{\cdot}$ stands for the complex conjugate of $\cdot$. For $w\in L^2(\mathbb{T})$, we can write $w(x)=\frac{1}{\sqrt{2\pi}}\sum_{\ell\in\mathbb{Z}} \hat{w}_{\ell} {\rm e}^{{\bf i}\ell x}$, being $\{\hat{w}_{\ell}\}_{\ell\in\mathbb{Z}}$ the Fourier coefficients of $w$. Moreover, we will denote by $w^{(k)}$ the function $w^{(k)}(x)=\frac{1}{\sqrt{2\pi}}\sum_{\ell\in\mathbb{Z}} |\ell|^k|\hat{w}_{\ell} |{\rm e}^{{\bf i}\ell x}$, for any $k\in\mathbb{N}$. Also, for $q\in\mathbb{N}$, we will denote by
$$
\|w\|_{\ell^q}:=\big(\sum_{\ell\in\mathbb{Z}} |\hat{w}_{\ell}|^q\big)^{\frac{1}{q}}.
$$

and we denote by $\star$ the discrete convolution operator with 
$$
\widehat{\big(w\star z\big)}_{\ell} = \sum_{k\in\mathbb{Z}} \hat{w}_{k} \hat{z}_{\ell-k}.
$$
Clearly, for any $q\in\mathbb{N}$, one has
$$
\|w\star z\|_{\ell^q}=\big(\sum_{\ell\in\mathbb{Z}}\big|\sum_{k\in\mathbb{Z}} \hat{w}_{k} \hat{z}_{\ell-k}\big|^q\big)^{\frac{1}{q}}.
$$
For any $\sigma\in\mathbb{R}$, we can define the Sobolev norm of the space $W^{\sigma,2}(\mathbb{T})=H^\sigma(\mathbb{T})=H^\sigma$ as
$$
\|w\|_{H^\sigma}:=\big(\sum_{\ell\in\mathbb{Z}} (1+|\ell|^{2\sigma})|\hat{w}_{\ell}|^2\big)^{\frac{1}{2}}.
$$
Clearly, $L^2(\mathbb{T})=H^0$ and we will denote $\|\cdot \|_{L^2}=\|\cdot\|.$ 
Also, for $\sigma\ge 0$, we will denote by $\mathcal{L}_2^\sigma$ the space of Hilbert-Schmidt operators $\Psi: L^2(\mathbb{T}) \to H^\sigma$ with norm
$$
\|\Psi\|_{\mathcal{L}_2^\sigma}:=\big(\sum_{\ell\in\mathbb{Z}} \|\Psi e_{\ell}\|^2_{H^\sigma}\big)^{\frac{1}{2}},
$$
being $\{e_{\ell}\}_{\ell\in\mathbb{Z}}$ any orthonormal basis of $L^2(\mathbb{T})$. By Assumption \ref{ass_q}, $Q^{\frac{1}{2}}\in\mathcal{L}_2^0$. 

For a $H^\sigma$-valued stochastic process $z(t)$, we will use the notation
$$
\|z(t)\|_{L^{p}(\Omega; H^\sigma)}:=\big(\mathbb{E}\big[\|z(t)\|^{p}_{H^\sigma}\big]\big)^{\frac{1}{p}}, \quad p\in\mathbb{N}^+, \sigma\in\mathbb{R}
$$
and the notation
$$
\|z\|_{L^p(\Omega; L^{\infty}([0,T];H^\sigma))}:=\big(\mathbb{E}\big[\displaystyle\sup_{t\in[0,T]}\|z(t)\|_{H^\sigma}^p\big]\big)^{\frac{1}{p}}, \qquad p\in\mathbb{N}^+, \quad \sigma\in\mathbb{R}.
$$
By a matter of simplicity, we take the following assumption on the initial data $\xi_0$.
\begin{assumption}
\label{ass_data}
    The initial data $\xi_0(x)$ is a deterministic function.
\end{assumption}
\begin{remark}
    All the results presented in the sequel of this work may be extended to the case of stochastic $\xi_0(x)$, under the assumption that $\xi_0\in L^{p}(\Omega; H^\sigma)$, for some $\sigma\in\mathbb{R}$ and for a sufficiently large $p\in\mathbb{N}^+$.
\end{remark}
When not specified, in this work we will always work under Assumption \ref{ass_q} and Assumption \ref{ass_data}. Moreover, in the estimates, we will denote by $C$ (or $C_p$) a generic constant that may change line by line.

\section{Linear SSE with rough potential}
In this section, we study the case in which the mapping $F$ in \eqref{SSE} assumes the form $F(\xi)=V(x)\xi$, for a deterministic potential $V\in L^2(\mathbb{T})$, i.e., we consider the linear SSE with rough potential
\begin{equation}
    \label{lin_eq}
du(t)={\bf i}\Delta u(t) dt-{\bf i}V(x)u(t) dt+Q^{\frac{1}{2}}dW(t), \quad t\in[0,T], \ x \in \mathbb{T}, 
\end{equation}
whose mild solution is given by
\begin{equation}
    \label{mild_lin_eq}
u(t)={\rm e}^{{\bf i}t\Delta } u_0 -{\bf i}\int_{0}^{t}{\rm e}^{{\bf i}(t-s)\Delta}V(x)u(s)\,ds+\int_{0}^{t}{\rm e}^{{\bf i}(t-s)\Delta}Q^{\frac{1}{2}}\,dW(s). 
\end{equation}
Note that, through this section, we will denote the solution of the linear equation with rough potential $u$ to distinguish it from the cubic case that will be treated in the next section.

From \cite{an_co}, we have the following well-posedness result for \eqref{lin_eq}.
\begin{theorem}
\label{wp_lin}
    Let us assume $u_0 \in H^\sigma$ and $Q^{\frac{1}{2}}\in\mathcal{L}_2^\sigma$, for some $\sigma\in\mathbb{N}$. Assume also $V\in H^\sigma \cap L^\infty$.
    Then, for any $T>0$, there exists an unique solution $u$ to the linear SSE with rough potential \eqref{lin_eq} on $[0,T]$, given by \eqref{mild_lin_eq}.  Moreover, $u\in L^{2p}(\Omega; L^{\infty}([0,T];H^\sigma))$, for any $p\in\mathbb{N}^+$.
\end{theorem}

\begin{remark}
    Note that, since $\mathbb{T}$ is the one-dimensional torus, one has that, for $\sigma\in\mathbb{N}^+$, $H^\sigma \cap L^\infty=H^\sigma$ while, for $\sigma=0$, $L^2\cap L^\infty = L^\infty$, i.e., for the result in $L^2$, an uniform boundedness of the potential in $\mathbb{T}$ is needed.
\end{remark}
\subsection{Construction of a low-regularity numerical integrator}
Here, we construct a low-regularity integrator for \eqref{lin_eq}. Following the idea in \cite{os_ka}, it is convenient to define the twisted variable
\begin{equation}
    \label{twis_def}
    v(t)={\rm e}^{-{\bf i}t \Delta} u(t).
\end{equation}
Then, $v$ solves the equation
$$
v(t) =  v_0-{\bf i}\int_{0}^{t}{\rm e}^{-{\bf i}s\Delta}\big(V(x){\rm e}^{{\bf i} s \Delta}v(s)\big)\,ds+\int_{0}^{t}{\rm e}^{-{\bf i}s \Delta}Q^{\frac{1}{2}}\,dW(s).
$$
Then, after a change of variable (see also \cite{os_ka,fe_ka}), we get 
$$
v(t_n+\tau)=v(t_n)-{\bf i}\int_{0}^{\tau}{\rm e}^{-{\bf i}(t_n+s)\Delta}\big(V(x){\rm e}^{{\bf i}(t_n+s)\Delta}v(t_n+s)\big)\,ds+\int_{t_n}^{t_n+\tau}{\rm e}^{-{\bf i}s \Delta}Q^{\frac{1}{2}}\,dW(s).
$$
First, we consider the following approximation
\begin{equation}
    \label{app1}
v(t_n+\tau) \approx v(t_n)-{\bf i}\int_{0}^{\tau}{\rm e}^{-{\bf i}(t_n+s)\Delta}\big(V(x){\rm e}^{{\bf i}(t_n+s)\Delta}v(t_n)\big)\,ds+\int_{t_n}^{t_n+\tau}{\rm e}^{-{\bf i}s \Delta}Q^{\frac{1}{2}}\,dW(s).
\end{equation}
Hence, we now provide an approximation of the first integral in the right-hand side of \eqref{app1}. Under Fourier expansion, we have
\begin{align}
    \int_{0}^{\tau}{\rm e}^{-{\bf i}(t_n+s)\Delta}\big(V(x){\rm e}^{{\bf i}(t_n+s)\Delta}v(t_n)\big)\,ds
    &=\frac{1}{2\pi} \int_{0}^{\tau}{\rm e}^{-{\bf i}(t_n+s)\Delta}\big(\displaystyle\sum_{\ell\in \mathbb{Z}}\hat{V}_\ell {\rm e}^{{\bf i}\ell x}\cdot {\rm e}^{{\bf i}(t_n+s)\Delta}\displaystyle\sum_{\ell\in\mathbb{Z}}\hat{v}_\ell(t_n) {\rm e}^{{\bf i}\ell x}\big)\,ds\notag \\
    &=\frac{1}{2\pi} \int_{0}^{\tau}{\rm e}^{-{\bf i}(t_n+s)\Delta}\big(\displaystyle\sum_{\ell\in\mathbb{Z}}\hat{V}_\ell {\rm e}^{{\bf i}\ell x}\cdot \displaystyle\sum_{k\in \mathbb{Z}}\hat{v}_k (t_n) {\rm e}^{-{\bf i}(t_n+s)k^2}{\rm e}^{{\bf i}k x}\big)\,ds\notag \\
    &=\frac{1}{2\pi} \int_{0}^{\tau}{\rm e}^{-{\bf i}(t_n+s)\Delta}\displaystyle\sum_{\ell\in\mathbb{Z}}{\rm e}^{{\bf i}\ell x}\displaystyle\sum_{\substack{\ell_1,\ell_2\in\mathbb{Z}\\ \ell_1+\ell_2=\ell}}\hat{V}_{\ell_1} \hat{v}_{\ell_2}(t_n) {\rm e}^{-{\bf i}(t_n+s)\ell_2^2}\,ds\notag \\
    &=\frac{1}{2\pi} \int_{0}^{\tau}\displaystyle\sum_{\ell\in\mathbb{Z}}{\rm e}^{{\bf i}\ell x}\displaystyle\sum_{\substack{\ell_1,\ell_2\in\mathbb{Z}\\ \ell_1+\ell_2=\ell}}\hat{V}_{\ell_1} \hat{v}_{\ell_2}(t_n) {\rm e}^{{\bf i}(t_n+s)(\ell^2-\ell_2^2)}\,ds\notag \\
&=\frac{1}{2\pi} \displaystyle\sum_{\ell\in\mathbb{Z}}{\rm e}^{{\bf i}\ell x}\displaystyle\sum_{\substack{\ell_1,\ell_2\in\mathbb{Z}\\ \ell_1+\ell_2=\ell}}\hat{V}_{\ell_1} \hat{v}_{\ell_2}(t_n){\rm e}^{{\bf i}t_n(\ell^2-\ell_2^2)}\int_{0}^{\tau} {\rm e}^{{\bf i}s(\ell_1^2+2\ell_1 \ell_2)}\,ds \label{app2},
\end{align}
where we have used the key relation $\ell^2-\ell_2^2=\ell_1^2+2\ell_1 \ell_2$. Then, we make the following approximation
\begin{equation}
    \label{lr_app0}
\int_0^\tau {\rm e}^{{\bf i}s(\ell_1^2+2\ell_1 \ell_2)}\,ds\approx \int_0^\tau {\rm e}^{{\bf i}s \ell_1^2}\,ds = 
\begin{cases}
    \tau, & \ell_1=0,\\
    \frac{{\rm e}^{{\bf i}\ell_1^2\tau}-1}{{\bf i}\ell_1^2}, &\ell_1\in\mathbb{Z}_{\ne 0}.
\end{cases}
\end{equation}
Hence, plugging \eqref{lr_app0} in \eqref{app2}, we get
\begin{equation}
\label{lr_app1}
\int_{0}^{\tau}{\rm e}^{-{\bf i}(t_n+s)\Delta}\big(V(x){\rm e}^{{\bf i}(t_n+s)\Delta}v(t_n)\big)\,ds\approx \frac{\tau}{2\pi}
\displaystyle\sum_{\ell\in\mathbb{Z}}{\rm e}^{{\bf i}\ell x}\displaystyle\sum_{\substack{\ell_1,\ell_2\in\mathbb{Z}\\ \ell_1+\ell_2=\ell}}\hat{V}_{\ell_1} \hat{v}_{\ell_2}(t_n){\rm e}^{{\bf i}t_n(\ell^2-\ell_2^2)}\varphi_1({\bf i} \tau \ell_1^2),
%&=\tau {\rm e}^{-{\bf i}t_n\Delta }\big[({\rm e}^{{\bf i}t_n\Delta }v(t_n)) \varphi_1(-{\bf i}\tau \Delta)V(x)\big].
\end{equation}
where 
$$
\varphi_1(z)=\begin{cases}
    1, & z=0,\\
    \frac{1}{z}({\rm e}^{z}-1), & z\in\mathbb{Z}_{\ne 0}.
\end{cases}
$$
Returning back to physical space from Fourier expansion in \eqref{lr_app1}, we finally get
$$
\int_{0}^{\tau}{\rm e}^{-{\bf i}(t_n+s)\Delta}\big(V(x){\rm e}^{{\bf i}(t_n+s)\Delta}v(t_n)\big)\,ds\approx \tau {\rm e}^{-{\bf i}t_n\Delta }\big[({\rm e}^{{\bf i}t_n\Delta }v(t_n)) \varphi_1(-{\bf i}\tau \Delta)V(x)\big],
$$
where we define the operator $\varphi_1(\Delta)$ through its action on the elements of Fourier basis as
$$
\varphi_1(\Delta){\rm e}^{{\bf i}\ell x}=\begin{cases}
    1, & \ell=0,\\
    \frac{{\rm e}^{-\ell^2}-1}{-\ell^2}, &\ell\in\mathbb{Z}_{\ne 0}.
\end{cases}
$$
Then, we get the following approximation for the twisted variable $v$
$$
v_{n+1}=v_n-{\bf i}\tau {\rm e}^{-{\bf i}t_n\Delta }\big[({\rm e}^{{\bf i}t_n\Delta }v_n )\varphi_1(-{\bf i}\tau \Delta)V(x)\big]+\int_{t_n}^{t_n+\tau}{\rm e}^{-{\bf i}s \Delta}Q^{\frac{1}{2}}\,dW(s).
$$
Twisting back the variable $v$, we finally obtain the following numerical integrator for \eqref{lin_eq}
\begin{equation}
    \label{stand_lr_met}
u_{n+1}={\rm e}^{{\bf i}\tau\Delta}\big[u_n-{\bf i}\tau u_n \big(\varphi_1(-{\bf i}\tau \Delta)V(x)\big)\big]
+\mathcal{W}_n,
\end{equation}
where the stochastic convolution 
$$\mathcal{W}_n:=\int_{t_n}^{t_{n+1}}{\rm e}^{{\bf i}(t_{n+1}-s)\Delta} Q^{\frac{1}{2}}\,dW(s)
$$
is supposed to be simulated exactly \cite{lord}. In the sequel, we may refer to the integrator \eqref{stand_lr_met} as (SLR1).

\subsection{Error analysis for (SLR1)}
\label{sec_lr}
Here, we provide an error analysis for (SLR1). We start by providing the following proposition.
\begin{proposition}
\label{prop_stab_lin}
    Let us consider $u_0, V\in H^{\sigma}$ and $Q^{\frac{1}{2}}\in\mathcal{L}_2^\sigma$, for some $\sigma\in\mathbb{N}^+$. Then, for any $p\in\mathbb{N}^+$, there exists a positive constant $C_p$ such that the numerical solution given by (SLR1) satisfies
    \begin{equation*}
    \sup_{N\in\mathbb{N}^+}\mathbb{E}\big[\max_{1\le n\le N}\|u_n\|^{2p}_{H^\sigma}\big]\le C_p.
    \end{equation*}
\end{proposition}
\begin{proof}
    By recursion of the iteration in \eqref{stand_lr_met}, we get
    $$
    u_n={\rm e}^{{\bf i}t_n\Delta}u_0-{\bf i}\tau \displaystyle\sum_{j=1}^{n}{\rm e}^{{\bf i}j\tau\Delta}\big(u_{n-j}\varphi_1(-{\bf i}\tau\Delta)V(x)\big)
    +\displaystyle\sum_{j=1}^{n}{\rm e}^{{\bf i}(j-1)\tau\Delta}\mathcal{W}_{n-j}.
    $$
    Then, using the fact that ${\rm e}^{{\bf i}t\Delta}$ is an isometry, that $H^\sigma$ is an algebra on $\mathbb{T}$, for $\sigma\in\mathbb{N}^+$, that $\varphi_1(-{\bf i}\tau\Delta)$ is bounded on $H^\sigma$, and the Hölder inequality, we get
    \begin{align*}
    \|u_n\|_{H^\sigma}^{2p}&\le C \big(\|u_0\|_{H^\sigma}^{2p}+\tau^{2p}\|\displaystyle\sum_{j=1}^{n}{\rm e}^{{\bf i}j\tau\Delta}\big(u_{n-j}\varphi_1(-{\bf i}\tau\Delta)V(x)\big)\|_{H^\sigma}^{2p}+\|\displaystyle\sum_{j=1}^{n}{\rm e}^{{\bf i}(j-1)\tau\Delta}\mathcal{W}_{n-j}\|_{H^\sigma}^{2p}\big)\\
    &\le C\big(\|u_0\|_{H^\sigma}^{2p}+\tau \displaystyle\sum_{j=1}^{n}\|u_{n-j}\varphi_1(-{\bf i}\tau\Delta)V(x)\|^{2p}_{H^\sigma}+\|\displaystyle\sum_{j=0}^{n-1}\int_{t_j}^{t_{j+1}}{\rm e}^{{\bf i}(t_n-s)\Delta}Q^{\frac{1}{2}}\,dW(s)\|^{2p}_{H^\sigma}\big)\\
    &\le C\big(\|u_0\|_{H^\sigma}^{2p}+\tau\|V\|^{2p}_{H^\sigma} \displaystyle\sum_{j=0}^{n-1}\|u_{j}\|^{2p}_{H^\sigma}+\|\displaystyle\sum_{j=0}^{n-1}\int_{t_j}^{t_{j+1}}{\rm e}^{{\bf i}(t_n-s)\Delta}Q^{\frac{1}{2}}\,dW(s)\|^{2p}_{H^\sigma}\big).
    \end{align*}
    Then, by BDG inequality, we get
    \begin{align}
    \mathbb{E}\big[\max_{1\le n\le N}\|u_n\|^{2p}_{H^\sigma}\big]&\le C \|u_0\|_{H^\sigma}^{2p}+C\tau \|V\|^{2p}_{H^\sigma} \displaystyle\sum_{j=0}^{N-1}\mathbb{E}\big[\max_{1\le k\le j}\|u_k\|_{H^\sigma}^{2p}\big]\notag \\
    &\qquad \qquad +C\mathbb{E}\big[\displaystyle\sup_{t\in[0,T]}\|\int_0^t {\rm e}^{{\bf i}(t-s)\Delta}Q^{\frac{1}{2}}\,dW(s)\|^{2p}_{H^\sigma}\big]\notag \\
    &\le  C \|u_0\|_{H^\sigma}^{2p}+C\tau \|V\|^{2p}_{H^\sigma} \displaystyle\sum_{j=0}^{N-1}\mathbb{E}\big[\max_{1\le k\le j}\|u_k\|_{H^\sigma}^{2p}\big]+CT^p \|Q^{\frac{1}{2}}\|^{2p}_{\mathcal{L}_2^\sigma}.
    \end{align}
    An application of the discrete Gronwall inequality then yields the result.
\end{proof}
\begin{remark}
\label{stab_lim_L2}
    It is possible to obtain the same result of Proposition \ref{prop_stab_lin} for $\sigma=0$, i.e., for $L^2(\mathbb{T})$, if one assumes $V\in H^{\frac{1}{2}^+}$. Indeed, for any $s>\frac{1}{2}$, Hölder inequality gives
    \begin{align*}
    \|\varphi_1(-{\bf i}\tau\Delta) V(x)\|_{L^\infty} &= \frac{1}{\sqrt{2\pi}}\displaystyle\sup_{x\in\mathbb{T}}|\hat{V}_0+\displaystyle\sum_{\ell\in\mathbb{Z}_{\ne 0}}{\rm e}^{{\bf i}\ell x} \hat{V}_\ell\frac{1}{{\bf i}\tau\ell^2}({\rm e}^{{\bf i}\tau \ell^2}-1)|\\
    &\le \frac{1}{\sqrt{2\pi}}\displaystyle\sum_{\ell\in\mathbb{Z}}|\hat{V}_{\ell}|(1+|\ell|)^{s}|(1+|\ell|)^{-s}\\
    &\le \frac{1}{\sqrt{2\pi}} \big(\displaystyle\sum_{\ell\in\mathbb{Z}}|\hat{V}_{\ell}|^2(1+|\ell|)^{2s}\big)^{\frac{1}{2}}\big(\displaystyle\sum_{\ell\in\mathbb{Z}}(1+|\ell|)^{-2s}\big)^{\frac{1}{2}}\\
    &\le C \|V\|_{H^s}.
    \end{align*}
    This gives the bound
    $$
    \|u_j \varphi_1(-{\bf i}\tau\Delta) V(x)\|\le C \|u_j\| \|V\|_{H^s},\qquad j=0,\dots,N, \quad s>\frac{1}{2}.
    $$
\end{remark}
We now state and prove strong and $\mathbb{P}-{\rm a.s.}$ convergence results for the integrator (SLR1).
\begin{theorem}
\label{conv_res_lin}
Let $\sigma \in \mathbb{N}^+$ and $\gamma\in(0,1]$. Assume $u_0, V\in H^{\sigma+\gamma}$ and $Q^{\frac{1}{2}}\in\mathcal{L}_2^{\sigma+\gamma}$.  
%assume $f\in H^{\gamma+\frac{1}{2}+s}$ if $\sigma=0$, with any small $s>0$, or
Then, for any $p\in\mathbb{N}^+$, the exists a constant $C_p$ such that we have the following error estimate
\begin{equation}
\label{strong_conv_lin}
\sup_{N\in\mathbb{N}^+}\big(\mathbb{E}\big[\max_{1\le n\le N}\|u(t_n)-u_n\|^{2p}_{H^\sigma}  \big]\big)^{\frac{1}{2p}} \le C_p\tau^\gamma.
\end{equation}
Moreover, it holds that for any $\delta<\gamma$ there exists a positive random variable $K_\delta(\omega)$ such that one has
    \begin{equation}
        \label{path_conv_lin}
        \max_{1\le n\le N}\|u(t_n)-u_n\|_{H^\sigma}\le K_\delta(\omega)\tau^\delta, \qquad \mathbb{P}-{\rm a.s.}
    \end{equation}
\end{theorem}
\begin{proof}
Since ${\rm e}^{{\bf i}t\Delta}$ is an isometry on $H^\sigma$, we prove the result for the twisted variable $v$. We have
$$
v(t_n)-v_n=-{\bf i}\displaystyle\sum_{j=0}^{n-1}\int_{0}^{\tau}{\rm e}^{-{\bf i}(t_j+s)\Delta} V(x){\rm e}^{{\bf i}(t_j+s)\Delta}v(t_j+s)\,ds
-\tau {\rm e}^{-{\bf i}t_j\Delta}\big[\big({\rm e}^{{\bf i}t_j \Delta} v_j\big)\varphi_1(-{\bf i}\tau\Delta)V(x)\big].
$$
Then, we see
$$
\|v(t_n)-v_n\|^{2p}_{H^\sigma}\le C(I_1+\tau^{2p} I_2),
$$ where
\begin{align*}
    &I_1:=\big\|\displaystyle\sum_{j=0}^{n-1}\int_{0}^{\tau}{\rm e}^{-{\bf i}(t_j+s)\Delta} V(x){\rm e}^{{\bf i}(t_j+s)\Delta}v(t_j+s)\,ds-\tau {\rm e}^{-{\bf i}t_j\Delta}\big[\big({\rm e}^{{\bf i}t_j \Delta} v(t_j)\big)\varphi_1(-{\bf i}\tau\Delta)V(x)\big]\big\|^{2p}_{H^\sigma},\\
    &I_2:=\big\|\displaystyle\sum_{j=0}^{n-1}{\rm e}^{-{\bf i}t_j\Delta}\big[\big({\rm e}^{{\bf i}t_j \Delta} v(t_j)\big)\varphi_1(-{\bf i}\tau\Delta)V(x)\big]-{\rm e}^{-{\bf i}t_j\Delta}\big[\big({\rm e}^{{\bf i}t_j \Delta} v_j\big)\varphi_1(-{\bf i}\tau\Delta)V(x)\big]\big\|^{2p}_{H^\sigma}.
\end{align*}
The term $I_2$ is estimated using the Hölder inequality, the algebra property of $H^\sigma$, the unitary property of ${\rm e}^{{\bf i}t\Delta}$ and the boundedness property of $\varphi_1$ in $H^\sigma$. One indeed has
$$
\begin{aligned}
    I_2&\le n^{2p-1} \displaystyle\sum_{j=0}^{n-1}\|{\rm e}^{{\bf i}t_j \Delta}(v(t_j)-v_j)\varphi_1(-{\bf i}\tau \Delta)V(x)\|^{2p}_{H^\sigma}\\
    &\le C \|V\|^{2p}_{H^\sigma} n^{2p-1} \displaystyle\sum_{j=0}^{n-1}\|v(t_j)-v_j\|^{2p}_{H^\sigma}.
\end{aligned}
$$
For the term $I_1$, we note that $I_1\le C( I_{1,1}+I_{1,2})$, where
$$
\begin{aligned}
  &  I_{1,1}:=\big\|\displaystyle\sum_{j=0}^{n-1}\int_{0}^{\tau}{\rm e}^{-{\bf i}(t_j+s)\Delta} V(x){\rm e}^{{\bf i}(t_j+s)\Delta}\big(v(t_j+s)-v(t_j)\big)\,ds\big\|^{2p}_{H^\sigma},\quad 
    I_{1,2}:=\big\|\displaystyle\sum_{j=0}^{n-1}\mathcal{R}_{t_j,\tau}(v(t_j))\big\|^{2p}_{H^\sigma},
\end{aligned}
$$
with
$$
\mathcal{R}_{t_j,\tau}(v(t_j)):=\displaystyle\sum_{\ell\in\mathbb{Z}}{\rm e}^{{\bf i}\ell x}\displaystyle\sum_{\substack{\ell_1,\ell_2\in\mathbb{Z}\\ \ell_1+\ell_2=\ell}}\hat{V}_{\ell_1}\hat{v}_{\ell_2}(t_j){\rm e}^{{\bf i}t_n(\ell^2-\ell_2^2)}\int_{0}^{\tau}{\rm e}^{{\bf i}s \ell_1^2}\big({\rm e}^{2{\bf i}s \ell_1 \ell_2}-1\big)\,ds.
$$
Using Fourier decomposition and an interpolation argument, we have
\begin{align*}
    \|\mathcal{R}_{t_j,\tau}(v(t_j))\|^2_{H^\sigma}&=\displaystyle\sum_{\ell\in\mathbb{Z}}(1+|\ell|)^{2\sigma}\big|\displaystyle\sum_{\substack{\ell_1,\ell_2\in\mathbb{Z}\\ \ell_1+\ell_2=\ell}}\hat{V}_{\ell_1}\hat{v}_{\ell_2}(t_j){\rm e}^{{\bf i}t_n(\ell^2-\ell_2^2)}\int_{0}^{\tau}{\rm e}^{{\bf i}s \ell_1^2}\big({\rm e}^{2{\bf i}s \ell_1 \ell_2}-1\big)\,ds\big|^2\\
    &\le \displaystyle\sum_{\ell\in\mathbb{Z}}(1+|\ell|)^{2\sigma}\big(\displaystyle\sum_{\substack{\ell_1,\ell_2\in\mathbb{Z}\\ \ell_1+\ell_2=\ell}}|\hat{V}_{\ell_1}| |\hat{v}_{\ell_2}(t_j)|\int_{0}^{\tau}|2s\ell_1 \ell_2|^\gamma\,ds\big)^2\\
    &\le C \tau^{2(\gamma+1)} \displaystyle\sum_{\ell\in\mathbb{Z}}(1+|\ell|)^{2\sigma}\big(\sum_{\substack{\ell_1,\ell_2\in\mathbb{Z}\\ \ell_1+\ell_2=\ell}}|\hat{V}_{\ell_1}| |\hat{v}_{\ell_2}(t_j)|(1+|\ell_1|)^\gamma(1+|\ell_2|)^\gamma\big)^2.
\end{align*}
This gives
$$
 \|\mathcal{R}_{t_j,\tau}(v(t_j))\|^2_{H^\sigma}\le C\tau^{2(\gamma+1)}\|V\|^2_{H^{\sigma+\gamma}}\|v(t_j)\|^2_{H^{\sigma+\gamma}}.
$$
This, combined with the Hölder inequality, then gives
$$
I_{1,2}\le C \tau^{2p\gamma+1}\|V\|^{2p}_{H^{\sigma+\gamma}}\displaystyle\sum_{j=0}^{n-1}\|v(t_j)\|^{2p}_{H^{\sigma+\gamma}}.
%&\le C \|f\|_{H^{\sigma+\gamma}} \|v\|_{C([0,T];L^{2p}(\Omega; H^{\sigma+\gamma}))}\tau^\gamma.
$$
It remains to estimate $I_{1,1}$. We use the mild formulation for $v(t_j+s)$, so we can decompose this term as $I_{1,1}\le C( I_{1,1,1}+I_{1,1,2})$, where
\begin{align*}
  &  I_{1,1,1}:=\big\|\displaystyle\sum_{j=0}^{n-1}\int_{0}^{\tau}{\rm e}^{-{\bf i}(t_j+s)\Delta} V(x){\rm e}^{{\bf i}(t_j+s)\Delta}\bigg(\int_{0}^{s}{\rm e}^{-{\bf i}(t_j+r)\Delta}V(x){\rm e}^{{\bf i}(t_j+r)\Delta}v(t_j+r)\,dr\bigg)\,ds\big\|^{2p}_{H^\sigma},\\
    &I_{1,1,2}:=\big\|\displaystyle\sum_{j=0}^{n-1}\int_{0}^{\tau}{\rm e}^{-{\bf i}(t_j+s)\Delta} V(x){\rm e}^{{\bf i}(t_j+s)\Delta}\big(\int_{t_j}^{t_j+s}{\rm e}^{-{\bf i}r \Delta}Q^{\frac{1}{2}}\,dW(r)\big)\,ds\big\|^{2p}_{H^\sigma}.
\end{align*}
For $I_{1,1,1}$, using similar techniques as above, we have
\begin{align*}
I_{1,1,1}&\le n^{2p-1}\displaystyle\sum_{j=0}^{n-1}\big(\int_0^\tau \|V(x){\rm e}^{{\bf i}(t_j+s)\Delta}\big(\int_{0}^{s}{\rm e}^{-{\bf i}(t_j+r)\Delta}V(x){\rm e}^{{\bf i}(t_j+r)\Delta}v(t_j+r)\,dr\big)\|_{H^\sigma}\,ds \big)^{2p}\\
&\le \displaystyle\sum_{j=0}^{n-1} \int_0^\tau \|V(x){\rm e}^{{\bf i}(t_j+s)\Delta}\big(\int_{0}^{s}{\rm e}^{-{\bf i}(t_j+r)\Delta}V(x){\rm e}^{{\bf i}(t_j+r)\Delta}v(t_j+r)\,dr\big)\|^{2p}_{H^\sigma}\,ds\\
&\le \tau^{2p-1}\|V\|^{2p}_{H^\sigma}\displaystyle\sum_{j=0}^{n-1}\int_0^\tau \int_0^s \|V(x)\|^{2p}_{H^\sigma} \|v(t_j+r)\|_{H^\sigma}^{2p}\,dr\,ds\\
&\le \tau^{2p+1}\|V\|^{4p}_{H^\sigma}\displaystyle\sum_{j=0}^{n-1}\displaystyle\sup_{t\in [t_j,t_j+\tau]}\|v(t)\|^{2p}_{H^\sigma}.
\end{align*}
Then, we achieve the following bound
\begin{align*}
&\mathbb{E}\big[\max_{1\le n\le N}\|v(t_n)-v_n\|^{2p}_{H^\sigma}\big] \le C\tau^{2p+1}\|V\|^{4p}_{H^\sigma}\mathbb{E}\big[\max_{1\le n \le N} \displaystyle\sum_{j=0}^{n-1}\displaystyle\sup_{t\in[t_j,t_j+\tau]}\|v(t_j)\|^{2p}_{H^\sigma}\big]\\
&\hspace{2cm}+C \tau^{2p\gamma+1}\|V\|^{2p}_{H^{\sigma+\gamma}}\mathbb{E}\big[\max_{1\le n\le N}\displaystyle\sum_{j=0}^{n-1}\|v(t_j)\|^{2p}_{H^{\sigma+\gamma}}\big]\\
& \hspace{2cm}+C\tau \|V\|^{2p}_{H^\sigma}\mathbb{E}\big[\max_{1\le n\le N}\displaystyle\sum_{j=0}^{n-1} \|v(t_j)-v_j\|^{2p}_{H^\sigma}\big]\\
&\hspace{2cm} +C\mathbb{E}\big[\max_{1\le n \le N}\|\displaystyle\sum_{j=0}^{n-1}\int_{0}^{\tau}{\rm e}^{-{\bf i}(t_j+s)\Delta} V(x){\rm e}^{{\bf i}(t_j+s)\Delta}\big(\int_{t_j}^{t_j+s}{\rm e}^{-{\bf i}r \Delta}Q^{\frac{1}{2}}\,dW(r)\big)\,ds\|^{2p}_{H^\sigma}\big]\\
& \hspace{0.5cm}\le C\tau^{2p}\|V\|^{4p}_{H^\sigma}\mathbb{E}\big[\displaystyle\sup_{t\in[0,T]}\|v(t)\|^{2p}_{H^\sigma}\big]+C \tau^{2p\gamma}\|V\|^{2p}_{H^{\sigma+\gamma}}\displaystyle\sup_{t\in[0,T]}\mathbb{E}\big[\|v(t)\|^{2p}_{H^{\sigma+\gamma}}\big]\\
&\hspace{2cm}+C \tau \|V\|^{2p}_{H^\sigma} \displaystyle\sum_{j=0}^{N-1}\mathbb{E}\big[\max_{0\le k\le j}\|v(t_k)-v_k\|^{2p}_{H^\sigma}\big]\\
&\hspace{2cm}+C\mathbb{E}\big[\max_{1\le n \le N}\|\displaystyle\sum_{j=0}^{n-1}\int_{0}^{\tau}{\rm e}^{-{\bf i}(t_j+s)\Delta} V(x){\rm e}^{{\bf i}(t_j+s)\Delta}\big(\int_{t_j}^{t_j+s}{\rm e}^{-{\bf i}r \Delta}Q^{\frac{1}{2}}\,dW(r)\big)\,ds\|^{2p}_{H^\sigma}\big].
\end{align*}
Finally, by stochastic Fubini Theorem (see also \cite{an_co}) and BDG inequality, we obtain the following estimate
\begin{align*}
&\mathbb{E}\big[\max_{1\le n\le N}\big\|\displaystyle\sum_{j=0}^{n-1}\int_{0}^{\tau}{\rm e}^{-{\bf i}(t_j+s)\Delta} V(x){\rm e}^{{\bf i}(t_j+s)\Delta}\bigg(\int_{t_j}^{t_j+s}{\rm e}^{-{\bf i}r \Delta}Q^{\frac{1}{2}}\,dW(r)\bigg)\,ds\big\|^{2p}_{H^\sigma}\big]\\
    &\hspace{2cm}=\mathbb{E}\big[\max_{1\le n\le N}\big\|\displaystyle\sum_{j=0}^{n-1}\int_{t_j}^{t_{j+1}}\int_{r}^{t_{j+1}}{\rm e}^{-{\bf i}s\Delta} V(x){\rm e}^{{\bf i}s\Delta}{\rm e}^{-{\bf i}r \Delta}Q^{\frac{1}{2}}\,ds\,dW(r)\big\|^{2p}_{H^\sigma}\big]\\
     &\hspace{2cm}=\mathbb{E}\big[\max_{1\le n\le N}\big\|\int_{0}^{t_{n}}\int_{r}^{[\frac{r}{\tau}+1]\tau}{\rm e}^{-{\bf i}s\Delta} V(x){\rm e}^{{\bf i}s\Delta}{\rm e}^{-{\bf i}r \Delta}Q^{\frac{1}{2}}\,ds\,dW(r)\big\|^{2p}_{H^\sigma}\big]\\
     &\hspace{2cm}\le \mathbb{E}\big[\displaystyle\sup_{t\in[0,T]}\big\|\int_{0}^{t}\int_{r}^{[\frac{r}{\tau}+1]\tau}{\rm e}^{-{\bf i}s\Delta} V(x){\rm e}^{{\bf i}s\Delta}{\rm e}^{-{\bf i}r \Delta}Q^{\frac{1}{2}}\,ds\,dW(r)\big\|^{2p}_{H^\sigma}\big]\\
     &\hspace{2cm} \le C \mathbb{E}\big[\big(\int_{0}^{T}\big\|\int_{r}^{[\frac{r}{\tau}+1]\tau}{\rm e}^{-{\bf i}s\Delta} V(x){\rm e}^{{\bf i}s\Delta}{\rm e}^{-{\bf i}r \Delta}Q^{\frac{1}{2}}\,ds\big\|^2_{\mathcal{L}_2^\sigma}\,dr\big)^p\big]\\
     &\hspace{2cm}= C \mathbb{E}\big[\big(\displaystyle\sum_{j=0}^{N-1}\int_{t_j}^{t_{j+1}}\big\|\int_{r}^{t_{j+1}}{\rm e}^{-{\bf i}s\Delta} V(x){\rm e}^{{\bf i}s\Delta}{\rm e}^{-{\bf i}r \Delta}Q^{\frac{1}{2}}\,ds\big\|^2_{\mathcal{L}_2^\sigma}\,dr\big)^p\big]\\
     &\hspace{2cm}\le C\mathbb{E}\big[\big(\displaystyle\sum_{j=0}^{N-1}\int_{t_j}^{t_{j+1}}\big(\int_{r}^{t_{j+1}}\big\|{\rm e}^{-{\bf i}s\Delta} V(x){\rm e}^{{\bf i}s\Delta}{\rm e}^{-{\bf i}r \Delta}Q^{\frac{1}{2}}\big\|_{\mathcal{L}_2^\sigma}\,ds\big)^2\,dr\big)^p\big]\\
     &\hspace{2cm}\le C\|V\|^{2p}_{H^\sigma}\big\|Q^{\frac{1}{2}}\big\|^{2p}_{\mathcal{L}_2^\sigma}\mathbb{E}\big[\big(\displaystyle\sum_{j=0}^{N-1}\int_{t_{j}}^{t_{j+1}}|t_{j+1}-r|^2\,dr\big)^p\big]\\
     &\hspace{2cm}\le  C\|V\|^{2p}_{H^\sigma}\big\|Q^{\frac{1}{2}}\big\|^{2p}_{\mathcal{L}_2^\sigma}\tau^{2p}.
\end{align*}
Then, gathering all the above estimates, thanks to Proposition \ref{wp_lin} and taking into account that $2p\gamma\le 2p$ for any $\gamma\le 1$, we get
$$
\mathbb{E}\big[\max_{1\le n\le N}\|v(t_n)-v_n\|^{2p}_{H^\sigma}\big] \le C\tau^{2p\gamma}+C \tau \|V\|^{2p}_{H^\sigma} \displaystyle\sum_{j=0}^{N-1}\mathbb{E}\big[\max_{0\le k\le j}\|v(t_k)-v_k\|^{2p}_{H^\sigma}\big].
$$
An application of the discrete Gronwall inequality yields the estimate in \eqref{strong_conv_lin}. Finally, the estimate in \eqref{path_conv_lin} directly follows by an application of Markov inequality, the estimate in \eqref{strong_conv_lin} and Borel-Cantelli Lemma. The proof is completed.
\end{proof}

\begin{remark}
    By Theorem \ref{conv_res_lin}, strong order one of convergence in $H^\sigma, \sigma\in\mathbb{N}^+,$ for (SLR1) is achieved for $u_0, V \in H^{\sigma+1}$ and $Q^{\frac{1}{2}}\in\mathcal{L}_2^{\sigma+1}$. This is a regularity improvement if compared with classical theory \cite{an_co} since it would require $u_0, V \in H^{\sigma+2}$ and $Q^{\frac{1}{2}}\in \mathcal{L}_2^{\sigma+2}$.
\end{remark}
\begin{remark} It is possible to obtain strong order $\gamma$ in $L^2$, i.e., for $\sigma=0$, by imposing $u_0\in H^\gamma, Q^{\frac{1}{2}}\in \mathcal{L}_2^\gamma$ and $V\in H^{\gamma+\frac{1}{2}+s}$, for any small $s>0$. The proof relies on similar computations than those shown in remark \ref{stab_lim_L2}. Here, we provide only the computation for the term $\|\mathcal{R}_{t_j,\tau}(v(t_j))\|$.
%It now remains to prove the result for the case $\sigma = 0$. For this case, the term $I_2$ is bounded as
%$$
%\begin{aligned}
 %   I_2&\le \displaystyle\sum_{j=0}^{n-1}\left\|{\rm e}^{-{\bf i}t_j \Delta}(v(t_j)-v_j)\varphi_1({\bf i}\tau \Delta)f(x)\right\|_{L^{2p}(\Omega;L^2)}\\
 %   &\le C \| \varphi_1({\bf i}\tau \Delta) f\|_{L^\infty} \displaystyle\sum_{j=0}^{n-1}\|v(t_j)-v_j\|_{L^{2p}(\Omega; L^2)}\\
  %  &\le C \| \varphi_1({\bf i}\tau \Delta) f\|_{H^{\frac{1}{2}+s}} \displaystyle\sum_{j=0}^{n-1}\|v(t_j)-v_j\|_{L^{2p}(\Omega; L^2)}\\
 %   &\le C \|  f\|_{H^{\frac{1}{2}+s}} \displaystyle\sum_{j=0}^{n-1}\|v(t_j)-v_j\|_{L^{2p}(\Omega; L^2)}.
%\end{aligned}
%$$
We have
\begin{align*}
    \|\mathcal{R}_{t_j,\tau}(v(t_j))\|^2&=\displaystyle\sum_{\ell\in\mathbb{Z}}\big|\displaystyle\sum_{\substack{\ell_1,\ell_2\in\mathbb{Z}\\ \ell_1+\ell_2=\ell}}\hat{V}_{\ell_1}\hat{v}_{\ell_2}(t_j){\rm e}^{{\bf i}t_n(\ell^2-\ell_2^2)}\int_{0}^{\tau}{\rm e}^{{\bf i}s \ell_1^2}\big({\rm e}^{2{\bf i}s \ell_1 \ell_2}-1\big)\,ds\big|^2\\
    &\le \displaystyle\sum_{\ell\in\mathbb{Z}}\big(\displaystyle\sum_{\substack{\ell_1,\ell_2\in\mathbb{Z}\\ \ell_1+\ell_2=\ell}}|\hat{V}_{\ell_1}| |\hat{v}_{\ell_2}(t_j)|\int_{0}^{\tau}|2s\ell_1 \ell_2|^\gamma\,ds\big)^2\\
    &\le C \tau^{2(\gamma+1)} \displaystyle\sum_{k\in\mathbb{Z}}\big(\sum_{\ell\in\mathbb{Z}}|\hat{V}_{k}| |\hat{v}_{k-\ell}(t_j)||k|^\gamma |k-\ell|^\gamma\big)^2.
\end{align*}
Then, the Young's inequality and the Hölder inequality then give
$$
\begin{aligned}
\|\mathcal{R}_{t_j,\tau}(v(t_j))\|&\le C \tau^{\gamma+1}\|V^{(\gamma)}(x)\star v^{(\gamma)}(t_j)\|_{\ell^2}\\
&\le C \tau^{\gamma+1}\|V^{(\gamma)}\|_{\ell^1}\|v^{(\gamma)}(t_j)\|_{\ell^2}\\
&\le C \tau^{\gamma+1}\|V\|_{H^{\gamma+\frac{1}{2}+s}}\|v(t_j)\|_{H^\gamma}.
\end{aligned}
$$
The rest of the proof is similar to the proof of theorem \ref{conv_res_lin}.

Finally, note that the strong order one in $L^2$ for (SLR1) is achieved with $u_0 \in H^1, Q^{\frac{1}{2}}\in\mathcal{L}_2^1$ and $V\in H^{\frac{3}{2}^+}$, i.e., less regularity is requested if compared with classical theory \cite{an_co} that would require $u_0, V\in H^2$ and $Q^{\frac{1}{2}}\in\mathcal{L}_2^2$.
\end{remark}

\begin{remark}
\label{app_W}
    If one approximates the stochastic convolution $\mathcal{W}_n$ with, e.g., exponential Euler method, i.e.,
    $$
    \mathcal{W}_n\approx {\rm e}^{{\bf i}\tau \Delta}Q^{\frac{1}{2}} \big(W(t_{n+1})-W(t_n)\big)
    $$
    then order $\gamma$ of convergence is maintained for $Q^{\frac{1}{2}}\in \mathcal{L}_2^{\sigma+2\gamma}$; see, e.g., \cite{an_co}. Thus, more regularity is required to the Covariance operator $Q$.
\end{remark}
As shown in \cite{an_co} for the standard stochastic exponential method, applied to \eqref{lin_eq}, also for (SLR1) we can show a nearly conservation of the trace equations for the expected mass and the energy computed along the solution to \eqref{lin_eq}. Indeed, we define the quantities
$$
M(u)=\int_{\mathbb{T}} |u|^2\,dx, \qquad H(u)=\frac{1}{2}\int_{\mathbb{T}}|\nabla u|^2\,dx-\frac{1}{2}\int_{\mathbb{T}} V(x)|u|^2\,dx.
$$
and provide the following results, whose proofs is omitted since it is similar to the proofs on Proposition 4.6 and Proposition 4.8 in \cite{an_co}.
\begin{proposition}
    Leu us assume $u_0\in L^2, V\in H^{\frac{1}{2}^+}$ and $Q^{\frac{1}{2}}\in\mathcal{L}_2^0$. Then, the numerical method given by (SLR1) satisfies
    $$
    \big|\mathbb{E}\big[\|u(t_n)\|^2\big]-\mathbb{E}\big[\|u_n\|^2\big]\big| \le C \tau, \qquad 0\le t_n=n\tau\le T.
    $$
    Moreover, if $u_0 \in H^1$, $V\in H^{\frac{3}{2}^+}$ and $Q^{\frac{1}{2}}\in\mathcal{L}_2^{1}$, then, the numerical integrator (SLR1) satisfies
    $$
    \big|\mathbb{E}\big[H(u(t_n))\big]-\mathbb{E}\big[H(u_n)\big]\big| \le C \tau, \qquad 0\le t_n=n\tau\le T.
    $$
\end{proposition}

\subsection{Long-time error estimates for small potential and noise}
Let $0<\varepsilon\ll1$. We consider the linear stochastic Schrodinger equation perturbed with a small rough  potential  and small noise, i.e., for a real valued function $\Phi$, we consider $V(x)=\varepsilon^2 \Phi(x)$ and the equation
\begin{equation}
    \label{lin_small}
d u^\varepsilon(t)={\bf i}\Delta u^\varepsilon(t) dt-{\bf i}\varepsilon^2 \Phi(x) u^\varepsilon(t)dt+\varepsilon Q^{\frac{1}{2}} dW.
\end{equation}
Equation \eqref{lin_small} describes small random fluctuations around the dynamics given by the classical deterministic linear Schrödinger equations.
The mild solution to Equation \eqref{lin_small} reads
$$
u^\varepsilon(t)={\rm e}^{{\bf i}t \Delta}u_0-{\bf i}\varepsilon^2 \int_0^t {\rm e}^{{\bf i}(t-s)\Delta}\Phi(x)u^\varepsilon(s)\,ds+\varepsilon\int_0^t {\rm e}^{{\bf i}(t-s)\Delta}Q^{\frac{1}{2}}\,dW(s).
$$
In the sequel of this work, we will use the following notation: $A\lesssim B$ denotes $A\le c B$, where $c$ is a positive deterministic constant independent on $\varepsilon$  and $\tau$.

Let us define $T_\varepsilon:=\frac{T}{\varepsilon^2}.$
We have the following result. 
\begin{proposition}
\label{smal_prop}
    Let $\sigma\in\mathbb{N}^+$ and let $u_0, \Phi \in H^\sigma$ and $Q^{\frac{1}{2}}\in\mathcal{L}_2^\sigma$. Then, the exact solution to \eqref{lin_small} satisfies
    $$
   \mathbb{E}\big[\sup_{t\in[0,T_\varepsilon]} \|u^\varepsilon(t)\|_{H^\sigma}^2\big]\lesssim 1.
    $$
    \end{proposition}

    \begin{proof}
By proposition \ref{wp_lin}, we have global well-posedness which implies the existence and uniqueness of solution to time $\frac{T}{\varepsilon^2}$. We only need to show the independence of the bound from $\varepsilon$ in the time window $[0,\frac{T}{\varepsilon^2}]$. Direct computation here gives 
$$
\begin{aligned}
\|u^\varepsilon(t)\|^{2}_{H^\sigma}&\lesssim \|u_0\|_{H^\sigma}^{2}+\varepsilon^{4}\|\int_{0}^{t}{\rm e}^{{\bf i}(t-s)\Delta}\Phi(x)u^\varepsilon(s)\,ds\|_{H^\sigma}^{2}+\varepsilon^{2}
\|\int_{0}^{t}{\rm e}^{{\bf i}(t-s)\Delta}Q^{\frac{1}{2}}\,dW(s)\|_{H^\sigma}^{2}\\
&\le \|u_0\|_{H^\sigma}^{2}+\varepsilon^{4} t \int_{0}^{t}\|\Phi(x)u^\varepsilon(s)\|^2_{H^\sigma}\,ds+\varepsilon^{2}
\|\int_{0}^{t}{\rm e}^{{\bf i}(t-s)\Delta}Q^{\frac{1}{2}}\,dW(s)\|_{H^\sigma}^{2}\\
&\le \|u_0\|_{H^\sigma}^{2}+\varepsilon^{4}t \|\Phi\|_{H^\sigma}^{2}\int_{0}^{t}\|u^\varepsilon(t)\|^{2}_{H^\sigma}\,dt+\varepsilon^{2}
\|\int_{0}^{t}{\rm e}^{{\bf i}(t-s)\Delta}Q^{\frac{1}{2}}\,dW(s)\|_{H^\sigma}^{2}\\
&\lesssim \|u_0\|_{H^\sigma}^{2}+\varepsilon^2 \|\Phi\|_{H^\sigma}^{2}\int_{0}^{t}\sup_{s\in[0,t]}\|u^\varepsilon(s)\|^{2}_{H^\sigma}\,dt+\varepsilon^{2}
\|\int_{0}^{t}{\rm e}^{{\bf i}(t-s)\Delta}Q^{\frac{1}{2}}\,dW(s)\|_{H^\sigma}^{2}.
\end{aligned}
$$
Taking the supremum and the expectation both sides and using BDG inequality we then get
$$
\mathbb{E}\big[\displaystyle\sup_{t\in [0,T_\varepsilon]}\|u^\varepsilon(t)\|^{2}_{H^\sigma}\big]\lesssim \|u_0\|^{2}_{ H^\sigma}
+\varepsilon^{2}\|\Phi\|_{H^\sigma}^{2} \int_{0}^{T_\varepsilon}\mathbb{E}\big[\sup_{s\in[0,t]}\|u^\varepsilon(s)\|^{2}_{H^\sigma}\big]\,dt+\|Q^{\frac{1}{2}}\|_{\mathcal{L}_2^\sigma}^2.
$$
An application of Gronwall inequality then gives the result.
\end{proof}

\begin{remark}
    The smallness assumption and the linearity of the Equation \eqref{lin_small} allows us to prove that the $L^2(\Omega; L^2)$ error of the integrator (SLR1) remains buounded and of size $\mathcal{O}(\varepsilon\tau)$ for time windows of length $\mathcal{O}(\varepsilon^{-1})$ and remains bounded and of size $\mathcal{O}(\tau)$ for time windows of length $\mathcal{O}(\varepsilon^{-2})$. However, the result in Proposition \ref{smal_prop} and the presence of an $\mathcal{O}(\varepsilon^2)$ potential suggests to provide an alternative numerical integrator to \eqref{lin_small} that exhibits a bounded error of size $\mathcal{O}(\varepsilon^2 \tau)$ up to intervals of length $\mathcal{O}(\varepsilon^{-2})$. This motivate the upcoming section.
\end{remark}
\subsection{A non-resonant low-regularity stochastic numerical integrator}
We recall that the derivation of the integrator (SLR1) relied on the following approximation
\begin{equation}
    \label{app0}
\int_{0}^{\tau}{\rm e}^{{\bf i}s(\ell^2-\ell_2^2)}\,ds \approx \int_{0}^{\tau}{\rm e}^{{\bf i}s\ell_1^2}\,ds 
\end{equation}
and we computed exactly the integral in the right hand-side of \eqref{app0}. Inspired to the idea in \cite{fe_ka}, the construction of a non-resonant integrator for  \eqref{lin_small} relies on the fact that the integral involved in \eqref{app2} can be computed exactly in the resonant case $\ell^2-\ell_2^2=0$. With this in hand and with approximation \eqref{app0} for the non resonant case $\ell^2-\ell_2^2\ne 0$, also taking into account $v^\varepsilon(t)={\rm e}^{-{\bf i}t\Delta}u^\varepsilon(t)$, we obtain
\begin{align*}
\int_{0}^{\tau}{\rm e}^{-{\bf i}(t_n+s)\Delta}\big(\Phi(x){\rm e}^{{\bf i}(t_n+s)\Delta} v^\varepsilon(t_n)\big)\,ds
&\approx \frac{1}{2\pi}\displaystyle\sum_{\ell \in \mathbb{Z}}{\rm e}^{{\bf i}\ell x}\displaystyle\sum_{\substack{\ell_1,\ell_2\in\mathbb{Z}\\ \ell_1+\ell_2=\ell \\ \ell^2-\ell_2^2\ne 0}}\hat{\Phi}_{\ell_1}\widehat{v^\varepsilon}_{\ell_2}(t_n){\rm e}^{{\bf i}t_n(\ell^2-\ell_2^2)}\int_{0}^{\tau}
{\rm e}^{{\bf i}s \ell_1^2}\,ds\\
&\hspace{2cm}+\frac{\tau}{2\pi} \displaystyle\sum_{\ell \in \mathbb{Z}}{\rm e}^{{\bf i}\ell x} \displaystyle\sum_{\substack{\ell_1,\ell_2\in\mathbb{Z}\\ \ell_1+\ell_2=\ell\\ \ell^2-\ell^2_2=0}} \hat{\Phi}_{\ell_1}\widehat{v^\varepsilon}_{\ell_2}(t_n).
\end{align*}
This can be equivalently written as
\begin{align*}
    \int_{0}^{\tau}{\rm e}^{-{\bf i}(t_n+s)\Delta}\big(\Phi(x){\rm e}^{{\bf i}(t_n+s)\Delta} v^\varepsilon(t_n)\big)\,ds & \approx \frac{1}{2\pi}
    \displaystyle\sum_{\ell \in \mathbb{Z}} {\rm e}^{{\bf i}\ell x}\displaystyle\sum_{\substack{\ell_1,\ell_2\in\mathbb{Z}\\ \ell_1+\ell_2=\ell}} {\rm e}^{{\bf i}t_n (\ell^2-\ell_2^2)}\hat{\Phi}_{\ell_1}\widehat{v^\varepsilon}_{\ell_2}(t_n) \tau \varphi_1({\bf i}\tau \ell_1^2)\\
    &\hspace{2cm}+\frac{1}{2\pi}\displaystyle\sum_{\ell\in\mathbb{Z}}{\rm e}^{{\bf i}\ell x}\displaystyle\sum_{\substack{\ell_1,\ell_2\in\mathbb{Z} \\ \ell_1+\ell_2=\ell\\ \ell^2-\ell_2^2=0}} \hat{\Phi}_{\ell_1}\widehat{v^\varepsilon}_{\ell_2}(t_n)\tau (1-\varphi_1({\bf i}\tau \ell_1^2))
\end{align*}
Since $\ell^2-\ell_2^2=(\ell+\ell_2)(\ell-\ell_2)$ and $\varphi_1(0)=1$, we rewrite the second term of the above equation as
\begin{align*}
\frac{1}{2\pi}\displaystyle\sum_{\ell\in\mathbb{Z}}{\rm e}^{{\bf i}\ell x}\displaystyle\sum_{\substack{\ell_1,\ell_2\in\mathbb{Z} \\ \ell_1+\ell_2=\ell\\ \ell^2-\ell_2^2=0}} \hat{\Phi}_{\ell_1}\widehat{v^\varepsilon}_{\ell_2}(t_n)\tau (1-\varphi_1({\bf i}\tau \ell_1^2))
=&\frac{1}{2\pi} \displaystyle\sum_{\ell\in\mathbb{Z}}{\rm e}^{{\bf i}\ell x}\displaystyle\sum_{\substack{\ell_1,\ell_2\in\mathbb{Z}\\ \ell_1=0\\ \ell_2=\ell}} \hat{\Phi}_{\ell_1}\widehat{v^\varepsilon}_{\ell_2}(t_n)\tau (1-\varphi_1({\bf i}\tau \ell_1^2))\\
&\hspace{1cm}+\frac{1}{2\pi}\displaystyle\sum_{\ell\in\mathbb{Z}}{\rm e}^{{\bf i}\ell x}\displaystyle\sum_{\substack{\ell_1,\ell_2\in\mathbb{Z} \\ \ell_1=2\ell\\ \ell_2=-\ell}} \hat{\Phi}_{\ell_1}\widehat{v^\varepsilon}_{\ell_2}(t_n)\tau (1-\varphi_1({\bf i}\tau \ell_1^2))\\
&=\frac{\tau}{\sqrt{2\pi}} (1-\varphi_1(0))\hat{\Phi}_0 v^\varepsilon(t_n)\\
&\hspace{1cm}+\frac{\tau}{2\pi} \displaystyle\sum_{\ell\in\mathbb{Z}}{\rm e}^{{\bf i}\ell x} \hat{\Phi}_{2 \ell}\widehat{v^\varepsilon}_{-\ell}(t_n)(1-\varphi_1(4{\bf i}\tau \ell^2))\\
&=\frac{\tau}{2\pi} \displaystyle\sum_{\ell\in\mathbb{Z}}{\rm e}^{{\bf i}\ell x} \hat{\Phi}_{2 \ell}\widehat{v^\varepsilon}_{-\ell}(t_n)(1-\varphi_1(4{\bf i}\tau \ell^2)).
\end{align*}
Then, defining $g_2(u,v)$ through its Fourier coefficients as $\widehat{(g_2(w,z))}_\ell = \frac{1}{\sqrt{2\pi}}\hat{w}_{2\ell} \hat{z}_{-\ell}$, then, we get the following integrator for the twisted variable $v$
\begin{align}
v^\varepsilon_{n+1}=v^\varepsilon_n-{\bf i}\varepsilon^2 \tau \big[{\rm e}^{-{\bf i}t_n\Delta }\big(({\rm e}^{{\bf i}t_n\Delta }v^\varepsilon_n )\varphi_1(-{\bf i}\tau \Delta)\Phi(x)\big)&+ (1-\varphi_1(-4{\bf i}\tau \Delta))g_2(\Phi(x),v^\varepsilon_n)\big]\notag \\
&+\varepsilon\int_{t_n}^{t_n+\tau}{\rm e}^{-{\bf i}s \Delta}Q^{\frac{1}{2}}\,dW(s)\label{twis_met_lin_small}.
\end{align}
Return back to the original variable $u$, we get 
\begin{align*}
u^\varepsilon_{n+1}&={\rm e}^{{\bf i}\tau \Delta}\big[u^\varepsilon_n-{\bf i} \varepsilon^2 \tau \big(u^\varepsilon_n\varphi_1(-{\bf i}\tau \Delta) \Phi(x)\\
&\hspace{2cm}+(1-\varphi_1(-4{\bf i}\tau \Delta)){\rm e}^{{\bf i}t_n \Delta} g_2(\Phi(x),{\rm e}^{-{\bf i}t_n\Delta}u^\varepsilon_n)\big)\big]+\varepsilon \mathcal{W}_n.
\end{align*}
Finally, since
\begin{align*}
    {\rm e}^{{\bf i}t_n \Delta} g_2(\Phi(x),{\rm e}^{-{\bf i}t_n\Delta}u^\varepsilon_n)& = \frac{1}{2\pi}
    {\rm e}^{{\bf i}t_n \Delta} \displaystyle\sum_{\ell\in\mathbb{Z}}{\rm e}^{{\bf i}\ell x}\hat{\Phi}_{2\ell}{\rm e}^{{\bf i}\ell^2 t_n}\widehat{u^\varepsilon}_{-\ell ,n}\\
    &= \frac{1}{2\pi}\displaystyle\sum_{\ell\in\mathbb{Z}}{\rm e}^{{\bf i}\ell x}\hat{\Phi}_{2\ell}\widehat{u^\varepsilon}_{-\ell,n}\\
    &=g_2(\Phi(x),u^\varepsilon_n),
\end{align*}
we get the following stochastic accelerated non-resonant low-regularity integrator for \eqref{lin_small}
\begin{align}
u^\varepsilon_{n+1}={\rm e}^{{\bf i}\tau\Delta} u^\varepsilon_n&-{\bf i}\varepsilon^2 \tau {\rm e}^{{\bf i}\tau \Delta} \big(u^\varepsilon_n\varphi_1(-{\bf i}\tau\Delta)\Phi(x)\big)\notag \\
&-{\bf i}\varepsilon^2 \tau (1-\varphi_1(-4{\bf i}\tau\Delta)){\rm e}^{{\bf i}\tau\Delta}g_2(\Phi(x),u^ \varepsilon_n)+\varepsilon \mathcal{W}_n \label{non_res1}.
\end{align}
We will call the integrator in \eqref{non_res1} (SNRLR1).
\begin{remark}
    It is possible to show that all the results obtained in Section \ref{sec_lr} for the integrator (SLR1) holds true also for (SNRLR1). In what follows, we provide a long-term error analysis for the latter.
\end{remark}

Recall that $v^\varepsilon$ stands for the twisted variable efined through \eqref{twis_def}. Let $\varphi^{\varepsilon,\tau}$ the flow map such that $v^\varepsilon(t_n+\tau)=\varphi^{\varepsilon,\tau}(v(t_n))$, i.e., 
\begin{align}
\varphi^{\varepsilon,\tau} (v^\varepsilon(t_n))=v^\varepsilon(t_n)&-{\bf i}\varepsilon^2\int_{0}^{\tau} {\rm e}^{-{\bf i}(t_n+s)\Delta} \big(\Phi(x){\rm e}^{{\bf i}(t_n+s)\Delta} v^\varepsilon(t_n+s)\big)\,ds\notag \\
&+\varepsilon \int_{t_n}^{t_n+\tau}{\rm e}^{-{\bf i}s \Delta} Q^{\frac{1}{2}}\,dW(s)     \label{twis_sol_lin_small]}
\end{align}
and let $\Psi_{t_n}^{\varepsilon,\tau}$ the map such that $v^\varepsilon_{n+1}=\Psi_{t_n}^{\varepsilon,\tau}(v^\varepsilon_n)$. 
%i.e., 
%$$
%\begin{aligned}
%\Phi_{t_n}^\tau(v(t_n))=v(t_n)-{\bf i}\varepsilon\tau \big[{\rm e}^{{\bf i}t_n\Delta }\big(({\rm e}^{{-\bf i}t_n\Delta }v(t_n))\varphi_1({\bf i}\tau \Delta)f(x)\big)&+ (1-\varphi_1(4{\bf i}\tau \Delta))g(f(x),v(t_n))\big]\\
%&-{\bf i}\varepsilon\int_{t_n}^{t_n+\tau}{\rm e}^{{\bf i}s \Delta}Q^{\frac{1}{2}}\,dW(s).
%\end{aligned}
%$$
We give the following theorem.
\begin{theorem}
\label{lobg_term_res1}
Let $\sigma\in\mathbb{N}^+$. Let $u_0\in H^\sigma$, $\Phi\in H^{\frac{3}{2}+}$ if $\sigma=1$ or $\Phi\in H^\sigma$ if $\sigma\ge 2$, and $Q^{\frac{1}{2}}\in\mathcal{L}_2^\sigma$. Then, fixed $\tau_0\in (0,1)$ independently from $\varepsilon$, there exists a sufficiently small $\tau$ (depending only on $\tau_0$) such that the numerical integrator (SNRLR1), applied to \eqref{lin_small}, satisfies the following error estimate 
\begin{equation}
    \label{re1}
\|u^\varepsilon(t_n)-u^\varepsilon_n\|_{L^2(\Omega; L^2)}\lesssim \tau_0^\sigma+\varepsilon^2 \tau, \qquad \qquad 0\le n \le \frac{T}{\varepsilon^2 \tau}.
\end{equation}
%If in addition $u_0\in L^{2p}(\Omega; H^2)$, for all $p\ge 1$, then, asymptotically, for any $t_n\in[0,\varepsilon^{-1}]$ and $q\in[0,\infty)$, one has
%\begin{equation}
 %   \label{res2}
%\|v(t_n)-v_n\|_{H^1}\le \varepsilon^q \tau^r, \qquad \mathbb{P}-{\rm a.s.}
%\end{equation}
%with any $r<1$.
\end{theorem}
\begin{proof}
   We define the quantities $\mathcal{E}^\varepsilon_n:=\Psi_{t_n}^{\varepsilon,\tau}(v^\varepsilon(t_n))-\varphi_{t_n}^{\varepsilon,\tau}(v^\varepsilon(t_n))$ and $e^\varepsilon_{n+1}:=v^\varepsilon_{n+1}-v^\varepsilon(t_{n+1})$. Then, we have
\begin{align*}
e^\varepsilon_{n+1} &= v^\varepsilon_{n+1}-\varphi_{t_{n}}^{\varepsilon,\tau} (v^\varepsilon(t_n))\\
&=v^\varepsilon_{n+1}-\Psi_{t_n}^{\varepsilon,\tau} (v^\varepsilon(t_n))+\mathcal{E}^\varepsilon_n\\
&=e^\varepsilon_n+Z^\varepsilon_n+\mathcal{E}^\varepsilon_n,
\end{align*}
where
\begin{align*}
    Z^\varepsilon_n&=-{\bf i}\varepsilon^2\tau\big[{\rm e}^{-{\bf i}t_n\Delta }\big(({\rm e}^{{\bf i}t_n\Delta }v^\varepsilon_n)\varphi_1(-{\bf i}\tau \Delta)\Phi(x)\big)+ (1-\varphi_1(-4{\bf i}\tau \Delta))g_2(\Phi(x),v^\varepsilon_n)\\
    &\hspace{2cm}-{\rm e}^{-{\bf i}t_n\Delta }\big(({\rm e}^{{\bf i}t_n\Delta }v^\varepsilon(t_n))\varphi_1(-{\bf i}\tau \Delta)\Phi(x)\big)+ (1-\varphi_1(-4{\bf i}\tau \Delta))g_2(\Phi(x),v^\varepsilon(t_n))\big]\\
    &= Z^1_{n}+Z^2_{n},
\end{align*}
where
\begin{align*}
&Z_n^{\varepsilon,1}=-{\bf i}\varepsilon^2\tau {\rm e}^{-{\bf i}t_n\Delta}\big(({\rm e}^{{\bf i}t_n \Delta}(v^\varepsilon_n-v^\varepsilon(t_n)))\varphi_1(-{\bf i}\tau\Delta)\Phi(x)\big)\\
&Z_n^{\varepsilon,2}=-{\bf i}\varepsilon^2\tau (1-\varphi_1(-4{\bf i}\tau \Delta))\big(g_2(\Phi(x),v^\varepsilon_n)-g_2(\Phi(x),v^\varepsilon(t_n))\big).
\end{align*}
First, direct computations show 
\begin{align*}
\|Z_n^{\varepsilon,1}\| &= \varepsilon^{2} \tau \|{\rm e}^{-{\bf i}t_n\Delta}\big(({\rm e}^{{\bf i}t_n \Delta}(v^\varepsilon_n-v^\varepsilon(t_n)))\varphi_1(-{\bf i}\tau\Delta)\Phi(x)\big)\|^{2}\\
&=\varepsilon^2 \tau \|({\rm e}^{{\bf i}t_n \Delta}(v^\varepsilon_n-v^\varepsilon(t_n)))\varphi_1(-{\bf i}\tau\Delta)\Phi(x)\|\\
&\le \varepsilon^{2} \tau \|\varphi_1(-{\bf i}\tau\Delta)\Phi(x)\|_{L^\infty}\|e^\varepsilon_n\|\\
&\lesssim \varepsilon^{2} \tau  \|\Phi\|_{H^{\frac{1}{2}+}}\|e^\varepsilon_n\|
\end{align*}
and
\begin{align*}
    \|Z_n^{\varepsilon,2}\|&\lesssim \varepsilon^{2} \tau \|g_2(\Phi(x),v^\varepsilon_n)-g_2(\Phi(x),v^\varepsilon(t_n))\|\\
    &= \varepsilon^{2} \tau  (\sum_{\ell\in\mathbb{Z}}|\hat{\Phi}_{2\ell}|^2 |\widehat{v^\varepsilon}_{-\ell,n}-\widehat{v^\varepsilon}_{-\ell}(t_n)|^2)^{\frac{1}{2}}\\
    &\le \varepsilon^{2}\tau \big(\max_{\ell\in\mathbb{Z}}|\Phi_{2\ell}|^2\big)^{\frac{1}{2}} \|e^\varepsilon_n\|\\
    &\le  \varepsilon^{2}\tau \|\Phi\| \|e^\varepsilon_n\|.
\end{align*}
This then yields
$$
\mathbb{E}\|Z^\varepsilon_n\|^{2}\lesssim \varepsilon^{4} \tau^{2} \mathbb{E}\|e^\varepsilon_n\|^{2}.
$$
We now have to estimate $\mathbb{E}\|e^\varepsilon_{n+1}\|^2.$ We note that $\mathcal{E}^\varepsilon_n$ can be decomposed as
$\mathcal{E}^\varepsilon_n=\mathcal{V}^\varepsilon_n+\mathcal{R}^\varepsilon(v^\varepsilon(t_n))$, where
\begin{align*}
&\mathcal{R}^\varepsilon(v^\varepsilon(t_n))={\bf i}\varepsilon^2 \displaystyle\sum_{\ell\in\mathbb{Z}}{\rm e}^{{\bf i}\ell x}\displaystyle\sum_{\substack{\ell_1,\ell_2\in\mathbb{Z}\\ \ell_1+\ell_2=\ell\\\ell^2-\ell_2^2\ne 0}}{\rm e}^{{\bf i}t_n (\ell^2-\ell_2^2)}\hat{\Phi}_{\ell_1}\widehat{v^\varepsilon}_{\ell_2}(t_n)\int_{0}^{\tau}{\rm e}^{{\bf i}s \ell_1^2}-{\rm e}^{{\bf i}s(\ell^2-\ell_2^2)}\,ds,\\
&\mathcal{V}_n^\varepsilon={\bf i}\varepsilon^2 \int_{0}^{\tau} {\rm e}^{-{\bf i}(t_n+s)\Delta} \big(\Phi(x){\rm e}^{{\bf i}(t_n+s)\Delta} v^\varepsilon(t_n+s)\big)\,ds-{\bf i}\varepsilon^2 \int_{0}^{\tau} {\rm e}^{-{\bf i}(t_n+s)\Delta} \big(\Phi(x){\rm e}^{{\bf i}(t_n+s)\Delta} v^\varepsilon(t_n)\big)\,ds.
\end{align*}
By recursion, we have
\begin{align*}
\mathbb{E}\|e^\varepsilon_{n+1}\|^{2}&=\mathbb{E}\|\sum_{k=0}^n Z^\varepsilon_k+\mathcal{V}^\varepsilon_k+\mathcal{R}^\varepsilon(v^\varepsilon(t_k))\|^{2}\\
&\lesssim \mathbb{E}\|\sum_{k=0}^n Z^\varepsilon_k\|^{2}+\mathbb{E}\|\sum_{k=0}^n \mathcal{V}^\varepsilon_k\|^{2}
+\mathbb{E}\|\sum_{k=0}^n \mathcal{R}^\varepsilon(v^\varepsilon(t_k))\|^{2}\\
&\lesssim (n+1)\varepsilon^{4} \tau^{2}\sum_{k=0}^{n}\mathbb{E}\|e^\varepsilon_k\|^{2}+\mathbb{E}\|\sum_{k=0}^n \mathcal{V}^\varepsilon_k\|^{2}+
\mathbb{E}\|\sum_{k=0}^n \mathcal{R}^\varepsilon(v^\varepsilon(t_k))\|^{2}\\
&\lesssim \varepsilon^2\tau \sum_{k=0}^{n}\mathbb{E}\|e^\varepsilon_k\|^{2}+\mathbb{E}\|\sum_{k=0}^n \mathcal{V}^\varepsilon_k\|^{2}+
\mathbb{E}\|\sum_{k=0}^n \mathcal{R}^\varepsilon(v^\varepsilon(t_k))\|^{2}.
\end{align*}
We now estimate the term $\mathbb{E}\|\displaystyle\sum_{k=0}^{n}\mathcal{V}^\varepsilon_k\|^{2}$. It is direct to see that
$$
\mathbb{E}\|\displaystyle\sum_{k=0}^{n}\mathcal{V}^\varepsilon_k\|^{2}\lesssim (n+1)\displaystyle\sum_{k=0}^{n}\mathbb{E}\|\mathcal{V}^\varepsilon_{k,1}\|^{2}+\mathbb{E}\|\displaystyle\sum_{k=0}^{n} \mathcal{V}^\varepsilon_{k,2}\|^{2},
$$
where
\begin{align*}
  &  \mathbb{E}\|\mathcal{V}^\varepsilon_{k,1}\|^{2}=\varepsilon^{8}\mathbb{E}\big\|\int_{0}^{\tau}{\rm e}^{-{\bf i}(t_k+s)\Delta} \Phi(x){\rm e}^{{\bf i}(t_k+s)\Delta}\big(\int_{0}^{s}{\rm e}^{-{\bf i}(t_k+r)\Delta}\Phi(x){\rm e}^{{\bf i}(t_k+r)\Delta}v^\varepsilon(t_k+r)\,dr\big)\,ds\big\|^{2},\\
&\mathbb{E}\|\displaystyle\sum_{k=0}^{n}\mathcal{V}^\varepsilon_{k,2}\|^{2}=\varepsilon^{6}\mathbb{E}\big\|\displaystyle\sum_{k=0}^n\int_{0}^{\tau}{\rm e}^{-{\bf i}(t_k+s)\Delta} \Phi(x){\rm e}^{{\bf i}(t_k+s)\Delta}\big(\int_{t_k}^{t_k+s}{\rm e}^{-{\bf i}r \Delta}Q^{\frac{1}{2}}\,dW(r)\big)\,ds\big\|^{2}.
\end{align*}
Using twice the Hölder inequality, the fact that ${\rm e}^{{\bf i}t\Delta}$ is an isometry on $L^2$ and Proposition \ref{smal_prop},  we have
\begin{align*}
\mathbb{E}\|\mathcal{V}^\varepsilon_{k,1}\|^{2} & \le \tau \varepsilon^8  \mathbb{E}\int_{0}^{\tau}\|\Phi(x){\rm e}^{{\bf i}(t_k+s)\Delta}\big(\int_{0}^{s}{\rm e}^{-{\bf i}(t_k+r)\Delta}\Phi(x){\rm e}^{{\bf i}(t_k+r)\Delta}v^\varepsilon(t_k+r)\,dr\big)\|^{2}\,ds\\
&\lesssim \tau \varepsilon^{8} \|\Phi\|_{L^\infty}^{2} \mathbb{E}\int_{0}^{\tau}\|\int_{0}^{s}{\rm e}^{-{\bf i}(t_k+r)\Delta}\Phi(x){\rm e}^{{\bf i}(t_k+r)\Delta}v^\varepsilon(t_k+r)\,dr\|^{2}\,ds\\
&\le \tau^{2} \varepsilon^{8} \|\Phi\|^{4}_{L^\infty} \mathbb{E}\int_{0}^{\tau}\int_{0}^{s}\|v^\varepsilon(t_k+r)\|^{2}\,dr\,ds\\
&\lesssim \tau^4 \varepsilon^8.
\end{align*}
For $\mathbb{E}\|\displaystyle\sum_{k=0}^n\mathcal{V}^\varepsilon_{k,2}\|^{2}$, the stochastic Fubini theorem and Itô isometry allow the following computation
\begin{align*}
\mathbb{E}\|\displaystyle\sum_{k=0}^n\mathcal{V}^\varepsilon_{k,2}\|^{2}&=\varepsilon^{6}\mathbb{E}\big\|\displaystyle\sum_{k=0}^{n}\int_{0}^{\tau}{\rm e}^{-{\bf i}(t_k+s)\Delta} \Phi(x){\rm e}^{{\bf i}(t_k+s)\Delta}\bigg(\int_{t_k}^{t_k+s}{\rm e}^{-{\bf i}r \Delta}Q^{\frac{1}{2}}\,dW(r)\bigg)\,ds\big\|^{2}\\
    &=\varepsilon^{6}\mathbb{E}\big\|\displaystyle\sum_{k=0}^{n}\int_{t_k}^{t_{k+1}}\int_{r}^{t_{k+1}}{\rm e}^{-{\bf i}s\Delta} \Phi(x){\rm e}^{{\bf i}s\Delta}{\rm e}^{-{\bf i}r \Delta}Q^{\frac{1}{2}}\,ds\,dW(r)\big\|^{2}\\
     &=\varepsilon^{6}\mathbb{E}\big\|\int_{0}^{t_{n+1}}\int_{r}^{[\frac{r}{\tau}+1]\tau}{\rm e}^{-{\bf i}s\Delta} \Phi(x){\rm e}^{{\bf i}s\Delta}{\rm e}^{-{\bf i}r \Delta}Q^{\frac{1}{2}}\,ds\,dW(r)\big\|^{2}\\
     &=\varepsilon^{6} \mathbb{E}\big[\displaystyle\sum_{k=0}^{n}\int_{t_k}^{t_{k+1}}\big\|\int_{r}^{t_{k+1}}{\rm e}^{-{\bf i}s\Delta} \Phi(x){\rm e}^{{\bf i}s\Delta}{\rm e}^{-{\bf i}r \Delta}Q^{\frac{1}{2}}\,ds\big\|^2_{\mathcal{L}_2^0}\,dr\big]\\
     &\lesssim \varepsilon^{6}\mathbb{E}\big[\displaystyle\sum_{k=0}^{n}\int_{t_k}^{t_{k+1}}\big(\int_{r}^{t_{k+1}}\big\|{\rm e}^{-{\bf i}s\Delta} \Phi(x){\rm e}^{{\bf i}s\Delta}{\rm e}^{-{\bf i}r \Delta}Q^{\frac{1}{2}}\big\|_{\mathcal{L}_2^0}\,ds\big)^2\,dr\big]\\
     &\lesssim \varepsilon^{6}\|\Phi\|^{2}_{L^\infty}\big\|Q^{\frac{1}{2}}\big\|^{2}_{\mathcal{L}_2^0}\displaystyle\sum_{k=0}^{n}\int_{t_{k}}^{t_{k+1}}|t_{k+1}-r|^2\,dr\\
     &\lesssim  \varepsilon^6 (n+1)\tau^3\\
     &\lesssim \varepsilon^4 \tau^2.
\end{align*}
The above computations finally give
$$
\mathbb{E}\|\displaystyle\sum_{k=0}^{n}\mathcal{V}^\varepsilon_k\|^{2}\lesssim \varepsilon^4\tau^2.
$$
Then, we achieve 
\begin{equation}
    \label{fe}
\mathbb{E}\|e^\varepsilon_{n+1}\|^{2}\lesssim\varepsilon^2 \tau \sum_{k=0}^{n}\mathbb{E}\|e^\varepsilon_k\|^{2}+\varepsilon^4 \tau^2 +\mathbb{E}\|\sum_{k=0}^n \mathcal{R}^\varepsilon(v^\varepsilon(t_k))\|^{2}.
\end{equation}
Inspired from \cite{bao,fe_ka}, we now exploit (RCO) technique to obtain an improved bound for $\mathbb{E}\|\displaystyle\sum_{k=0}^n \mathcal{R}^\varepsilon(v^\varepsilon(t_k))\|^2$. Fix a cut-off parameter $\tau_0\in(0.1)$. Let $N_0=2[1/\tau_0]$ and let $\mathcal{P}_{N_0}$ the Projection operator in the space spanned by frequencies $[-\frac{N_0}{2},\frac{N_0}{2}]$. We first estimate $\|\sum_{k=0}^n \mathcal{R}^\varepsilon(v^\varepsilon(t_k))-\mathcal{P}_{N_0}\mathcal{R}^\varepsilon(v^\varepsilon(t_k))\|^2$. Using the algebra property of $H^\sigma$ and a direct calculation, we have
\begin{align*}
\|\sum_{k=0}^n \mathcal{R}^\varepsilon(v^\varepsilon(t_k))-\mathcal{P}_{N_0}\mathcal{R}^\varepsilon(v^\varepsilon(t_k))\|^2&\le (n+1)\displaystyle\sum_{k=0}^{n}\|\mathcal{R}^\varepsilon(v^\varepsilon(t_k))-\mathcal{P}_{N_0}\mathcal{R}^\varepsilon(v^\varepsilon(t_k))\|^2\\
&\le (n+1)\varepsilon^4\displaystyle\sum_{k=0}^{n}\displaystyle\sum_{|\ell|\ge \frac{N_0}{2}+1}\big|\displaystyle\sum_{\substack{\ell_1,\ell_2\in\mathbb{Z}\\ \ell_1+\ell_2=\ell\\\ell^2-\ell_2^2\ne 0}}{\rm e}^{{\bf i}t_k (\ell^2-\ell_2^2)}\hat{\Phi}_{\ell_1}\widehat{v^\varepsilon}_{\ell_2}(t_k)\\
&\hspace{2cm}\times \int_{0}^{\tau}{\rm e}^{{\bf i}s \ell_1^2}-{\rm e}^{{\bf i}s(\ell^2-\ell_2^2)}\,ds\big|^2\\
&\le (n+1)\varepsilon^4\tau^2 \displaystyle\sum_{k=0}^{n} \displaystyle\sum_{|\ell|\ge \frac{N_0}{2}+1}\big(\displaystyle\sum_{\substack{\ell_1,\ell_2\in\mathbb{Z}\\ \ell_1+\ell_2=\ell}}
|\hat{\Phi}_{\ell_1}||\widehat{v^\varepsilon}_{\ell_2}(t_k)|\big)^2\\
&\le (n+1)\varepsilon^4\tau^2 \tau_0^{2\sigma} \displaystyle\sum_{k=0}^n \|\Phi(x)v^\varepsilon(t_k)\|^2_{H^\sigma}\\
&\le (n+1)\varepsilon^4\tau^2 \tau_0^{2\sigma} \|\Phi\|^2_{H^\sigma}\displaystyle\sum_{k=0}^n \|v^\varepsilon(t_k)\|^2_{H^\sigma}.
\end{align*}
Then, thanks to Proposition \ref{smal_prop}, we get
\begin{align}
    \mathbb{E}\big[\|\sum_{k=0}^n \mathcal{R}^\varepsilon(v^\varepsilon(t_k))-\mathcal{P}_{N_0}\mathcal{R}^\varepsilon(v^\varepsilon(t_k))\|^2\big]&\lesssim (n+1)^2\varepsilon^4\tau^2 \tau_0^{2\sigma}\notag \\
    &\lesssim  \tau_0^{2\sigma}\label{fe0}.
\end{align}
We now estimate the term $\mathbb{E}\big[\|\sum_{k=0}^n \mathcal{P}_{N_0}\mathcal{R}^\varepsilon(v^\varepsilon(t_k))\|^2\big]$. By the definition of $\mathcal{R}^\varepsilon(v^\varepsilon(t_k))$, we see
$$
\sum_{k=0}^n \mathcal{P}_{N_0}\mathcal{R}^\varepsilon(v^\varepsilon(t_k))={\bf i}\varepsilon^2 \sum_{k=0}^n \displaystyle\sum_{\ell \in \mathcal{T}_{N_0}}{\rm e}^{{\bf i}\ell x}\displaystyle\sum_{(\ell_1,\ell_2)\in\mathcal{I}_\ell^{N_0}} \Lambda^\varepsilon_{k,\ell,\ell_1,\ell_2},
$$
where
$
\mathcal{T}_{N_0}:=\{\ell: -\frac{N_0}{2},\dots,\frac{N_0}{2}\}
$
and, for $\ell \in \mathcal{T}_{N_0}$, 
$$
\mathcal{I}_\ell^{N_0}:=\{(\ell_1,\ell_2): \ell_1+\ell_2=\ell, \ell_1,\ell_2\in \mathcal{T}_{N_0}\}.
$$
Denoting $\delta_{\ell,\ell_1,\ell_2}=\ell^2-\ell_2^2$, note that $\Lambda^\varepsilon_{k,\ell,\ell_1,\ell_2}=0$ if $\delta_{\ell,\ell_1,\ell_2}=0$ and 
$$
\Lambda^\varepsilon_{k,\ell,\ell_1,\ell_2}=r_{\ell,\ell_1,\ell_2} {\rm e}^{{\bf i}t_k\delta_{\ell,\ell_1,\ell_2}} c^\varepsilon_{k,\ell,\ell_1,\ell_2}
$$
with
$$
r_{\ell,\ell_1,\ell_2}=\int_0^\tau {\rm e}^{{\bf i}s \ell_1^2}\,ds-\int_0^\tau {\rm e}^{{\bf i}s \delta_{\ell,\ell_1,\ell_2}}\,ds, \quad c^\varepsilon_{k,\ell,\ell_1,\ell_2}
=\hat{\Phi}_{\ell_1} \widehat{v^\varepsilon}_{\ell_2}(t_k).
$$
Note that a standard interpolation argument gives the bound
$$
r_{\ell,\ell_1,\ell_2}=\mathcal{O}(\tau^2 |\delta_{\ell,\ell_1,\ell_2}-\ell_1^2|).
$$
Moreover, for $\ell_\in \mathcal{T}_{N_0}$ and $\ell_1,\ell_2 \in \mathcal{I}_{\ell}^{N_0}$, one has
$$
|\delta_{\ell,\ell_1,\ell_2}|=|\ell^2-\ell_2^2|\le 2 \frac{N_0^2}{4}\le \frac{2(1+\tau_0)^2}{\tau_0^2}.
$$
This implies that for sufficiently small $\tau$ (depending on $\tau_0 \in (0,1)$) one has 
$$
\frac{\tau}{2}|\delta_{\ell,\ell_1,\ell_2}|<\pi.
$$
Then, denoting $S_{n,\ell,\ell_1,\ell_2}=\sum_{k=0}^n {\rm e}^{{\bf i}t_k\delta_{\ell,\ell_1,\ell_2}}$ using summation-by-parts formula (see also \cite{fe_ka}), one achieves
\begin{align*}
\sum_{k=0}^{n}\Lambda^\varepsilon_{k,\ell_1,\ell_2}&=r_{\ell,\ell_1,\ell_2}\sum_{k=0}^{n-1}S_{k,\ell,\ell_1,\ell_2}\hat{\Phi}_{\ell_1}\big(\widehat{v^\varepsilon}_{\ell_2}(t_k)-\widehat{v^\varepsilon}_{\ell_2}(t_{k+1})\big)\\
&\qquad\qquad +S_{n\ell,\ell_1,\ell_2}r_{\ell,\ell_1,\ell_2}\hat{\Phi}_{\ell_1}\widehat{v^\varepsilon}_{\ell_2}(t_n)\\
%&\lesssim \tau \frac{|\delta_{\ell,\ell_1,\ell_2}-\ell_1^2|}{|\delta_{\ell,\ell_1,\ell_2}|}\sum_{k=0}^{n-1}\big|\hat{\Phi}_{\ell_1}\big| \big|\hat{v}_{\ell_2}(t_k)-\hat{v}_{\ell_2}(t_{k-1})\big|\\
%&\qquad\qquad+\tau \frac{|\delta_{\ell,\ell_1,\ell_2}-\ell_1^2|}{|\delta_{\ell,\ell_1,\ell_2}|}\big|\hat{\Phi}_{\ell_1}\big|\big|\hat{v}_{\ell_2}(t_n)\big|.
&={\bf i}\varepsilon^2 r_{\ell,\ell_1,\ell_2}\sum_{k=0}^{n-1}S_{k,\ell,\ell_1,\ell_2}\hat{\Phi}_{\ell_1}\widehat{\mathcal{A}^\varepsilon}_{1,\ell_2}(t_k)-\varepsilon 
r_{\ell,\ell_1,\ell_2}\sum_{k=0}^{n-1}S_{k,\ell,\ell_1,\ell_2}\hat{\Phi}_{\ell_1}\hat{\mathcal{A}}_{2,\ell_2}(t_k)\\
&\qquad\qquad +S_{n\ell,\ell_1,\ell_2}r_{\ell,\ell_1,\ell_2}\hat{\Phi}_{\ell_1}\widehat{v^\varepsilon}_{\ell_2}(t_n),
\end{align*}
where $\widehat{\mathcal{A}^\varepsilon}_{1,\ell}$ denotes the Fourier coefficient of
$$
\mathcal{A^\varepsilon}_1=\int_0^\tau {\rm e}^{-{\bf i}(t_k+s)\Delta}\big( \Phi(x){\rm e}^{{\bf i}(t_k+s)\Delta}v^\varepsilon(t_k+s)\big) \,ds
$$
for the mode $\ell\in\mathbb{Z}$ and 
$$
\mathcal{\hat{A}}_{2,\ell}(t_k)=\int_{t_k}^{t_k+\tau}{\rm e}^{{\bf i}s \ell^2} q_\ell^{\frac{1}{2}}\,d\beta_\ell(s), \quad \ell\in\mathbb{Z}.
$$ 
\begin{comment}
Moreover, note that
\begin{align*}
\sum_{k=0}^{n-1} S_{k,\ell,\ell_1,\ell_2}\mathcal{\hat{A}}_{2,\ell_2}(t_k) &= \sum_{k=0}^{n-1} \sum_{j=0}^k {\rm e}^{{\bf i}t_j \delta_{\ell,\ell_1,\ell_2}}
\int_{t_k}^{t_k+\tau}{\rm e}^{{\bf i}s \ell_2^2} q_{\ell_2}^{\frac{1}{2}}\,d\beta_{\ell_2}(s)\\
&=\sum_{k=0}^{n-1}\int_{t_k}^{t_k+\tau} \sum_{j=0}^{[s/\tau]} {\rm e}^{{\bf i}j\tau \delta_{\ell,\ell_1,\ell_2}} {\rm e}^{{\bf i}s \ell_2^2} q_{\ell_2}^{\frac{1}{2}}\,d\beta_{\ell_2}(s)\\
&=\int_0^{t_n} \sum_{j=0}^{[s/\tau]} {\rm e}^{{\bf i}j\tau \delta_{\ell,\ell_1,\ell_2}} {\rm e}^{{\bf i}s \ell_2^2} q_{\ell_2}^{\frac{1}{2}}\,d\beta_{\ell_2}(s):=\mathcal{\hat{A}}_{\ell_2}(t_n),
\end{align*}
where
$$
\mathcal{A}(t_n)=\int_0^{t_n}\sum_{j=0}^{[s/\tau]} {\rm e}^{{\bf i}j\tau \delta_{\ell,\ell_1,\ell_2}} {\rm e}^{-{\bf i}s \Delta} Q^{\frac{1}{2}}\,dW(s).
$$
Then, we get
\begin{align*}
\sum_{k=0}^{n}\Lambda_{k,\ell_1,\ell_2,\ell_3}&={\bf i}\varepsilon^2 r_{\ell,\ell_1,\ell_2}\sum_{k=0}^{n-1}S_{k,\ell,\ell_1,\ell_2}\hat{\Phi}_{\ell_1}\mathcal{\hat{A}}_{1,\ell_2}(t_k)-\varepsilon 
r_{\ell,\ell_1,\ell_2}\hat{\Phi}_{\ell_1} \mathcal{\hat{A}}_{\ell_2}(t_n)\\
&\qquad\qquad +S_{n\ell,\ell_1,\ell_2}r_{\ell,\ell_1,\ell_2}\hat{\Phi}_{\ell_1}\hat{v}_{\ell_2}(t_n).
\end{align*}
\end{comment}
Then, we have
\begin{align*}
\mathbb{E}\big[\|\sum_{k=0}^{n}\mathcal{P}_{N_0}\mathcal{R}^\varepsilon(v^\varepsilon(t_k))\|^2\big]&= \varepsilon^4 \sum_{\ell \in \mathcal{T}_{N_0}}\mathbb{E}\big[\big|\sum_{(\ell_1,\ell_2)\in\mathcal{I}_{N_0}^{\ell}} \sum_{k=0}^n\Lambda^\varepsilon_{k,\ell,\ell_1,\ell_2}\big|^2\big]\\
&\lesssim \mathcal{R}^\varepsilon_1+\mathcal{R}^\varepsilon_2+\mathcal{R}^\varepsilon_3,
\end{align*}
where 
\begin{align*}
    &\mathcal{R}^\varepsilon_1=\varepsilon^8 \sum_{\ell\in \mathcal{T}_{N_0}}\mathbb{E}\big[\big|\sum_{(\ell_1,\ell_2)\in\mathcal{I}_{N_0}^\ell}r_{\ell,\ell_1,\ell_2} \sum_{k=0}^{n-1} S_{k,\ell,\ell_1,\ell_2}\hat{\Phi}_{\ell_1}\widehat{\mathcal{A}^\varepsilon}_{1,\ell_2}(t_k)\big|^2\big],\\
    &\mathcal{R}^\varepsilon_2=\varepsilon^6 \sum_{\ell\in \mathcal{T}_{N_0}}\mathbb{E}\big[\big|\sum_{(\ell_1,\ell_2)\in\mathcal{I}_{N_0}^\ell}r_{\ell,\ell_1,\ell_2} \sum_{k=0}^{n-1} S_{k,\ell,\ell_1,\ell_2}\hat{\Phi}_{\ell_1}\mathcal{\hat{A}}_{2,\ell_2}(t_k)\big|^2\big],\\
    &\mathcal{R}^\varepsilon_3=\varepsilon^4 \sum_{\ell\in \mathcal{T}_{N_0}}\mathbb{E}\big[\big|\sum_{(\ell_1,\ell_2)\in\mathcal{I}_{N_0}^\ell}
    S_{n,\ell,\ell_1,\ell_2}r_{\ell,\ell_1,\ell_2} \hat{\Phi}_{\ell_1}\widehat{v^\varepsilon}_{\ell_2}(t_n)\big|^2\big].
\end{align*}
By similar manner than \cite{fe_ka}, we observe that for any $n\ge 0$ and for $\ell\in\mathcal{T}_{N_0}$ and $\ell_1,\ell_2 \in \mathcal{I}_{N_0}^{\ell}$, one has the bound
\begin{align}
\big|r_{l,\ell_1,\ell_,2} S_{n,\ell,\ell_1,\ell_2}\big|& =\lesssim \tau \frac{|\delta_{\ell,\ell_1,\ell_2}-\ell_1^2|}{|\delta_{\ell,\ell_1,\ell_2}|}\notag \\
&\lesssim \tau \frac{|\ell_1 \ell_2|}{|\ell^2-\ell_2^2|}\lesssim \tau |\ell_2| \label{fe_bound}.
\end{align}
We are now ready to estimate the terms $\mathcal{R}^\varepsilon_j$, for $j=1,2,3$. For the term $\mathcal{R}_1$, exploiting the bound in \eqref{fe_bound}, the Young inequality, the fact that ${\rm e}^{{\bf i}t\Delta}$ si an isometry on $H^1$, the algebra property of $H^1$ and Proposition \ref{smal_prop}, we have
\begin{align*}
    \mathcal{R}^\varepsilon_1&\le \varepsilon^8 \sum_{\ell\in\mathcal{T}_{N_0}} \mathbb{E}\big[\big(\sum_{k=0}^{n-1} \sum_{\ell_2\ne \pm \ell} |r_{\ell,\ell-\ell_2,\ell_2}| |S_{k,\ell,\ell-\ell_2,\ell_2}| |\hat{\Phi}_{\ell-\ell_2}| |\widehat{\mathcal{A}^\varepsilon}_{1,\ell_2}(t_k)|\big)^2\big]\\
    &\le  \varepsilon^8 n \tau^2 \sum_{k=0}^{n-1} \mathbb{E} \big[\sum_{\ell\in\mathcal{T}_{N_0}}\big(\sum_{\ell_2\ne \pm \ell}|\ell_2| |\hat{\Phi}_{\ell-\ell_2}| |\widehat{\mathcal{A}^\varepsilon}_{1,\ell_2}(t_k)|\big)^2\big]\\
    &\le \varepsilon^8 \tau^2 n \sum_{k=0}^{n-1}\mathbb{E}\big[\|\Phi\star {\mathcal{A}^\varepsilon}_1^{(1)}(t_k)\|^2_{\ell^2}]\\
    &\lesssim  \varepsilon^8 \tau^2 n \|\Phi\|_{\ell^1}^2\sum_{k=0}^{n-1} \mathbb{E}\big[\|\mathcal{A}^\varepsilon_1(t_k)\|^2_{H^1}\big]\\
    &\lesssim \varepsilon^8 \tau^4 n^2\\
    &\lesssim \varepsilon^4 \tau^2.
\end{align*}
For the term $\mathcal{R}^\varepsilon_2$, an application of the Ito isometry and \eqref{fe_bound} give
\begin{align*}
\mathcal{R}^\varepsilon_2&=\varepsilon^6  \sum_{\ell\in \mathcal{T}_{N_0}}\mathbb{E}\big[\big|\displaystyle\sum_{\ell_2 \ne \pm\ell} \int_0^{t_n} r_{\ell,\ell-\ell_2,\ell_2}S_{[s/\tau],\ell,\ell-\ell_2,\ell_2}\hat{\Phi}_{\ell-\ell_2} {\rm e}^{{\bf i}s \ell_2^2}q_{\ell_2}^{\frac{1}{2}}\,d\beta_{\ell_2}(s)\big|^2\big]\\
&=\varepsilon^6  \sum_{\ell\in \mathcal{T}_{N_0}}\displaystyle\sum_{\ell_2 \ne \pm\ell} \mathbb{E}\big[\big| \int_0^{t_n} r_{\ell,\ell-\ell_2,\ell_2}S_{[s/\tau],\ell,\ell-\ell_2,\ell_2}\hat{\Phi}_{\ell-\ell_2} {\rm e}^{{\bf i}s \ell_2^2}q_{\ell_2}^{\frac{1}{2}}\,d\beta_{\ell_2}(s)\big|^2\big]\\
&\lesssim\varepsilon^4 \tau^2 \sum_{\ell\in \mathcal{T}_{N_0}}\displaystyle\sum_{\ell_2 \ne \pm\ell} \frac{|\delta_{\ell,\ell-\ell_2,\ell_2}-\ell_1^2|^2}{|\delta_{\ell,\ell-\ell_2,\ell_2}|^2}
|\hat{\Phi}_{\ell-\ell_2}|^2|q_{\ell_2}|\\
&\lesssim \varepsilon^4 \tau^2 \sum_{\ell\in \mathcal{T}_{N_0}}\displaystyle\sum_{\ell_2 \ne \pm\ell} |\hat{\Phi}_{\ell-\ell_2}|^2 |\ell_2|^2|q_{\ell_2}|\\
&\lesssim \varepsilon^4 \tau^2 \sum_{\ell \in\mathbb{Z}}|\hat{\Phi}_{\ell}|^2 \sum_{\ell\in\mathbb{Z}}|q_{\ell}||\ell|^2\\
&\lesssim \varepsilon^4 \tau^2 \|\Phi\| \|Q^{\frac{1}{2}}\|_{\mathcal{L}_2^1}.
\end{align*}
Finally, for the term $\mathcal{R}^\varepsilon_3$, similarly to the case of $\mathcal{R}^\varepsilon_1$, we get
\begin{align*}
\mathcal{R}^\varepsilon_3 &\le \varepsilon^4 \tau^2 \mathbb{E}\big[\sum_{\ell\in\mathcal{T}_{N_0}}\big(\sum_{\ell_2\ne \pm \ell} |\ell_2||\hat{\Phi}_{\ell-\ell_2}| |\widehat{v^\varepsilon}_{\ell_2}(t_n)|\big)^2\big] \\
&\lesssim \varepsilon^4 \tau^2 \|\Phi\|^2_{\ell^1} \mathbb{E}\big[\|v^\varepsilon(t_n)\|^2_{H^1}\big]\\
&\lesssim \varepsilon^4 \tau^2.
\end{align*}
This finally gives the estimate
\begin{equation}
    \label{fe1}
\mathbb{E}\big[\|\sum_{k=0}^{n}\mathcal{P}_{N_0}\mathcal{R}^\varepsilon(v^\varepsilon(t_k))\|^2\big]\lesssim\varepsilon^4 \tau^2.
\end{equation}
Using \eqref{fe0} and \eqref{fe1} into \eqref{fe}, we finally achieve
$$
\begin{aligned}
    \mathbb{E}\|e^\varepsilon_{n+1}\|^{2} \lesssim \tau_0^2+\varepsilon^2 \tau \sum_{k=0}^{n}\mathbb{E}\|e^\varepsilon_k\|^2 +\varepsilon^4 \tau^2.
\end{aligned}
$$
The result then follows by an application of the discrete Gronwall inequality and taking the square root both sides.
\begin{comment}
An application of the discrete Gronwall inequality here gives
$$
\begin{aligned}
    \mathbb{E}\|e_{n+1}\|^{2p}_{H^1}&\lesssim \varepsilon^{2p(1-\alpha)}\tau^{2p} {\rm exp}(\varepsilon^{(2p-1)(1-\alpha)}\varepsilon \tau (n+1))\\
    &\lesssim \varepsilon^{2p(1-\alpha)}\tau^{2p} {\rm exp}(\varepsilon^{(2p-1)(1-\alpha)}\varepsilon^{1-\alpha}).
\end{aligned}
$$
This implies
$$
\mathbb{E}\|e_{n+1}\|^{2p}_{H^1}\lesssim \varepsilon^{2p(1-\alpha)}\tau^{2p},
$$
i.e.,
$$
(\mathbb{E}\|e_{n+1}\|^{2p}_{H^1})^{\frac{1}{2p}}\lesssim \varepsilon^{1-\alpha}\tau,
$$
which gives the result. in \eqref{re1}. To obtain \eqref{res2}, for a given $N\in\mathbb{N}$, for any $t_n\in[0,\varepsilon^{-1}]$, any $r<1$ and any $q\in[0,\infty)$, we define the event
$$
\mathcal{A}_N=\{\|v(t_n)-v_n\|_{H^1}>\varepsilon^q \tau^r\},
$$
with $\tau=T/N.$ Applying Markov inequality and \eqref{re1}, we get
$$
\begin{aligned}
\mathbb{P}\big(\mathcal{A}_N\big)&\le \varepsilon^{-2pq}\tau^{-2pr}\displaystyle\max_{1\le n\le N}\mathbb{E}\|v(t_n)-v_n\|^{2p}_{H^1}\\
&\le c  \varepsilon^{-2pq}\tau^{2p(1-r)} \\
&=c \varepsilon^{-2pq} T^{2p(1-r)}N^{2p(r-1)}\\
&\le c  \varepsilon^{-2p(q+r-1)} N^{2p(r-1)},
\end{aligned}
$$
where $c$ is independent on $\varepsilon,T$ and $\tau$. Then, we have
$$
\displaystyle\sum_{N=0}^{\infty}\mathbb{P}\big(\mathcal{A}_N\big)\le  c  \varepsilon^{-2p(q+r-1)}\displaystyle\sum_{N=0}^\infty N^{2p(r-1)}<\infty,
$$
for all $2p>1/|1-r|$. The result \eqref{res2} then follows by an application of Borel-Cantelli lemma. The proof is completed.
\end{comment}
\end{proof}

\begin{remark}
    Theorem \ref{lobg_term_res1} states that the long-term $L^2(\Omega; L^2)$ error for method (SNRLR1) remains bounded and of size $\mathcal{O}(\varepsilon^2 \tau)$ up to time $\mathcal{O}(\varepsilon^{-2})$ provided $\tau$ sufficiently small and $\sigma$ reasonable large. Indeed, a bias of size $\mathcal{O}(\tau_0^\sigma)$ appears also in this case (e.g., see \cite{fe_ka}), ensuring an excellent long-term behaviour for smooth enough solutions. Indeed, in this last case one finally obtains
    $$
    \|u^\varepsilon(t_n)-u^\varepsilon_n\|_{L^2(\Omega; L^2)}\lesssim\varepsilon^2 \tau, \qquad \qquad 0\le n \le \frac{T}{\varepsilon^2 \tau}.
    $$
\end{remark}

\begin{remark}
    If one approximates $\mathcal{W}_n$ as in Remark \ref{app_W}, then, in the estimate of the error $e^\varepsilon_{n+1}$  of Theorem \ref{lobg_term_res1}, one has to take into account the contribution of this kind of approximation, i.e., we end up with the addition of a term
    $$
    \varepsilon^2 \mathbb{E}\big[\|\sum_{k=0}^{n}\int_{t_k}^{t_{k+1}} \big({\rm e}^{-{\bf i}s\Delta} -{\rm e}^{-{\bf i}t_k \Delta} \big)Q^{\frac{1}{2}}\,dW(s)\|^2\big].
    $$
    Using Lemma 2.1 in \cite{an_co}, the Itô isometry and the unitary property of ${\rm e}^{{\bf i}t\Delta}$, for $Q^{\frac{1}{2}}\in\mathcal{L}_2^2$, we get that this term is can be bounded as 
    \begin{align*}
    &\varepsilon^2 \mathbb{E}\big[\|\sum_{k=0}^{n}\int_{t_k}^{t_{k+1}} \big({\rm e}^{-{\bf i}s\Delta} -{\rm e}^{-{\bf i}t_k \Delta} \big)Q^{\frac{1}{2}}\,dW(s)\|^2\big]\\
    &\qquad = \varepsilon^2 
    \mathbb{E}\big[\|\int_0^{t_{n+1}} ({\rm e}^{-{\bf i}s\Delta}-{\rm e}^{-{\bf i}([s/\tau]\tau)\Delta})Q^{\frac{1}{2}}\,dW(s)\|^2\big]\\
    &\qquad = \varepsilon^2 \sum_{k=0}^{n}\int_{t_k}^{t_{k+1}}\|{\big(\rm e}^{{\bf i}(t_k-s)\Delta}-I\big)Q^{\frac{1}{2}}\|^2_{\mathcal{L}_2^0}\,ds\\
    &\qquad \lesssim \varepsilon^2 \sum_{k=0}^{n}\int_{t_k}^{t_{k+1}}\big|t_k-s\big|^2 \|Q^{\frac{1}{2}}\|^2_{\mathcal{L}_2^2}\,ds\\
    &\qquad \lesssim \varepsilon^2 \tau^3 (n+1)\\
    &\qquad \lesssim \tau^2.
    \end{align*}
    This shows that we lost the bound $\mathcal{O}(\varepsilon^2 \tau)$ for the error.
\end{remark}

\section{Nonlinear SSE with cubic nonlinearity}
In this section, we study Equation \eqref{SSE_2} with $F(\xi)=|\xi|^2\xi$, i.e., we consider the stochastic cubic Schrödinger equation
\begin{equation}
\label{cub_eq}
d\psi(t) = {\bf i}\Delta \psi(t) dt-{\bf i}|\psi(t)|^2 \psi(t)   dt+Q^{\frac{1}{2}}dW, \qquad t\in[0,T],\  x \in \mathbb{T},
\end{equation}
under  periodic boundary conditions. Note that we change notation for the solution to distinguish for the previous section. A mild solution for \eqref{cub_eq} reads
\begin{equation}
    \label{mild_cub}
    \psi(t)={\rm e}^{{\bf i}t \Delta }\psi_0-{\bf i}\int_0^t {\rm e}^{{\bf i}(t-s)\Delta} |\psi(s)|^2 \psi(s)\,ds
    +\int_0^t {\rm e}^{{\bf i}(t-s)\Delta} Q^{\frac{1}{2}}\,dW(s).
 \end{equation}
From \cite{cui}, we get the following well-posedness result for \eqref{cub_eq}.
\begin{theorem}
\label{wp_cub}
    Let $\psi_0\in H^{\sigma}$ and $Q^{\frac{1}{2}}\in\mathcal{L}_2^\sigma$, for some $\sigma\in \mathbb{N}^+$. Then, there exists a unique solution to \eqref{cub_eq} for $T>0$ given by \eqref{mild_cub} satisfying
    $$
    \mathbb{E}\big[\sup_{t\in[0,T]}\|\psi(t)\|^p_{H^\sigma}\big]\le C,
    $$
    for any $p\in\mathbb{N}^+$.
\end{theorem}
\subsection{A low-regularity integrator for cubic equation}
%For Equation \eqref{cub_eq}, we will derive directly a non-resonant low-regularity integrator, following the idea presented in the previous section and given also in \cite{fe_ka}. 
Let $\eta(t)={\rm e}^{-{\bf i}t\Delta} \psi(t)$ the twisted variable as in \eqref{twis_def}. It is direct to get (see also e.g. \cite{os_ka}) that $\eta(t)$ solves the equation
$$
\eta(t)=\eta_0-{\bf i}\int_0^t {\rm e}^{-{\bf i}s\Delta}|{\rm e}^{{\bf i}s\Delta} \eta(s)|^2 {\rm e}^{{\bf i}s\Delta} \eta(s) \,ds+\int_0^t {\rm e}^{-{\bf i}s\Delta} Q^{\frac{1}{2}}\,dW(s),
$$
with $\eta_0=u_0$. Then, we get 
\begin{align}
\eta(t_n+\tau)=\eta(t_n)&-{\bf i}\int_0^\tau {\rm e}^{-{\bf i}(t_n+s)\Delta}|{\rm e}^{{\bf i}(t_n+s)\Delta}\eta(t_n+s)|^2 {\rm e}^{{\bf i}(t_n+s)\Delta}\eta(t_n+s)\,ds\notag \\
& +\int_{t_n}^{t_n+\tau} {\rm e}^{-{\bf i}s\Delta} Q^{\frac{1}{2}}\,dW(s)  \label{mild_cub_twis}.
\end{align}
Then, an approximation of $\eta(t_n+\tau)$ can be constructed by the following 
$$
\begin{aligned}
\int_0^\tau {\rm e}^{-{\bf i}(t_n+s)\Delta}|{\rm e}^{{\bf i}(t_n+s)\Delta}\eta(t_n+s)|^2 & {\rm e}^{{\bf i}(t_n+s)\Delta}\eta(t_n+s)\,ds  \\
&\approx \int_0^\tau {\rm e}^{-{\bf i}(t_n+s)\Delta}|{\rm e}^{{\bf i}(t_n+s)\Delta}\eta(t_n)|^2 {\rm e}^{{\bf i}(t_n+s)\Delta}\eta(t_n)\,ds
\end{aligned}
$$
and then approximating the integral in the right hand side of the above last line. Although such discretization relies on same idea given in \cite{os_ka}, we here briefly report the principle for the sake of completeness. Expanding $\eta(t_n)$ through its Fourier series, we have
\begin{align*}
&|{\rm e}^{{\bf i}(t_n+s)\Delta}\eta(t_n)|^2 {\rm e}^{{\bf i}(t_n+s)\Delta}\eta(t_n)\\
&\qquad =\bigg(\frac{1}{\sqrt{2\pi}}\bigg)^{3}\sum_{\ell\in\mathbb{Z}}{\rm e}^{{\bf i}\ell x}
\sum_{\substack{\ell_1,\ell_2,\ell_3\in\mathbb{Z}\\ \ell=-\ell_1+\ell_2+\ell_3}}{\rm e}^{{\bf i}(\ell_1^2-\ell_2^2-\ell_3^2)(t_n+s)}\overline{\hat{\eta}_{\ell_1}}(t_n)
\hat{\eta}_{\ell_2}(t_n)\hat{\eta}_{\ell_3}(t_n)
\end{align*}
so that
\begin{align}
    &\int_0^\tau {\rm e}^{-{\bf i}(t_n+s)\Delta}|{\rm e}^{{\bf i}(t_n+s)\Delta}\eta(t_n)|^2 {\rm e}^{{\bf i}(t_n+s)\Delta}\eta(t_n)\,ds\notag \\
    &\qquad = \bigg(\frac{1}{\sqrt{2\pi}}\bigg)^{3} \sum_{\ell\in\mathbb{Z}} {\rm e}^{{\bf i}\ell x} \sum_{\substack{\ell_1,\ell_2,\ell_3\in\mathbb{Z}\\ \ell=-\ell_1+\ell_2+\ell_3}}
    {\rm e}^{{\bf i}(\ell^2-\ell_1^2-\ell_2^2-\ell_3^2) t_n}\overline{\hat{\eta}_{\ell_1}}(t_n)
\hat{\eta}_{\ell_2}(t_n)\hat{\eta}_{\ell_3}(t_n)\notag \\
&\hspace{4cm}\times \int_0^\tau {\rm e}^{{\bf i}(\ell^2-\ell_1^2-\ell_2^2-\ell_3^2) s}\,ds. \label{four_met}
\end{align}
Using the key-relation
$$
\ell^2-\ell_1^2-\ell_2^2-\ell_3^2=2\ell_1^2-2\ell_1(\ell_2+\ell_3)+2\ell_2\ell_3,
$$
we take the following approximation
\begin{align*}
    \int_0^\tau {\rm e}^{{\bf i}(\ell^2-\ell_1^2-\ell_2^2-\ell_3^2) s}\,ds &= \int_0^\tau {\rm e}^{{\bf i}( 2\ell_1^2-2\ell_1(\ell_2+\ell_3)+2\ell_2\ell_3) s}\,ds\\
    &\approx  \int_0^\tau {\rm e}^{2{\bf i}\ell_1^2 s}\,ds\\
    &= \tau \varphi_1(2{\bf i}\ell_1^2 \tau).
\end{align*}
Hence, after returning back to the physical space, we end up with the following approximation
\begin{align*}
    &\int_0^\tau {\rm e}^{-{\bf i}(t_n+s)\Delta}|{\rm e}^{{\bf i}(t_n+s)\Delta}\eta(t_n)|^2 {\rm e}^{{\bf i}(t_n+s)\Delta}\eta(t_n)\,ds\\
    &\qquad \approx \tau \bigg(\frac{1}{\sqrt{2\pi}}\bigg)^{3} \sum_{\ell\in\mathbb{Z}} {\rm e}^{{\bf i}\ell x} \sum_{\substack{\ell_1,\ell_2,\ell_3\in\mathbb{Z}\\ \ell=-\ell_1+\ell_2+\ell_3}}
    {\rm e}^{{\bf i}(\ell^2-\ell_1^2-\ell_2^2-\ell_3^2) t_n}\overline{\hat{\eta}_{\ell_1}}(t_n)
\hat{\eta}_{\ell_2}(t_n)\hat{\eta}_{\ell_3}(t_n)  \varphi_1(2{\bf i}\ell_1^2 \tau)\\
&\qquad =\tau {\rm e}^{-{\bf i}t_n\Delta}\big(({\rm e}^{{\bf i}t_n\Delta}\eta(t_n))^2\varphi_1(-2{\bf i}\tau\Delta){\rm e}^{{-\bf i}t_n\Delta}\overline{\eta(t_n)}\big).
\end{align*}
Then, we obtain the following approximation for the twisted variable $\eta$
$$
\eta_{n+1}=\eta_n-{\bf i}\tau {\rm e}^{-{\bf i}t_n\Delta}\big(({\rm e}^{{\bf i}t_n\Delta}\eta_n)^2\varphi_1(-2{\bf i}\tau\Delta){\rm e}^{{-\bf i}t_n\Delta}\overline{\eta_n}\big)
+\int_{t_n}^{t_n+\tau} {\rm e}^{-{\bf i}s\Delta} Q^{\frac{1}{2}}\,dW(s).
$$
%The following approximation will be here adopted
%\begin{align*}
 %  \int_0^\tau {\rm e}^{-{\bf i}(t_n+s)\Delta}|{\rm e}^{{\bf i}(t_n+s)\Delta}v(t_n)|^2 & {\rm e}^{{\bf i}(t_n+s)\Delta}v(t_n)\,ds  \\ 
  % & \approx \tau\big[{\rm e}^{-{\bf i}t_n\Delta}\big(({\rm e}^{{\bf i}t_n\Delta} v_n)^2 \varphi_1(-2{\bf i}\tau\Delta){\rm e}^{-{\bf i}t_n\Delta} \overline{v_n}\big)
   %+2\widehat{\big(g(v_n)\big)}_0 v_n-h(v_n)\big],
%\end{align*}
%where $g(u)=u(1-\varphi_1(-2{\bf i}\tau \Delta))\overline{u}$, $\widehat{(h(u))}_\ell=(1-\varphi_1(2{\bf i}\ell^2 \tau))\overline{\hat{u}_\ell}\hat{u}_\ell \hat{u}_\ell$.
%and
%$$
%\varphi_1(z)=\frac{{\rm e}^z-1}{z}, \qquad \qquad \big(\widehat{\varphi_1({\bf i}t\Delta )}\big)_\ell
%=\begin{cases}
  %  \varphi_1(-{\bf i}t \ell^2), & \ell \ne 0,\\
   % 1, & \ell=0.
%\end{cases}
%$$

%Then, we get the following approximation for the twisted variable $v$
%\begin{align*}
%v_{n+1}=v_n&-{\bf i} \tau\big[{\rm e}^{-{\bf i}t_n\Delta}\big(({\rm e}^{{\bf i}t_n\Delta} v_n)^2 \varphi_1(-2{\bf i}\tau\Delta){\rm e}^{-{\bf i}t_n\Delta} \overline{v_n}\big)
 %  +2\widehat{\big(g(v_n)\big)}_0 v_n-h(v_n)\big]\\
  % &+\int_{t_n}^{t_n+\tau} {\rm e}^{-{\bf i}s\Delta}Q^{\frac{1}{2}}\,dW(s).
  % \end{align*}
Returning back to the original variable $\psi$, we get the following stochastic low-regularity integrator (SLR2) for \eqref{cub_eq}
   % \begin{align}
   % u_{n+1}&={\rm e}^{{\bf i}\tau\Delta}\big[u_n-{\bf i}\tau(u_n)^2\varphi_1(-2{\bf i}\tau\Delta)\overline{u_n}\big]+\mathcal{W}_n\notag\\
   % &\hspace{2cm}-2{\bf i}\tau \widehat{\big(g(u_n)\big)}_0 {\rm e}^{{\bf i}\tau\Delta} u_n+{\bf i}\tau {\rm e}^{{\bf i}\tau\Delta} h(u_n)    \label{SLNRI} .
    %    \end{align}

         \begin{align}
         \label{SLR2}
    \psi_{n+1}&={\rm e}^{{\bf i}\tau\Delta}\big[\psi_n-{\bf i}\tau\psi_n^2\varphi_1(-2{\bf i}\tau\Delta)\overline{\psi_n}\big]+\mathcal{W}_n.
    %&\hspace{2cm}-2{\bf i}\tau \widehat{\big(g(u_n)\big)}_0 {\rm e}^{{\bf i}\tau\Delta} u_n+{\bf i}\tau {\rm e}^{{\bf i}\tau\Delta} h(u_n)    \label{SLR2} .
        \end{align}
\begin{comment}
Let $\varphi^t$ a flow map such that $\varphi^t(v_0)=v(t)$ and let
$$
\begin{aligned}
\Phi_t^\tau(w)=w&-{\bf i} \tau\big[{\rm e}^{-{\bf i}t_n\Delta}\big(({\rm e}^{{\bf i}t_n\Delta} w)^2 \varphi_1(-2{\bf i}\tau\Delta){\rm e}^{-{\bf i}t_n\Delta} \overline{w}\big)
   +2\widehat{\big(g(w)\big)}_0 w-h(w)\big]\\
   &+\int_{t_n}^{t_n+\tau} {\rm e}^{-{\bf i}s\Delta}Q^{\frac{1}{2}}\,dW(s)
   \end{aligned}
$$
such that one has $v_{n+1}=\Phi_{t_n}^\tau(v_n)$.
\end{comment}

As remarked in \cite{cui}, in contrast to what done for the method (SLR1) in the linear case, it is not possible to derive a strong convergence result for several explicit exponential methods applied to the cubic equation \eqref{cub_eq}, due to the lack of exponential integrability conditions needed to obtain such a convergence result. This is also an issue for the derived integrator (SLR2). However, with similar ideas than those developed in \cite{deb0}, we will be able to get $\mathbb{P}-{\rm a.s.}$ convergence result and order, using an auxiliary truncated equation for \eqref{cub_eq}.
\subsection{Pathwise convergence}
In order to investigate the pathwise convergence and order of the method (SLR2), applied to the stochastic cubic Schrödinger equation \eqref{cub_eq}, we introduce the follow cut-off function $\theta\in C^\infty([0,\infty))$, such that $\theta(x)=1$ for $x\in[0,1]$ and $\theta(x)=0$ for $x\ge 2$. Then, for $\sigma\in\mathbb{N}^+$ and $R>0$, we define $\theta_R(\cdot)=\theta(\|\cdot\|_{H^\sigma}/R)$ and we consider the truncated equation
\begin{align}
\psi^R(t)={\rm e}^{{\bf i}t\Delta}\psi_0&-{\bf i}\int_0^t {\rm e}^{{\bf i}(t-s)\Delta} \theta_R(\psi^R(s))|\psi^R(s)|^2 \psi^R(s)\,ds\notag \\
&+\int_0^t {\rm e}^{{\bf i}(t-s)\Delta} Q^{\frac{1}{2}}\,dW(s) \label{trunc_cub}.
\end{align}
\begin{remark}
    It is classical (see, e.g., \cite{deb0}) that Theorem \ref{wp_cub} applies also to \eqref{trunc_cub} since it has globally Lipschitz nonlinearity. In particular, one has that if $\psi_0\in H^\sigma$ and $Q^{\frac{1}{2}}\in\mathcal{L}_2^\sigma$, for $\sigma\in\mathbb{N}^+$, then $\psi^R \in L^{2p}(\Omega; L^\infty([0,T];H^\sigma))$, for any $p\in\mathbb{N}^+$. Moreover, if one assumes $\psi_0\in H^{\sigma+\gamma}$ and $Q^{\frac{1}{2}}\in\mathcal{L}_2^{\sigma+\gamma}$, with $\gamma\in[0,1]$, then $u^R \in L^{2p}(\Omega; L^\infty([0,T];H^{\sigma+\gamma}))$ and
\begin{equation}
    \label{time_reg_trunc}
    \mathbb{E}\big[\|\psi^R(t_2)-\psi^R(t_1)\|^{2p}_{H^\sigma}\big]\le K_{R,p}|t_2-t_1|^{p\gamma}\mathbb{E}\big[\displaystyle\sup_{t\in[t_1,t_2]}\|\psi^R(t)\|_{H^{\sigma+\gamma}}^{2p}\big],
\end{equation}
for any $p\in\mathbb{N}^+$ and any $t_1\le t_2$. 
\end{remark}
The method (SLR2), applied to the truncated equation \eqref{trunc_cub}, reads
 \begin{align}
 \label{trunc_met}
    \psi^R_{n+1}&={\rm e}^{{\bf i}\tau\Delta}\big[\psi_n^R-{\bf i}\tau \theta_R(\psi^R_n) (\psi_n^R)^2\varphi_1(-2{\bf i}\tau\Delta)\overline{\psi_n^R}+\big]+\mathcal{W}_n\notag
    %\\
    %&&\hspace{2cm}-2{\bf i}\theta_R(u^R_n)\tau \widehat{\big(g(u_n^R)\big)}_0 {\rm e}^{{\bf i}\tau\Delta} u_n^R+{\bf i}\theta_R(u^R_n)\tau {\rm e}^{{\bf i}\tau\Delta} h(u_n^R) \label{trunc_met}.
      \end{align}
      We give the following lemma.
      \begin{lemma}
         Let us assume $\psi_0\in H^\sigma, Q^{\frac{1}{2}}\in\mathcal{L}_2^\sigma, \ \sigma\in\mathbb{N}^+$. Then, the numerical integrator (SLR2), applied to \eqref{trunc_cub}, satisfies
         $$
         \sup_{N\in\mathbb{N}^+}\mathbb{E}\big[\max_{1\le n\le N}\|\psi^R_n\|_{H^\sigma}^{2p}\big]\le C_{R,p},
         $$
         for any $p\in\mathbb{N}^+$.
      \end{lemma}
      \begin{proof}
          Similarly to the proof of Lemma \ref{prop_stab_lin}, by recursion, we get
          $$
          \psi_n^R={\rm e}^{{\bf i}t_n\Delta} \psi_0-{\bf i}\tau \sum_{j=1}^{n}{\rm e}^{{\bf i}j\tau \Delta}\theta_R(\psi_{n-j}^R) (\psi_{n-j}^R)^2 \varphi_1(-2{\bf i}\tau\Delta)\overline{\psi_{n-j}^R}+\sum_{j=1}^{n}{\rm e}^{{\bf i}(j-1)\tau\Delta}\mathcal{W}_{n-j}.
          $$
          So, we have only to bound 
          $$
          \mathbb{E}\big[\|{\bf i}\tau\sum_{j=1}^{n}{\rm e}^{{\bf i}j\tau \Delta}\theta_R(\psi_{n-j}^R) (\psi_{n-j}^R)^2 \varphi_1(-2{\bf i}\tau\Delta)\overline{\psi_{n-j}^R}\|^{2p}_{H^\sigma}\big].
          $$
          It is direct to see (see, e.g., \cite{os_ka}) that the algebra property of $H^\sigma$ gives
          $$
          \| (\psi_{n-j}^R)^2 \varphi_1(-2{\bf i}\tau\Delta)\overline{\psi_{n-j}^R}\|_{H^\sigma}\lesssim \|\psi_{n-j}^R\|^3_{H^\sigma}.
          $$
          Then, using this, the discrete Hölder inequality and the definition of the function $\theta_R$, we get
          \begin{align*}
             & \mathbb{E}\big[\|{\bf i}\tau\sum_{j=1}^{n}{\rm e}^{{\bf i}j\tau \Delta}\theta_R(\psi_{n-j}^R) (\psi_{n-j}^R)^2 \varphi_1(-2{\bf i}\tau\Delta)\overline{\psi_{n-j}^R}\|^{2p}_{H^\sigma}\big]\\
             &\qquad \le \tau^{2p} n^{2p-1} \sum_{j=1}^{n}\mathbb{E}\big[\theta_R(\psi_{n-j}^R)^{2p}\|(\psi_{n-j}^R)^2 \varphi_1(-2{\bf i}\tau\Delta)\overline{\psi_{n-j}^R}\|^{2p}_{H^\sigma}\big]\\
             &\qquad \lesssim \tau^{2p} n^{2p-1} \sum_{j=1}^{n}\mathbb{E}\big[\theta_R(\psi_{n-j}^R)^{2p}\|\psi_{n-j}^R\|^{6p}_{H^\sigma}\big]\\
             &\qquad \lesssim \tau^{2p} n^{2p-1} \sum_{j=1}^{n} (2R)^{6p}\\
             &\qquad \lesssim T^{2p}.
          \end{align*}
          The rest of the proof works in the same way of the proof of Lemma \ref{prop_stab_lin}.
      \end{proof}

The following lemma provides a first-order strong convergence result for the numerical solution $\psi^R_n$ given by method (SLR2), applied to \eqref{trunc_cub}.
\begin{lemma}
\label{strong_lem_trunc}
    Let us assume $\psi_0\in H^{\sigma+\gamma}$, $Q^{\frac{1}{2}}\in \mathcal{L}_2^{\sigma+\gamma}$, for $\sigma\in\mathbb{N}^+$ and $\gamma\in(0,1]$. Then, for any $p\in\mathbb{N}^+$, we have
    $$
    \big(\mathbb{E}\big[\max_{1\le n\le N}\|\psi^R(t_n)-\psi^R_n\|^{2p}_{H^\sigma}\big]\big)^{\frac{1}{2p}}\lesssim c_R \tau^{\gamma}
    $$
\end{lemma}
\begin{proof}
   Also here, it is convenient to introduce the twisted variable $\eta^R(t)={\rm e}^{-{\bf i}t \Delta}\psi^R(t)$. It is possible to directly check that the following holds
    \begin{align*}
        \eta^R(t_n)-\eta^R_n&=-{\bf i}\sum_{j=0}^{n-1}\int_0^\tau {\rm e}^{-{\bf i}(t_j+s)\Delta}\theta_R(\eta^R(t_j+s))|{\rm e}^{{\bf i}(t_j+s)\Delta}\eta^R(t_j+s)|^2 {\rm e}^{{\bf i}(t_j+s)\Delta}\eta^R(t_j+s)\,ds\\
        &\hspace{2cm}- \tau \theta_R(\eta^R_j){\rm e}^{-{\bf i}t_j\Delta}\big(({\rm e}^{{\bf i}t_j\Delta} \eta^R_j)^2 \varphi_1(-2{\bf i}\tau\Delta){\rm e}^{-{\bf i}t_j\Delta} \overline{\eta_j^R}\big)\\ 
   %+2\widehat{\big(g(\_j^R)\big)}_0 v^R_j-h(v_j^R)\big]\\
   &=-{\bf i}\sum_{j=0}^{n-1}\mathcal{I}_{1,j}+\tau \mathcal{I}_{2,j},
    \end{align*}
    where 
    \begin{align*}
    \mathcal{I}_{1,j}&:=\int_0^\tau {\rm e}^{-{\bf i}(t_j+s)\Delta}\theta_R(\eta^R(t_j+s))|{\rm e}^{{\bf i}(t_j+s)\Delta}\eta^R(t_j+s)|^2 {\rm e}^{{\bf i}(t_j+s)\Delta}\eta^R(t_j+s)\,ds\\
    &\hspace{2cm}- \tau \theta_R(\eta^R(t_j)){\rm e}^{-{\bf i}t_j\Delta}\big(({\rm e}^{{\bf i}t_j\Delta} \eta^R(t_j))^2 \varphi_1(-2{\bf i}\tau\Delta){\rm e}^{-{\bf i}t_j\Delta} \overline{\eta^R(t_j)}\big)
   %+2\widehat{\big(g(v^R(t_j))\big)}_0 v^R(t_j)-h(v^R(t_j))\big]
     \end{align*}
    and
    \begin{align*}
    \mathcal{I}_{2,j}&:=\theta_R(\eta^R(t_j)) {\rm e}^{-{\bf i}t_j\Delta}\big(({\rm e}^{{\bf i}t_j\Delta} \eta^R(t_j))^2 \varphi_1(-2{\bf i}\tau\Delta){\rm e}^{-{\bf i}t_j\Delta} \overline{\eta^R(t_j)}\big)\\
   %+2\widehat{\big(g(v^R(t_j))\big)}_0 v^R(t_j)-h(v^R(t_j))\big]\\
   &\hspace{2cm} - \theta_R(\eta^R_j){\rm e}^{-{\bf i}t_j\Delta}\big(({\rm e}^{{\bf i}t_j\Delta} \eta^R_j)^2 \varphi_1(-2{\bf i}\tau\Delta){\rm e}^{-{\bf i}t_j\Delta} \overline{\eta_j^R}\big)
%   +2\widehat{\big(g(v_j^R)\big)}_0 v^R_j-h(v_j^R)\big.
     \end{align*}
    Hence, we get
    \begin{align*}
        \mathbb{E}\big[\max_{1\le n\le N}\|\eta^R(t_n)-\eta^R_n\|^{2p}_{H^\sigma}\big]& = \mathbb{E}\big[\max_{1\le n\le N}\|\sum_{j=0}^{n-1}\mathcal{I}_{1,j}+\tau \mathcal{I}_{2,j}\|^{2p}_{H^\sigma}\big]\\
    &\lesssim \mathbb{E}\big[\max_{1\le n\le N}\|\sum_{j=0}^{n-1}\mathcal{I}_{1,j}\|^{2p}_{H^\sigma}\big]+\tau^{2p}\mathbb{E}\big[\max_{1\le n\le N}\|\sum_{j=0}^{n-1}\mathcal{I}_{2,j}\|^{2p}_{H^\sigma}\big].
    \end{align*}
    We start by estimating $\mathbb{E}\big[\max_{1\le n\le N}\|\sum_{j=0}^{n-1}\mathcal{I}_{1,j}\|^{2p}_{H^\sigma}\big]$. Denoting $$F_{t}(u)=\theta_R(u)|{\rm e}^{{\bf i}t\Delta}u|^2 {\rm e}^{{\bf i}t\Delta}u,$$ we have
    $$
    \mathcal{I}_{1,j}=\mathcal{I}^a_{1,j}+\mathcal{I}^b_{1,j},
    $$
    where
    \begin{align*}
    &\mathcal{I}^a_{1,j}=\int_0^\tau {\rm e}^{-{\bf i}(t_j+s)\Delta}\big[ F_{t_j+s}(\eta^R(t_j+s))-F_{t_j+s}(\eta^R(t_j))\big]\,ds,\\
    &\mathcal{I}^b_{1,j}=\int_0^\tau {\rm e}^{-{\bf i}(t_j+s)\Delta}F_{t_j+s}(\eta^R(t_j))\,ds\\
    &\hspace{2cm}-\tau\theta_R(\eta^R(t_j)){\rm e}^{-{\bf i}t_j\Delta}\big(({\rm e}^{{\bf i}t_j\Delta} \eta^R(t_j))^2 \varphi_1(-2{\bf i}\tau\Delta){\rm e}^{-{\bf i}t_j\Delta} \overline{\eta^R(t_j)}\big)
   %+2\widehat{\big(g(v^R(t_j))\big)}_0 v^R(t_j)-h(v^R(t_j))\big].
        \end{align*}
    Moreover, using the fact that $F_t$ is bounded with bounded first and second derivatives (see, e.g., \cite{deb0}), by Taylor formula, we get
    \begin{align*}
        \mathcal{I}^a_{1,j}&=\int_0^\tau {\rm e}^{-{\bf i}(t_j+s)\Delta} DF_{t_j+s}(\eta^R(t_j)) \int_{t_j}^{t_j+s}{\rm e}^{-{\bf i}r\Delta}Q^{\frac{1}{2}}\,dW(r)\,ds\\
        &\ -{\bf i}\int_0^\tau {\rm e}^{-{\bf i}(t_j+s)\Delta} F'_{t_j+s}(\eta^R(t_j)) \int_0^s {\rm e}^{-{\bf i}(t_j+r)\Delta} F_{t_j+r}(\eta^R(t_j+r))\,dr\,ds\\
        &\ +\int_0^\tau {\rm e}^{-{\bf i}(t_j+s)\Delta} \int_0^1 D^2 F_{t_j+s}(\mu \eta^R(t_j+s)+(1-\mu)\eta^R(t_j))\\
        &\hspace{5cm}(\eta^R(t_j+s)-\eta^R(t_j),\eta^R(t_j+s)-\eta^R(t_j))\,d\mu \,ds\\
        &=\mathcal{I}_{1,j}^{a,1}+\mathcal{I}_{1,j}^{a,2}+\mathcal{I}_{1,j}^{a,3}.
    \end{align*}
    Similarly to the proof of Theorem \ref{strong_conv_lin}, the stochastic Fubini Theorem and BDG inequality here yield
    \begin{align*}
       &\mathbb{E}\big[\max_{1\le n\le N}\|\sum_{j=0}^{n-1} \mathcal{I}_{1,j}^{a,1}\|^{2p}_{H^\sigma}\big]= \mathbb{E}\big[\max_{1\le n\le N}\|\sum_{j=0}^{n-1}\int_0^\tau {\rm e}^{-{\bf i}(t_j+s)\Delta} DF_{t_j+s}(\eta^R(t_j)) \int_{t_j}^{t_j+s}{\rm e}^{-{\bf i}r\Delta}Q^{\frac{1}{2}}\,dW(r)\,ds \|^{2p}_{H^\sigma}\big]\\
       &\qquad =\mathbb{E}\big[\max_{1\le n\le N}\|\sum_{j=0}^{n-1}\int_{t_j}^{t_{j+1}}\int_{r}^{t_{j+1}} {\rm e}^{-{\bf i}s\Delta} DF_{s}(\eta^R(t_j)) {\rm e}^{-{\bf i}r\Delta}Q^{\frac{1}{2}}\,ds\,dW(r) \|^{2p}_{H^\sigma}\big]\\
       &\qquad =\mathbb{E}\big[\max_{1\le n\le N}\|\int_{0}^{t_{n}}\int_{r}^{[\frac{r}{\tau}+1]\tau} {\rm e}^{-{\bf i}s\Delta} DF_{s}(\eta^R([r/\tau]\tau)) {\rm e}^{-{\bf i}r\Delta}Q^{\frac{1}{2}}\,ds\,dW(r) \|^{2p}_{H^\sigma}\big]\\
       &\qquad \le \mathbb{E}\big[\sup_{0\le t\le T}\|\int_{0}^{t}\int_{r}^{[\frac{r}{\tau}+1]\tau} {\rm e}^{-{\bf i}s\Delta} DF_{s}(\eta^R([r/\tau]\tau)) {\rm e}^{-{\bf i}r\Delta}Q^{\frac{1}{2}}\,ds\,dW(r) \|^{2p}_{H^\sigma}\big]\\
       &\qquad \lesssim \mathbb{E}\big[\big(\int_0^T \|\int_r^{[\frac{r}{\tau}+1]\tau}{\rm e}^{-{\bf i}s\Delta} DF_{s}(\eta^R([r/\tau]\tau)) {\rm e}^{-{\bf i}r\Delta} Q^{\frac{1}{2}}\,ds \|_{\mathcal{L}_2^\sigma}^2\,dr\big)^p\big]\\
       &\qquad \le \|Q^{\frac{1}{2}}\|_{\mathcal{L}_2^\sigma}^{2p} \mathbb{E}\big[\big(\sum_{j=0}^{N-1} \int_{t_j}^{t_{j+1}}|t_{j+1}-r|^2\,dr\big)^p\big]\\
       &\qquad \lesssim \tau^{2p}.
    \end{align*}
    Moreover, applying twice the Hölder inequality and exploiting the boundedness property of the cut-off function and of $DF$, we also see
    \begin{align*}
        \max_{1\le n\le N}\|\sum_{j=0}^{n-1} \mathcal{I}_{1,j}^{a,2}\|^{2p}_{H^\sigma}
        &= \max_{1\le n\le N}\|\sum_{j=0}^{n-1}\int_0^\tau {\rm e}^{-{\bf i}(t_j+s)\Delta} DF_{t_j+s}(\eta^R(t_j))\\
        &\hspace{3cm}\times \int_0^s {\rm e}^{-{\bf i}(t_j+r)\Delta} F_{t_j+r}(\eta^R(t_j+r))\,dr\,ds\|^{2p}_{H^\sigma}\\
        & \le N^{2p-1}\mathbb{E}\big[\sum_{j=0}^{N-1} \big(\int_0^\tau \int_0^s  \|F_{t_j+r}(\eta^R(t_j+r))\|_{H^\sigma}\,dr\,ds\big)^{2p}\big]\\
        & \lesssim N^{2p-1}\mathbb{E}\big[\sum_{j=0}^{N-1} \big(\int_0^\tau \int_0^s 8 R^{3}\,dr\,ds\big)^{2p}\big]\\
        &\lesssim 8 R^3 N^{2p}\tau^{4p}\\
        &\lesssim 8 R^3 \tau^{2p}
    \end{align*}
    and then
    $$
    \mathbb{E}\big[\max_{1\le n\le N}\|\sum_{j=0}^{n-1} \mathcal{I}_{1,j}^{a,2}\|^{2p}_{H^\sigma}\big]\lesssim \tau^{2p}.
    $$
    For the term $\mathcal{I}_{1,j}^{a,3}$, applying several time the discrete and continuous Hölder inequalities, due to the uniform boundedness of $D^2 F_t$ and thanks to \eqref{time_reg_trunc}, we directly get
    \begin{align*}
    \mathbb{E}\big[\max_{1\le n\le N}\|\sum_{j=0}^{n-1} \mathcal{I}_{1,j}^{a,3}\|^{2p}_{H^\sigma}\big]& \lesssim \sum_{j=0}^{N-1} \int_0^\tau \int_0^1 \mathbb{E}\big[\|\eta^R(t_j+s)-\eta^R(t_j)\|^{4p}_{H^\sigma}\big]\,d\mu\,ds\\
    &\lesssim \tau^{2p\gamma}.
    \end{align*}
    Then, as long as $\psi_0\in H^{\sigma+\gamma}$ and $Q^{\frac{1}{2}}\in \mathcal{L}_2^{\sigma+\gamma}$, we get
    $$
    \mathbb{E}\big[\max_{1\le n\le N} \|\sum_{j=0}^{n-1}\mathcal{I}_{1,j}^a\|^{2p}_{H^\sigma}\big]\lesssim \tau^{2p\gamma}, \qquad \gamma\in(0,1].
    $$
    For the term $\mathcal{I}_{1,j}^b$, by the construction of the method, using Fourier expansion, one has
    %we may exploit Fourier techniques as done, e.g., in \cite{fe_ka,os_ka} and also here in previous section. We indeed get
    \begin{align*}
    \mathcal{I}_{1,j}^b&=\theta_R(\eta^R(t_j))\sum_{\ell \in \mathbb{Z}}{\rm e}^{{\bf i}\ell x}\sum_{\substack{\ell_1,\ell_2,\ell_3 \in \mathbb{Z}\\ \ell=-\ell_1+\ell_2+\ell_3\\ }}{\rm e}^{{\bf i}t_j (\ell^2+\ell_1^2-\ell_2^2-\ell_3^2)}\big(\widehat{\eta^R}_{\ell_1}(t_j)\big)^*\widehat{\eta^R}_{{\ell_2}}(t_j)\widehat{\eta^R}_{{\ell_3}}(t_j)\\
    &\hspace{5cm}\times \int_0^\tau {\rm e}^{2{\bf i}s \ell_1}({\rm e}^{{\bf i}\beta}-1)\,ds,
        \end{align*}
    with $\beta=s(-2\ell_1(\ell_2+\ell_3)+2\ell_2\ell_3)$ and being $(\cdot)^*=\overline{\cdot}$. This term can be bounded by using an interpolation argument and the algebra property of $H^\sigma$. Such as estimation is given in \cite{os_ka} and we skip it not to report tedious computations.
    Thus, we obtain
    $$
    \|\mathcal{I}_{1,j}^b\|_{H^\sigma}\lesssim \tau^{1+\gamma}\|\eta^R(t_j)\|^3_{H^{\sigma+\gamma}}.
    $$
    This gives
    \begin{align*}
    \mathbb{E}\big[\max_{1\le n\le N} \|\sum_{j=0}^{n-1}\mathcal{I}_{1,j}^b\|^{2p}_{H^\sigma}\big]&\le N^{2p-1}
    \mathbb{E}\big[\sum_{j=0}^{N-1}\|\mathcal{I}_{1,j}^b\|^{2p}_{H^\sigma}\big]\\
    &\lesssim \tau^{1+2p\gamma}\sum_{j=0}^{N-1}\|\eta^R(t_j)\|^{2p}_{H^{\sigma+\gamma}}\\
    &\lesssim \tau^{2p\gamma}.
     \end{align*}
    Then, we finally obtain
    $$
    \mathbb{E}\big[\max_{1\le n\le N} \|\sum_{j=0}^{n-1}\mathcal{I}_{1,j}\|^{2p}_{H^\sigma}\big]\lesssim \tau^{2p\gamma},
    $$
    provided $\psi_0\in H^{\sigma+\gamma}$ and $Q^{\frac{1}{2}}\in\mathcal{L}_2^{\sigma+\gamma}$, with $\gamma\in(0,1]$.

    The term $\mathcal{I}_{2,j}$ is bounded in the following manner. Following similar arguments of \cite{fe_ka,os_ka} and of those presented in the previous section, by the unitary property of the operator ${\rm e}^{{\bf i}t\Delta}$ and by using the $H^\sigma$ is an algebra for $\sigma\in\mathbb{N}^+$, we have that the map
    $\tilde{F}:H^\sigma\to H^\sigma$ defined by
    $$
    \tilde{F}(u)={\rm e}^{-{\bf i}t_j\Delta}\big(({\rm e}^{{\bf i}t_j\Delta} u)^2 \varphi_1(-2{\bf i}\tau\Delta){\rm e}^{-{\bf i}t_j\Delta} \overline{u}\big)
   %+2\widehat{\big(g(u)\big)}_0 u-h(u)
    $$
    is locally Lipschitz. Then, the map $\mathcal{F}_R:=\theta_R(u)\tilde{F}(u):H^\sigma\to H^\sigma$ is a globally Lipschitz map. This then yields
    \begin{align*}
    \|\mathcal{I}_{2,j}\|_{H^\sigma} &= \|\mathcal{F}_R(\eta^R(t_j))-\mathcal{F}_R(\eta^R_j)\|_{H^\sigma}\\
    &\lesssim \|\eta^R(t_j)-\eta^R_j\|_{H^\sigma}.
    \end{align*}
    Then, we get
        \begin{align*}
            \mathbb{E}\big[\max_{1\le n\le N}\|\sum_{j=0}^{n-1}\mathcal{I}_{2,j}\|^{2p}_{H^\sigma}\big]&\lesssim N^{2p-1}\sum_{j=0}^{N-1}\mathbb{E}\big[\|\eta^R(t_j)-\eta^R_j\|^{2p}_{H^\sigma}\big].
        \end{align*}
    Then, we finally get the following estimate
    \begin{align*}
    \mathbb{E}\big[\max_{1\le n\le N}\|\eta^R(t_n)-\eta_j^R\|^{2p}_{H^\sigma}\big]&\lesssim \tau^{2p\gamma}+\tau^{2p}N^{2p-1}\sum_{j=0}^{N-1}\mathbb{E}\big[\|\eta^R(t_j)-\eta^R_j\|^{2p}_{H^\sigma}\big]\\
    &\lesssim \tau^{2p\gamma}+\tau \sum_{j=0}^{N-1}\mathbb{E}\big[\max_{0\le k\le j} \|\eta^R(t_k)-\eta^R_k\|^{2p}_{H^\sigma}\big].
     \end{align*}
    An application of the discrete Gronwall inequality then completes the proof.
\end{proof}

 The following lemma provides a pathwise error estimate for the truncated numerical solution $\psi^R_n$, given by method (SLR2) applied to \eqref{trunc_cub}.
    \begin{lemma}
    \label{path_trunc_lem}
        Under the assumptions of Lemma \ref{strong_lem_trunc}, for any $\delta<\gamma$, there exists a random variable $K_\delta(\omega)$ such that one has
        $$
        \max_{1\le n\le N} \|\psi^R(t_n)-\psi^R_n\|_{H^\sigma}\lesssim K_\delta(\omega) \tau^{\delta},\qquad \mathbb{P}-{\rm a.s.}
        $$
    \end{lemma}

    \begin{proof}
        The proof exploits Borel-Cantelli lemma, several $L^{2p}$-estimate presented in the proof of Lemma \ref{strong_lem_trunc} and similar arguments of the proof of Proposition 4.5 in \cite{deb0}. Here, we highlight only the differences. We give the result for the twisted variable $v^R$. Following the decomposition given in the proof of Lemma \ref{strong_lem_trunc} and the pathwise estimate given there for the term $\mathcal{I}_{2,j}$, we have
        \begin{align*}
        \max_{1\le n\le N}\|\eta^R(t_n)-\eta^R_n\|_{H^\sigma}&\le \max_{1\le n\le N}\|\sum_{j=0}^{n-1}\mathcal{I}_{1,j}^{a,1}+\mathcal{I}_{1,j}^{a,2}+\mathcal{I}_{1,j}^{a,3}+\mathcal{I}_{1,j}^{b}\|_{H^\sigma}+\max_{1\le n\le N}\|\sum_{j=0}^{n-1} \mathcal{I}_{2,j}\|_{H^\sigma}\\
        &\le \max_{1\le n\le N}\|\sum_{j=0}^{n-1}\mathcal{I}_{1,j}^{a,1}+\mathcal{I}_{1,j}^{a,2}+\mathcal{I}_{1,j}^{a,3}+\mathcal{I}_{1,j}^{b}\|_{H^\sigma}
        +\tau\sum_{j=0}^{N-1}\max_{0\le k\le j}\|\eta^R(t_k)-\eta^R_k\|_{H^\sigma}.
        \end{align*}
       A pathwise bound for $\displaystyle\max_{1\le n\le N} \|\displaystyle\sum_{j=0}^{n-1}\mathcal{I}_{1,j}^{a,1}\|_{H^\sigma}$ can be achieved by using the $L^{2p}$ estimated given in the proof of Lemma \ref{strong_lem_trunc}, combined with Borel-Cantelli Lemma and classical arguments \cite{deb0}. The estimate for the term $\displaystyle\max_{1\le n\le N} \|\displaystyle\sum_{j=0}^{n-1}\mathcal{I}_{1,j}^{a,2}\|_{H^\sigma}$  can be obtained directly by the computation given in the proof of Lemma \ref{strong_lem_trunc}, since it is already done pathwise. The pathwise estimate for the sum containing $\mathcal{I}_{j,1}^{a,3}$ can be obtained by using the standard fact that $\eta^R(t)$ is a.s. Hölder continuous with rate less than $\frac{1}{2}$ and, finally, the pathwise estimate for the sum containing $\mathcal{I}_{j,1}^b$ arises from the pathwise computation given in the proof of Lemma \ref{strong_lem_trunc}, noting that $\displaystyle\sup_{t\in[0,T]} \|\eta^R(t)\|_{H^\sigma}<\infty, a.s.$ These arguments and an application of the discrete Gronwall inequality completes the proof.
    \end{proof}

    With the aid of Lemma \ref{strong_lem_trunc} and Lemma \ref{path_trunc_lem}, we are in a position to state and prove the following result ensuring pathwise convergence of our integrator (SLR2) applied to Equation \eqref{cub_eq}.
    \begin{proposition}
    \label{conv_path1}
        Let us assume $\psi_0\in H^{\sigma+\gamma}$ and $Q^{\frac{1}{2}}\in\mathcal{L}_2^{\sigma+\gamma}$, for $\sigma\in \mathbb{N}^+$ and $\gamma\in(0,1]$. Then, for any $\delta<\gamma$, there exists a random variable $K_\delta(\omega)$ such that one has
        $$
        \max_{1\le n\le N}\|\psi(t_n)-\psi_n\|_{H^\sigma}\le K_\delta(\omega) \tau^\delta,\qquad \mathbb{P}-{\rm a.s.}
        $$
    \end{proposition}
    \begin{proof}
    First, we prove almost sure convergence. For $Q^{\frac{1}{2}}\in\mathcal{L}_2^{\sigma}$, one has that the process
    $$
    X(t)=\int_0^t {\rm e}^{{\bf i}(t-s)\Delta}Q^{\frac{1}{2}}\,dW(s)
    $$
    is a.s. Hölder continuous in $H^\sigma$ with rate $\frac{1}{2}^-.$ This, in particular, shows that 
    $$
    \max_{1\le n\le N}\|\int_{t_{n-1}}^{t_{n}} {\rm e}^{{\bf i}(t_n-s)\Delta}Q^{\frac{1}{2}}\,dW(s)\|_{H^\sigma} \le C(\omega), \qquad \mathbb{P}-{\rm a.s.},
    $$
    independently on $\tau$. Now, by Theorem \ref{wp_cub}, we have that $\sup_{t\in[0,T]} \|\psi(t)\|_{H^\sigma}<\infty, {\rm a.s.}$ Let $R_0$ a random variable such that
    $$
    \max_{1\le n\le N}\|\int_{t_{n-1}}^{t_{n}} {\rm e}^{{\bf i}(t-s)\Delta}Q^{\frac{1}{2}}\,dW(s)\|_{H^\sigma} +\sup_{t\in[0,T]} \|\psi(t)\|_{H^\sigma} \le R_0.
    $$
    Moreover, let $\tau\le C_{GN}R_0^{-2}$ and $\epsilon\in(0,1)$ and assume
    $$
    \max_{1\le n\le N}\|\psi(t_n)-\psi_n\|_{H^\sigma}\ge \epsilon.
    $$
    Set $n_\epsilon=\min\{n:\|\psi(t_n)-\psi_n\|_{H^\sigma}\ge \epsilon \}$. Then, for $0\le n\le n_\epsilon-1$, one has
    $$
    \begin{aligned}
        \|u_{n}\|_{H^\sigma}&\le \|\psi(t_n)\|_{H^\sigma}+\|\psi(t_n)-\psi_n\|_{H^\sigma}\\
        &\le  \|\psi(t_n)\|_{H^\sigma}+\epsilon\\
        &\le R_0.
    \end{aligned}
    $$
    Hence, By the triangular inequality, the unitary property of the operator ${\rm e}^{{\bf i}t\Delta}$ and the algebra property of $H^\sigma$, one has
    \begin{align*}
\| \psi_{n+1}\|_{H^\sigma} \le 2 R_0&+ \tau \|(\psi^n)^2(\varphi_1(-2 {\bf i} \tau \Delta)\overline{\psi_n})\|_{H^\sigma}\\
&\le 2 R_0+C \tau \|\psi_n\|_{H^\sigma}^3\\
&\le 2 R_0+C \tau R_0^2\|\psi_n\|_{H^\sigma}\\
&\le 3 R_0.
    \end{align*}
   It follows that $\psi_n=\psi^{3R_0}_n$, for $1\le n\le n_\epsilon$. Then, one get
    $$
    \max_{1\le n\le N}\|\psi^{3 R_0}(t_n)-\psi^{3 R_0}_n\|_{H^\sigma}\ge \epsilon, 
    $$
    i.e., an absurd holds since, by Lemma \ref{path_trunc_lem}, this is impossible for $\tau\to 0$. Then, we get
    $$
   \lim_{\tau \to 0} \max_{1\le n\le N}\|\psi(t_n)-\psi_n\|_{H^\sigma}=0, \quad \mathbb{P}-{\rm a.s.}
    $$
    The rest of the proof follows in a similar arguments than those presented in \cite{deb0} and then is omitted.
    \end{proof}

\begin{remark}
    Also for the cubic SSE \eqref{cub_eq}, we get order one convergence in $H^\sigma$ for solutions in $H^{\sigma+1}$, instead of $H^{\sigma+2}$ requested by standard theory.
\end{remark}

\subsection{long-term numerical integration for cubic SSE}
In this section, similarly as done for the linear SSE with potential \eqref{lin_eq}, we aim to provide a numerical integrator suitable for the cubic SSE \eqref{cub_eq} up to large time scales, in presence of $\mathcal{O}(\varepsilon)$ stochastic perturbation of the correspondent deterministic cubic equation with $\mathcal{O}(\varepsilon)$ initial data. We start by recalling that in \cite{fe_ka} the authors highlighted that if the deterministic cubic schrodinger equation is equipped with an $\mathcal{O}(\varepsilon)$ initial data, i.e., $\psi_0(x)=\varepsilon \phi(x)$, then, upon change of variable $w(t,x)=\psi(t,x)/\varepsilon$, one has that $w$ solves the equation
$$
\dot{w}={\bf i}\Delta w-{\bf i}\varepsilon^2 |w|^2 w, \qquad w(0)=\phi(x).
$$
%we here aim to study the long-term strong error of the integrator (SNRLR2) when applied to the sto stochastic cubic Schrodinger equation given by a stochastic perturbation by means of $\mathcal{O}(\varepsilon)$ noise, i.e., 
Then, we here consider the equation
%is governed by an  $\mathcal{O}(\varepsilon)$ noise and by an $\mathcal{O}(\varepsilon^{2})$ nonlinearity, i.e., we consider
\begin{equation}
    \label{cub_small}
    d\psi^\varepsilon(t)={\bf i}\Delta \psi^\varepsilon(t) dt-{\bf i}\varepsilon^2 |\psi^\varepsilon(t)|^2\psi^\varepsilon(t)dt+\varepsilon Q^{\frac{1}{2}}dW(t), \qquad \psi^\varepsilon(0)=\phi(x).
\end{equation}
%By simplicity of treatment, we assume $\phi(x)$ to be deterministic.

We make the following assumption.
\begin{assumption}
\label{smal_cub_ass}
For $\phi\in H^\sigma$ and $Q^{\frac{1}{2}}\in\mathcal{L}_2^\sigma$, for $\sigma\ge 2$, we assume the following bound for the exact solution to  \eqref{cub_small}
$$
\big(\mathbb{E}\big[\displaystyle\sup_{t\in[0,T_\varepsilon]}\|\psi^\varepsilon(t)\|^{2p}_{H^\sigma}\big]\big)^{\frac{1}{2p}}\lesssim 1, \qquad T_\varepsilon=\frac{T}{\varepsilon^2},
$$
for any $p\in\mathbb{N}^+$.
\end{assumption}
%\begin{remark}
 %   Assumption \ref{smal_cub_ass} is reasonable since it holds true in the deterministic setting \cite{fe_ka} and due to the a.s. Hölder continuity in $H^\sigma$ of the stochastic convolution in the mild solution \eqref{mild_cub} with rate $\frac{1}{2}^-$. One indeed has
 %   $$
  %  \varepsilon\| \int_0^{t}{\rm e}^{{\bf i}(t-s)\Delta} Q^{\frac{1}{2}}dW\|_{H^\sigma}\le k(\omega) \varepsilon t^{\frac{1}{2}}\lesssim k(\omega), \quad 0\le t\le \frac{T}{\varepsilon^2}, \ \mathbb{P}-{\rm a.s.},
   % $$
    %provided $Q^{\frac{1}{2}}\in\mathcal{L}_2^\sigma$.
%\end{remark}
Inspired from the analysis given in the previous section and from \cite{fe_ka}, we here slightly modified the integrator (SLR2), integrating exactly the zero mode term in \eqref{four_met}, i.e., we consider the following approximation for the twisted variable $\eta^\varepsilon(t)={\rm e}^{-{\bf i}t\Delta}\psi^\varepsilon(t)$
\begin{align*}
    &\int_0^\tau {\rm e}^{-{\bf i}(t_n+s)\Delta}|{\rm e}^{{\bf i}(t_n+s)\Delta}\eta^\varepsilon(t_n)|^2 {\rm e}^{{\bf i}(t_n+s)\Delta}\eta^\varepsilon(t_n)\,ds\notag \\
    &\quad \approx \bigg(\frac{1}{\sqrt{2\pi}}\bigg)^{3} \sum_{\ell\in\mathbb{Z}} {\rm e}^{{\bf i}\ell x} \sum_{\substack{\ell_1,\ell_2,\ell_3\in\mathbb{Z}\\ \ell=-\ell_1+\ell_2+\ell_3\\ \ell^2+\ell_1^2-\ell_2^2-\ell_3^2\ne 0}}
    {\rm e}^{{\bf i}(\ell^2-\ell_1^2-\ell_2^2-\ell_3^2) t_n}\big(\widehat{\eta^\varepsilon}_{\ell_1}(t_n)\big)^*
\widehat{\eta^\varepsilon}_{\ell_2}(t_n)\widehat{\eta^\varepsilon}_{\ell_3}(t_n)\int_0^\tau {\rm e}^{2{\bf i}s\ell_1^2}\,ds \\
&\hspace{4cm}+\tau \bigg(\frac{1}{\sqrt{2\pi}}\bigg)^{3} \sum_{\ell\in\mathbb{Z}} {\rm e}^{{\bf i}\ell x}  \sum_{\substack{\ell_1,\ell_2,\ell_3\in\mathbb{Z}\\ \ell=-\ell_1+\ell_2+\ell_3\\ \ell^2+\ell_1^2-\ell_2^2-\ell_3^2= 0}}  \big(\widehat{\eta^\varepsilon}_{\ell_1}(t_n)\big)^*
\widehat{\eta^\varepsilon}_{\ell_2}(t_n)\widehat{\eta^\varepsilon}_{\ell_3}(t_n).
\end{align*}
After same computations than those already presented in \cite{fe_ka} and not reported here for the sake of brevity, we get the following non-resonant low-regularity method for the twisted variable $\eta^\varepsilon$
\begin{align*}
\eta^\varepsilon_{n+1}=\eta^\varepsilon_n&-{\bf i} \tau\varepsilon^2\big[{\rm e}^{-{\bf i}t_n\Delta}\big(({\rm e}^{{\bf i}t_n\Delta} \eta^\varepsilon_n)^2 \varphi_1(-2{\bf i}\tau\Delta){\rm e}^{-{\bf i}t_n\Delta} \overline{\eta^\varepsilon_n}\big)
  +2\widehat{\big(g(\eta^\varepsilon_n)\big)}_0 \eta^\varepsilon_n-h(\eta^\varepsilon_n)\big]\\
   &+\varepsilon\int_{t_n}^{t_n+\tau} {\rm e}^{-{\bf i}s\Delta}Q^{\frac{1}{2}}\,dW(s),
   \end{align*}
   with the functions $g$ and $h$ defined as 
   $$g(u)=u(1-\varphi_1(-2{\bf i}\tau\Delta))\overline{u}, \quad (\widehat{h(u)})_{\ell}=(1-\varphi_1(2{\bf i}\tau \ell^2))\overline{\hat{u}_{\ell}}\hat{u}_{\ell}\hat{u}_{\ell}, \ \ell\in\mathbb{Z}.
   $$ 
Returning back to the original variable $\psi^\varepsilon$, we get the following integrator for \eqref{cub_small}, named (SNRLR2),
    \begin{align}
    \psi^\varepsilon_{n+1}&={\rm e}^{{\bf i}\tau\Delta}\big[\psi^\varepsilon_n-{\bf i}\tau\varepsilon^2 (\psi^\varepsilon_n)^2\varphi_1(-2{\bf i}\tau\Delta)\overline{\psi^\varepsilon_n}\big]+\varepsilon\mathcal{W}_n\notag\\
    &\hspace{3cm}-2{\bf i}\tau\varepsilon^2 \widehat{\big(g(\psi^\varepsilon_n)\big)}_0 {\rm e}^{{\bf i}\tau\Delta} \psi^\varepsilon_n+{\bf i}\tau\varepsilon^2 {\rm e}^{{\bf i}\tau\Delta} h(\psi^\varepsilon_n)    \label{SLNR2} .
        \end{align}
    \begin{remark}
        It is not difficult to (see, e.g., also \cite{fe_ka}) that the convergence results shown for the methods (SLR2) also hold true for the method (SNRLR2), i.e., it converges in $L^{2p}(\Omega; H^\sigma), \ \sigma\in\mathbb{N}^+$, with order $\gamma\in(0,1]$, for solutions in $H^{\sigma+\gamma}$, when applied to the truncated equation \eqref{trunc_cub}, and,
        by Proposition \ref{conv_path1}, it convergences $\mathbb{P}-{\rm a.s.}$ with order $\gamma^-$ in $H^\sigma$ for solutions in $H^{\sigma+\gamma}$, when applied to the original cubic equation \eqref{cub_eq}. However, for the non-resonant method (SNRLR2), we can exploit (RCO) technique \cite{bao,fe_ka} and get improved long-term error bounds under Assumption \ref{smal_cub_ass}. This is the aim of the next subsection.
        \end{remark}
        \begin{remark}
            As announced in the introduction, the application of (RCO) technique in our stochastic setting presents more difficulties with respect to the deterministic case \cite{fe_ka}, if one wants to avoid the fact that the presence of the stochastic convolution forces the exact solutions to \eqref{cub_small} to be Hölder continuous with rate one-half.
        \end{remark}
We start by providing an auxiliary result on the strong long term error estimation of the numerical method (SNRLR2) applied to the truncated version of Equation \eqref{cub_small}, i.e., to the equation
\begin{equation}
    \label{trunc_cub2}
 d\psi^{\varepsilon,R}={\bf i}\Delta \psi^{\varepsilon,R} dt-{\bf i}\varepsilon^2 \theta_R(\psi^{\varepsilon,R})|\psi^{\varepsilon,R}|^2\psi^{\varepsilon,R} dt+\varepsilon Q^{\frac{1}{2}}dW, \qquad \psi^{\varepsilon,R}(0)=\phi(x),
\end{equation}
under Assumption \ref{smal_cub_ass}, with $\theta_R(u)=\theta(\|u\|_{H^1}/R), \ R>0$.

Let us denote $\varphi^{\varepsilon,R,\tau}$ the flow map such that $\eta^{\varepsilon,R}(t_n+\tau)=\varphi^{\varepsilon,R,\tau}(\eta^{\varepsilon,R}(t_n))$ and $\Psi^{\varepsilon,R,\tau}_{t_n}$ the map such that $\eta^{\varepsilon,R}_{n+1}=\Psi_{t_n}^{\varepsilon,R,\tau}(\eta^{\varepsilon,R}_n)$, with $\eta^{\varepsilon,R}$ twisted variable.
We prove the following results.
\begin{lemma}
\label{long_term_cub_trunc}
    Let us assume Assumption \ref{smal_cub_ass}. For any fixed $\tau_0\in(0,1)$ independent on $\varepsilon$, there exists a constant $K(\tau_0)$ such that, for any $\tau\in(0,K(\tau_0))$, the numerical integrator (SNRLR2), applied to the truncated cubic equation \eqref{trunc_cub2}, satisfies the following error estimate
    \begin{equation}
        \label{long_term_mean_trunc_cub}
    \big(\mathbb{E}\big[\max_{1\le n\le N_\varepsilon}\|\psi^{\varepsilon,R}(t_n)-\psi^{\varepsilon,R}_n\|^2_{H^1}\big]\big)^{\frac{1}{2}}\lesssim\tau_0^{\sigma-1}+\varepsilon^2\tau, \qquad N_\varepsilon:=\frac{T}{\varepsilon^2 \tau}.
    \end{equation}
    In addition, if we assume a setting in which the term $\tau_0^{\sigma-1}$ in \eqref{long_term_mean_trunc_cub} is neglectable, there exists a sequence $\{\tau_k\}_{k\in\mathbb{N}^+}$, with $\tau_k\to 0, \ k\to \infty$, such that, for any $R>0$ and $\delta\in(0,1)$, one has
    \begin{equation}
        \label{long_term_path_trunc_cub}
        \max_{1\le n\le N^k_\varepsilon}\|\psi^{\varepsilon,R}(t_n)-\psi^{\varepsilon,R}_n\|_{H^1}<\varepsilon^2 \tau_k^\delta, \quad \forall k>k_0(\omega), \qquad \mathbb{P}-{\rm a.s.},
    \end{equation}
    where $N_\varepsilon^k=\frac{T}{\varepsilon^2 \tau_k}$, with $k_0(\omega)$ independent on $\varepsilon$.
\end{lemma}
\begin{proof}
    %Some lines of the proof will follows similar techniques of the proof of Theorem \ref{lobg_term_res1}, the proof of Lemma \ref{strong_lem_trunc} and of the proof of Theorem 3.1 in \cite{fe_ka}. Here, we only spot the differences. 
    We first give the result \eqref{long_term_mean_trunc_cub} for the twisted variable $\eta^{\varepsilon,R}$. We define $\mathcal{E}^{\varepsilon,R}_n=\Psi_{t_n}^{\varepsilon,R,\tau} (\eta^{\varepsilon,R}(t_n))-\varphi^{\varepsilon,R,\tau}(\eta^{\varepsilon,R}(t_n))$ and $e^{\varepsilon,R}_n=\eta^{\varepsilon,R}_n-\eta^{\varepsilon,R}(t_n)$. Then, we have
    $$
    e^{\varepsilon,R}_{n+1}=e^{\varepsilon,R}_n+Z^{\varepsilon,R}_n+\mathcal{E}^{\varepsilon,R}_n,
    $$
    where 
    \begin{align*}
        Z^{\varepsilon,R}_n=-{\bf i}\varepsilon^2 \tau & \big[\theta_R(\eta^{\varepsilon,R}_n){\rm e}^{-{\bf i}t_n\Delta}\big(({\rm e}^{{\bf i}t_n\Delta}\eta^{\varepsilon,R}_n)^2 \varphi_1(-2{\bf i}\tau\Delta){\rm e}^{-{\bf i}t_n\Delta}\overline{\eta^{\varepsilon,R}_n}\big) \\
        &+\theta_R(\eta^{\varepsilon,R}_n)\big(2(\widehat{g(\eta^{\varepsilon,R}_n)})_0 \eta^{\varepsilon,R}_n-h(\eta^{\varepsilon,R}_n)\big)\\
        & - \theta_R(\eta^{\varepsilon,R}(t_n)){\rm e}^{-{\bf i}t_n\Delta}\big(({\rm e}^{{\bf i}t_n\Delta}\eta^{\varepsilon,R}(t_n))^2 \varphi_1(-2{\bf i}\tau\Delta){\rm e}^{-{\bf i}t_n\Delta}\overline{\eta^{\varepsilon,R}(t_n)}\big) \\
        &-\theta_R(\eta^{\varepsilon,R}(t_n))\big(2(\widehat{g(\eta^{\varepsilon,R}(t_n))})_0 \eta^{\varepsilon,R}(t_n)+h(\eta^{\varepsilon,R}(t_n))\big)\big]
    \end{align*}
    and $\mathcal{E}^{\varepsilon,R}_{n}=\mathcal{V}^{\varepsilon,R}_n+\mathcal{R}^{\varepsilon,R}(\eta^{\varepsilon,R}(t_n))$ with
    \begin{align*}
        &\mathcal{R}^{\varepsilon,R}(\eta^{\varepsilon,R}(t_n))={\bf i}\varepsilon^2\theta_R(\eta^{\varepsilon,R}(t_n)) \displaystyle\sum_{\ell \in\mathbb{Z}}{\rm e}^{{\bf i}\ell x}\displaystyle\sum_{\substack{\ell_1,\ell_2,\ell_3 \in\mathbb{Z} \\ \ell=-\ell_1+\ell_2+\ell_3\\ \ell^2+\ell_1^2-\ell_2^2-\ell_3^2\ne 0}}{\rm e}^{{\bf i}t_n (\ell^2+\ell_1^2-\ell_2^2-\ell_3^2)}\\ 
        &\hspace{2cm}\times \big(\widehat{v^{\varepsilon,R}}_{\ell_1}(t_n)\big)^*\widehat{v^{\varepsilon,R}}_{\ell_2}(t_n)\widehat{v^{\varepsilon,R}}_{\ell_3}(t_n)
         \int_0^\tau  {\rm e}^{2{\bf i}s \ell_1^2}-{\rm e}^{{\bf i}s (\ell^2+\ell_1^2-\ell_2^2-\ell_3^2)}\,ds
    \end{align*}
and
\begin{align*}
   & \mathcal{V}^{\varepsilon,R}_n = {\bf i}\varepsilon^2 \int_0^\tau {\rm e}^{-{\bf i}(t_n+s)\Delta}\theta_R(\eta^{\varepsilon,R}(t_n+s))\big|{\rm e}^{{\bf i}(t_n+s)\Delta}\eta^{\varepsilon,R}(t_n+s)\big|^2 {\rm e}^{{\bf i}(t_n+s)\Delta}\eta^{\varepsilon,R}(t_n+s)\,ds\\
    &\hspace{2cm}-{\bf i}\varepsilon^2 \int_0^\tau {\rm e}^{-{\bf i}(t_n+s)\Delta}\theta_R(\eta^{\varepsilon,R}(t_n))\big|{\rm e}^{{\bf i}(t_n+s)\Delta}\eta^{\varepsilon,R}(t_n)\big|^2 {\rm e}^{{\bf i}(t_n+s)\Delta}\eta^{\varepsilon,R}(t_n)\,ds.
\end{align*}
Then, by recursion, we get
$$
e^{\varepsilon,R}_{n}=\sum_{k=0}^{n-1}Z^{\varepsilon,R}_k+\sum_{k=0}^{n-1}\mathcal{R}^{\varepsilon,R}(\eta^{\varepsilon,R}(t_k))+\sum_{k=0}^{n-1}\mathcal{V}^{\varepsilon,R}_k.
$$
Then, we get
\begin{align*}
\mathbb{E}\big[\max_{1\le n\le N_\varepsilon}\|e^{\varepsilon,R}_{n}\|^{2}_{H^1}\big] \lesssim \mathbb{E}\big[\max_{1\le n\le N_\varepsilon}\|\sum_{k=0}^{n-1}Z^{\varepsilon,R}_{k}\|^{2}_{H^1}\big]&+\mathbb{E}\big[\max_{1\le n\le N_\varepsilon}\|\sum_{k=0}^{n-1}\mathcal{V}^{\varepsilon,R}_k\|^{2}_{H^1}\big]\\
&+\mathbb{E}\big[\max_{1\le n \le N_\varepsilon}\|\sum_{k=0}^{n-1}\mathcal{R}^{\varepsilon,R}(\eta^{\varepsilon,R}(t_k))\|^{2}_{H^1}\big].
\end{align*}
In the same manner of the estimation of the term $\mathcal{I}_{2,j}$ in the proof of Lemma \ref{strong_lem_trunc}, we get
\begin{align*}
\mathbb{E}\big[\max_{1\le n\le N_\varepsilon}\|\sum_{k=0}^{n-1}Z^{\varepsilon,R}_{k}\|^{2}_{H^1}\big]& = (\varepsilon^2\tau)^{2} \mathbb{E}\big[\max_{1\le n\le N_\varepsilon}\|\sum_{k=0}^{n-1}\tilde{F}(\eta^{\varepsilon,R}(t_k))-\tilde{F}(\eta^{\varepsilon,R}_k)\|^{2}_{H^1}\big]\\
&\lesssim (\varepsilon^2\tau)^{2}N_\varepsilon\sum_{k=0}^{N_\varepsilon-1}\mathbb{E}\big[\|\eta^{\varepsilon,R}(t_k)-\eta^{\varepsilon,R}_k\|^{2}_{H^1}\big]\\
&\lesssim \varepsilon^2 \tau \sum_{k=0}^{N_\varepsilon-1}\mathbb{E}\big[\max_{0\le l\le k}\|e^{\varepsilon,R}_l\|^{2}_{H^1}\big].
\end{align*}
Moreover, by similar computation than those shown for the estimation of the term $\mathcal{I}_{1,j}^{a}$ in the proof of Lemma \ref{strong_lem_trunc} and using Assumption \ref{smal_cub_ass} with $\sigma=2$ (note that $\sigma \ge 2$) one has
\begin{align*}
    \mathbb{E}\big[\max_{1\le n\le N_\varepsilon}\|\sum_{k=0}^{n-1}\mathcal{V}^{\varepsilon,R}_k\|^{2}_{H^1}\big]&= \varepsilon^4\mathbb{E}\big[\max_{1\le n\le N_\varepsilon}\|\sum_{k=0}^{n-1}\mathcal{I}_{1,k}^a\|^{2}_{H^1}\big]\\
    &\lesssim \varepsilon^4\big(\varepsilon^{2}\mathbb{E}\big[\max_{1\le n \le N_\varepsilon}\|\sum_{k=0}^{n-1}\mathcal{I}_{1,k}^{a,1}\|^{2}_{H^1}\big]+\varepsilon^{4}\mathbb{E}\big[\max_{1\le n\le N_\varepsilon}\|\sum_{k=0}^{n-1}\mathcal{I}_{1,k}^{a,2}\|^{2}_{H^1}\big]\\
    &\hspace{3cm}+\mathbb{E}\big[\max_{1\le n\le N_\varepsilon}\|\sum_{k=0}^{n-1}\mathcal{I}_{1,k}^{a,3}\|^{2}_{H^1}\big]\big)\\
    &\lesssim  \varepsilon^{8}\|Q^{\frac{1}{2}}\|^{2}_{\mathcal{L}_2^2}\sum_{k=0}^{N_\varepsilon-1} \int_{t_k}^{t_{k+1}} |t_{k+1}-r|^2\,dr
    + \varepsilon^{8}(N_\varepsilon)^{2} \tau^{4}\\
    &\hspace{1cm}+\varepsilon^{4}N_\varepsilon \tau\sum_{k=0}^{N_\varepsilon-1}\int_0^\tau \int_0^1 \mathbb{E}\big[\|\eta^{\varepsilon,R}(t_k+s)-\eta^{\varepsilon,R}(t_k)\|^{4}_{H^1}\big]\,d\mu\,ds.
\end{align*}
It is direct to verify that the following estimate holds true
\begin{equation}
    \label{t_reg_cub}
\|\eta^{\varepsilon,R}(t_k+s)-\eta^{\varepsilon,R}(t_k)\|_{L^{2q}(\Omega; H^\sigma)}\lesssim \varepsilon s^{\frac{1}{2}}, \quad q,\sigma\in\mathbb{N}^+
\end{equation}
as long as $\eta^{\varepsilon,R}$ satisfies Assumption \ref{smal_cub_ass}. Using this, we then get
\begin{align*}
    \mathbb{E}\big[\max_{1\le n\le N_\varepsilon}\|\sum_{k=0}^{n-1}\mathcal{V}^{\varepsilon,R}_k\|^{2}_{H^1}\big]&\lesssim \varepsilon^{6} N_\varepsilon \tau^3+(\varepsilon^2 \tau)^{2}
    +\varepsilon^4(N_\varepsilon\tau)^{2}\varepsilon^{4}\tau^{2}\\
    &\lesssim (\varepsilon^2 \tau)^{2}.
\end{align*}
Then, we get
\begin{equation*}
\mathbb{E}\big[\max_{1\le n \le N_\varepsilon}\|e^{\varepsilon,R}_{n}\|^{2}_{H^1}\big] \lesssim  (\varepsilon^2 \tau)^{2}+\varepsilon^2 \tau \sum_{k=0}^{N_\varepsilon-1} \mathbb{E}\big[\max_{0\le l\le k}\|e^{\varepsilon,R}_l\|^{2}_{H^1}\big]+
\mathbb{E}\big[\max_{1\le n\le N_\varepsilon}\|\sum_{k=0}^{n-1}\mathcal{R}^{\varepsilon,R}(\eta^{\varepsilon,R}(t_k))\|^{2}_{H^1}\big].
\end{equation*}
We have to estimate the term $\mathbb{E}\big[\displaystyle\max_{1\le n \le N_\varepsilon}\|\sum_{k=0}^{n-1}\mathcal{R}^{\varepsilon,R}(\eta^{\varepsilon,R}(t_k))\|^{2}_{H^1}\big]$. Let $\tau_0$, $N_0$ and $\mathcal{P}_{N_0}$ defined as in the proof of Theorem \ref{lobg_term_res1}. Then
\begin{align*}
\mathbb{E}\big[\max_{1\le n\le N_\varepsilon}\|\sum_{k=0}^{n-1}\mathcal{R}^{\varepsilon,R}(\eta^{\varepsilon,R}(t_k))\|^{2}_{H^1}\big]&\lesssim \mathbb{E}\big[\max_{1\le n\le N_\varepsilon}\|\sum_{k=0}^{n-1}\mathcal{R}(\eta^{\varepsilon,R}(t_k))-\mathcal{P}_{N_0}\mathcal{R}^{\varepsilon,R}(\eta^{\varepsilon,R}(t_k))\|^{2}_{H^1}\big]\\
&\hspace{3.1cm}+\mathbb{E}\big[\max_{1\le n\le N_\varepsilon}\|\sum_{k=0}^{n-1}\mathcal{P}_{N_0}\mathcal{R}^{\varepsilon,R}(\eta^{\varepsilon,R}(t_k))\|^{2}_{H^1}\big].
\end{align*}
By the definition of $\mathcal{R}^{\varepsilon,R}(\eta^{\varepsilon,R}(t_k))$ and $\mathcal{P}_{N_0}$, thanks to the algebra property of $H^\sigma$ and due to Assumption \ref{smal_cub_ass}, we have
\begin{align*}
    &\mathbb{E}\big[\max_{1\le n\le N_\varepsilon}\|\sum_{k=0}^{n-1}\mathcal{R}^{\varepsilon,R}(\eta^{\varepsilon,R}(t_k))-\mathcal{P}_{N_0}\mathcal{R}^{\varepsilon,R}(\eta^{\varepsilon,R}(t_k))\|^{2}_{H^1}\big]\\
    &\hspace{0.2cm}\lesssim \varepsilon^{4}N_\varepsilon\sum_{k=0}^{N_\varepsilon-1}\mathbb{E}\big[\displaystyle\sum_{|\ell|\ge \frac{N_0}{2}+1}(1+|\ell|^2)\big|\displaystyle\sum_{\substack{\ell_1,\ell_2,\ell_3 \in\mathbb{Z} \\ \ell=-\ell_1+\ell_2+\ell_3\\ \ell^2+\ell_1^2-\ell_2^2-\ell_3^2\ne 0}}{\rm e}^{{\bf i}t_k (\ell^2+\ell_1^2-\ell_2^2-\ell_3^2)}\big(\widehat{\eta^{\varepsilon,R}}_{\ell_1}(t_k)\big)^*\widehat{\eta^{\varepsilon,R}}_{\ell_2}(t_k)\widehat{\eta^{\varepsilon,R}}_{\ell_3}(t_k)\\
        &\hspace{6cm}\times \int_0^\tau  {\rm e}^{2{\bf i}s \ell_1^2}-{\rm e}^{{\bf i}s (\ell^2+\ell_1^2-\ell_2^2-\ell_3^2)}\,ds\big|^2\big]\\
        &\hspace{0.2cm}\lesssim \varepsilon^{4}N_\varepsilon \tau_0^{2(\sigma-1)}\tau^{2}\sum_{k=0}^{N_\varepsilon-1}\mathbb{E}\big[\displaystyle\sum_{|\ell|\ge \frac{N_0}{2}+1}(1+|\ell|^{2\sigma})\big(\displaystyle\sum_{\substack{\ell_1,\ell_2,\ell_3 \in\mathbb{Z} \\ \ell=-\ell_1+\ell_2+\ell_3\\ \ell^2+\ell_1^2-\ell_2^2-\ell_3^2\ne 0}}|\widehat{\eta^{\varepsilon,R}}_{\ell_1}(t_k)| |\widehat{\eta^{\varepsilon,R}}_{\ell_2}(t_k)| |\widehat{\eta^{\varepsilon,R}}_{\ell_3}(t_k)|)^2\big]\\
        &\hspace{0.2cm}\lesssim \varepsilon^{4}N_\varepsilon\tau_0^{2(\sigma-1)}\tau^{2}\sum_{k=0}^{N_\varepsilon-1}\mathbb{E}\big[\|\eta^{\varepsilon,R}(t_k)\|^{6}_{H^\sigma}\big]\\
        &\hspace{0.2cm}\lesssim \varepsilon^{4}N_\varepsilon^{2}\tau_0^{2(\sigma-1)}\tau^{2}\\
        &\hspace{0.2cm}\lesssim \tau_0^{2(\sigma-1)}.
\end{align*}
Then, we end up with 
\begin{align}
\mathbb{E}\big[\max_{1\le n\le N_\varepsilon}\|e^{\varepsilon,R}_{n}\|^{2}_{H^1}\big] \lesssim  \tau_0^{2(\sigma-1)}+(\varepsilon^2 \tau)^{2}&+\varepsilon^2 \tau\sum_{k=0}^{N_\varepsilon-1} \mathbb{E}\big[\max_{0\le l \le k}\|e^{\varepsilon,R}_l\|^{2}_{H^1}\big]\notag \\
&+\mathbb{E}\big[\max_{1\le n\le N_\varepsilon}\|\sum_{k=0}^{n-1}\mathcal{P}_{N_0}\mathcal{R}^{\varepsilon,R}(\eta^{\varepsilon,R}(t_k))\|^{2}_{H^1}\big]  \label{in_es0}.
\end{align}
 Similarly to proof of Theorem \ref{lobg_term_res1},  by the definition of $\mathcal{R}(v^R(t_k))$, we see
$$
\sum_{k=0}^{n-1} \mathcal{P}_{N_0}\mathcal{R}^{\varepsilon,R}(\eta^{\varepsilon,R}(t_k))={\bf i}\varepsilon^2 \sum_{k=0}^{n-1} \displaystyle\sum_{\ell \in \mathcal{T}_{N_0}}{\rm e}^{{\bf i}\ell x}\displaystyle\sum_{(\ell_1,\ell_2,\ell_3)\in\mathcal{I}_\ell^{N_0}} \Lambda^{\varepsilon,R}_{k,\ell,\ell_1,\ell_2,\ell_3},
$$
where
$
\mathcal{T}_{N_0}:=\{\ell: -\frac{N_0}{2},\dots,\frac{N_0}{2}\}
$
and, for $\ell \in \mathcal{T}_{N_0}$, 
$$
\mathcal{I}_\ell^{N_0}:=\{(\ell_1,\ell_2,\ell_2): -\ell_1+\ell_2+\ell_3=\ell,\  \ell_1,\ell_2,\ell_3\in \mathcal{T}_{N_0}\}.
$$
Without lost of generality, we consider $n\ge 2$. The bound for $n=1$, i.e., at single step, arises similarly.
Denoting $\delta_{\ell,\ell_1,\ell_2,\ell_3}=\ell^2+\ell_1^2-\ell_2^2-\ell_3^2$, note that $\Lambda^{\varepsilon,R}_{k,\ell,\ell_1,\ell_2,\ell_3}=0$ if $\delta_{\ell,\ell_1,\ell_2,\ell_3}=0$ and 
$$
\Lambda^{\varepsilon,R}_{k,\ell,\ell_1,\ell_2,\ell_3}=r_{\ell,\ell_1,\ell_2,\ell_3} {\rm e}^{{\bf i}t_k\delta_{\ell,\ell_1,\ell_2,\ell_3}} c^{\varepsilon,R}_{k,\ell,\ell_1,\ell_2,\ell_3}
$$
with
\begin{align*}
&r_{\ell,\ell_1,\ell_2,\ell_3}=\int_0^\tau {\rm e}^{2{\bf i}s \ell_1^2}\,ds-\int_0^\tau {\rm e}^{{\bf i}s \delta_{\ell,\ell_1,\ell_2,\ell_3}}\,ds,\\
& c^{\varepsilon,R}_{k,\ell,\ell_1,\ell_2,\ell_3}
=\theta_R(\eta^{\varepsilon,R}(t_k))\big(\widehat{\eta^{\varepsilon,R}}_{\ell_1}(t_k)\big)^*\widehat{\eta^{\varepsilon,R}}_{\ell_2}(t_k) \widehat{\eta^{\varepsilon,R}}_{\ell_3}(t_k).
\end{align*}
Note that a standard interpolation argument gives the bound
$$
r_{\ell,\ell_1,\ell_2,\ell_3}=\mathcal{O}(\tau^2 |\delta_{\ell,\ell_1,\ell_2,\ell_3}-2\ell_1^2|).
$$
Moreover, also here we have that, for $\ell_\in \mathcal{T}_{N_0}$ and $\ell_1,\ell_2, \ell_3 \in \mathcal{I}_{\ell}^{N_0}$, one has
$$
|\delta_{\ell,\ell_1,\ell_2,\ell_3}|=\le  \frac{N_0^2}{2}\le \frac{2(1+\tau_0)^2}{\tau_0^2}.
$$
This implies that for sufficiently small $\tau$ (depending on $\tau_0 \in (0,1)$) one has 
$$
\frac{\tau}{2}|\delta_{\ell,\ell_1,\ell_2,\ell_3}|<\pi.
$$
Hence, again denoting $S_{n,\ell,\ell_1,\ell_2.\ell_3}=\sum_{k=0}^n {\rm e}^{{\bf i}t_k\delta_{\ell,\ell_1,\ell_2.\ell_3}}$ and using summation-by-parts formula, we get
\begin{align*}
\sum_{k=0}^{n-1}\Lambda^{\varepsilon,R}_{k,\ell_1,\ell_2,\ell_3}&=r_{\ell,\ell_1,\ell_2,\ell_3}\sum_{k=0}^{n-2}S_{k,\ell,\ell_1,\ell_2,\ell_3}\big(c^{\varepsilon,R}_{k,\ell,\ell_1,\ell_2,\ell_3}-c^{\varepsilon,R}_{k+1,\ell,\ell_1,\ell_2,\ell_3}\big)\\
&\qquad\qquad +S_{n-1,\ell,\ell_1,\ell_2,\ell_3}r_{\ell,\ell_1,\ell_2,\ell_3}c^{\varepsilon,R}_{n-1,\ell,\ell_1,\ell_2,\ell_3},
\end{align*}
where
\begin{align*}
    &c^{\varepsilon,R}_{k,\ell,\ell_1,\ell_2,\ell_3}-c^{\varepsilon,R}_{k+1,\ell,\ell_1,\ell_2,\ell_3} = \theta_R(\eta^{\varepsilon,R}(t_k)) \big(\widehat{\eta^{\varepsilon,R}}_{\ell_1}(t_k)\big)^*(\widehat{\eta^{\varepsilon,R}}_{\ell_2}(t_k)-\widehat{\eta^{\varepsilon,R}}_{\ell_2}(t_{k+1}))\widehat{\eta^{\varepsilon,R}}_{\ell_3}(t_k)\\
   & \hspace{2cm}+ \theta_R(\eta^R(t_k)) (\widehat{\eta^{\varepsilon,R}}_{\ell_1}(t_k)-\widehat{\eta^{\varepsilon,R}}_{\ell_1}(t_{k+1}))^*\widehat{\eta^{\varepsilon,R}}_{\ell_2}(t_{k+1})\widehat{\eta^{\varepsilon,R}}_{\ell_3}(t_k)\\
   &\hspace{2cm}+(\widehat{\eta^{\varepsilon,R}}_{\ell_1}(t_{k+1}))^*\widehat{\eta^{\varepsilon,R}}_{\ell_2}(t_{k+1})\theta_R(\eta^{\varepsilon,R}(t_k))(\widehat{\eta^{\varepsilon,R}}_{\ell_3}(t_k)-\widehat{\eta^{\varepsilon,R}}_{\ell_3}(t_{k+1}))\\
   &\hspace{2cm}+(\widehat{\eta^{\varepsilon,R}}_{\ell_1}(t_{k+1}))^*\widehat{\eta^{\varepsilon,R}}_{\ell_2}(t_{k+1})\widehat{\eta^{\varepsilon,R}}_{\ell_3}(t_{k+1})(\theta_R(\eta^{\varepsilon,R}(t_k))-\theta_R(\eta^{\varepsilon,R}(t_{k+1}))).
\end{align*}
Using this last formula, one achieves
\begin{equation}
    \label{p_0}
    \sum_{k=0}^{n-1} \mathcal{P}_{N_0} \mathcal{R}^{\varepsilon,R}(\eta^{\varepsilon,R}(t_k)) = \mathcal{I}^{n-1}_1+\mathcal{I}^{n-1}_2+\mathcal{I}^{n-1}_3+\mathcal{I}^{n-1}_4+\mathcal{I}^{n-1}_5,
\end{equation}
where
\begin{align*}
    &\mathcal{I}_1^{n-1}:={\bf i}\varepsilon^2 \sum_{\ell\in\mathcal{T}_{N_0}}{\rm e}^{{\bf i}\ell x}\sum_{(\ell_1,\ell_2,\ell_3)\in\mathcal{I}_{\ell}^{N_0}}
    r_{\ell,\ell_1,\ell_2,\ell_3}\sum_{k=0}^{n-2}S_{k,\ell,\ell_1,\ell_2,\ell_3}\\
    &\hspace{4cm}\times \theta_R(\eta^{\varepsilon,R}(t_k)) (\widehat{\eta^{\varepsilon,R}}_{\ell_1}(t_k))^*(\widehat{\eta^{\varepsilon,R}}_{\ell_2}(t_k)-\widehat{\eta^{\varepsilon,R}}_{\ell_2}(t_{k+1}))\widehat{\eta^{\varepsilon,R}}_{\ell_3}(t_k),\\
    &\mathcal{I}_2^{n-1}:={\bf i}\varepsilon^2 \sum_{\ell\in\mathcal{T}_{N_0}}{\rm e}^{{\bf i}\ell x}\sum_{(\ell_1,\ell_2,\ell_3)\in\mathcal{I}_{\ell}^{N_0}}
    r_{\ell,\ell_1,\ell_2,\ell_3}\sum_{k=0}^{n-2}S_{k,\ell,\ell_1,\ell_2,\ell_3}\\
    &\hspace{4cm}\times \theta_R(\eta^{\varepsilon,R}(t_k)) (\widehat{\eta^{\varepsilon,R}}_{\ell_1}(t_k)-\widehat{\eta^{\varepsilon,R}}_{\ell_1}(t_{k+1}))^*\widehat{\eta^{\varepsilon,R}}_{\ell_2}(t_{k+1})\widehat{\eta^{\varepsilon,R}}_{\ell_3}(t_k),\\
    &\mathcal{I}_3^{n-1}:={\bf i}\varepsilon^2 \sum_{\ell\in\mathcal{T}_{N_0}}{\rm e}^{{\bf i}\ell x}\sum_{(\ell_1,\ell_2,\ell_3)\in\mathcal{I}_{\ell}^{N_0}}
    r_{\ell,\ell_1,\ell_2,\ell_3}\sum_{k=0}^{n-2}S_{k,\ell,\ell_1,\ell_2,\ell_3}\\
    &\hspace{4cm}\times (\widehat{\eta^{\varepsilon,R}}_{\ell_1}(t_{k+1}))^*\widehat{\eta^{\varepsilon,R}}_{\ell_2}(t_{k+1})\theta_R(\eta^{\varepsilon,R}(t_k))(\widehat{\eta^{\varepsilon,R}}_{\ell_3}(t_k)-\widehat{\eta^{\varepsilon,R}}_{\ell_3}(t_{k+1})),\\
    &\mathcal{I}_4^{n-1}:={\bf i}\varepsilon^2 \sum_{\ell\in\mathcal{T}_{N_0}}{\rm e}^{{\bf i}\ell x}\sum_{(\ell_1,\ell_2,\ell_3)\in\mathcal{I}_{\ell}^{N_0}}
    r_{\ell,\ell_1,\ell_2,\ell_3}\sum_{k=0}^{n-2}S_{k,\ell,\ell_1,\ell_2,\ell_3}\\
    &\hspace{4cm}\times (\widehat{\eta^{\varepsilon,R}}_{\ell_1}(t_{k+1}))^*\widehat{\eta^{\varepsilon,R}}_{\ell_2}(t_{k+1})\widehat{\eta^{\varepsilon,R}}_{\ell_3}(t_{k+1})(\theta_R(\eta^{\varepsilon,R}(t_k))-\theta_R(\eta^{\varepsilon,R}(t_{k+1}))),\\ 
    &\mathcal{I}_5^{n-1}:={\bf i}\varepsilon^2 \sum_{\ell\in\mathcal{T}_{N_0}}{\rm e}^{{\bf i}\ell x}\sum_{(\ell_1,\ell_2,\ell_3)\in\mathcal{I}_{\ell}^{N_0}}
    S_{n-1,\ell,\ell_1,\ell_2,\ell_3}r_{\ell,\ell_1,\ell_2,\ell_3}\\
    &\hspace{4cm}\times \theta_R(\eta^{\varepsilon,R}(t_{n-1}))(\widehat{\eta^{\varepsilon,R}}_{\ell_1}(t_{n-1}))^*\widehat{\eta^{\varepsilon,R}}_{\ell_2}(t_{n-1}) \widehat{\eta^{\varepsilon,R}}_{\ell_3}(t_{n-1}).
\end{align*}
Based on the above decomposition, we have
\begin{align}
    \mathbb{E}\big[\max_{1\le n\le N_\varepsilon}\|\sum_{k=0}^{n-1}\mathcal{P}_{N_0}\mathcal{R}^{\varepsilon,R}(\eta^{\varepsilon,R}(t_k))\|^{2}_{H^1}\big]&\lesssim \mathbb{E}\big[\max_{1\le n\le N_\varepsilon}\|\sum_{k=0}^{n-1}\mathcal{I}_{1}^{n-1}\|^{2}_{H^1}\big]+\mathbb{E}\big[\max_{1\le n\le N_\varepsilon}\|\sum_{k=0}^{n-1}\mathcal{I}_{2}^{n-1}\|^{2}_{H^1}\big]\notag \\
    &\quad + \mathbb{E}\big[\max_{1\le n\le N_\varepsilon}\|\sum_{k=0}^{n-1}\mathcal{I}_{3}^{n-1}\|^{2}_{H^1}\big]+\mathbb{E}\big[\max_{1\le n\le N_\varepsilon}\|\sum_{k=0}^{n-1}\mathcal{I}_{4}^{n-1}\|^{2}_{H^1}\big]\notag \\
    &\quad +\mathbb{E}\big[\max_{1\le n\le N_\varepsilon}\|\sum_{k=0}^{n-1}\mathcal{I}_{5}^{n-1}\|^{2}_{H^1}\big] \label{start_dec}.
\end{align}
From \cite{fe_ka}, we will use the fact that for any $\ell\in\mathcal{N_0}$, $(\ell_1,\ell_2,\ell_3)\in\mathcal{I}_{\ell}^{N_0}$ and $n\ge 0$ one has
\begin{equation}
    \label{b_fe0}
   (1+|\ell|) |r_{\ell,\ell_1,\ell_2,\ell_3} S_{n,\ell,\ell_1,\ell_2,\ell_3}|\lesssim \tau \frac{|\delta_{\ell,\ell_1,\ell_2,\ell_3}-2\ell_1^2|}{\delta_{\ell,\ell_1,\ell_2,\ell_3}} \lesssim \tau \prod_{j=1}^{3}(1+|\ell_j|).
   \end{equation}
 We start by the term $\mathcal{I}_5^n$. By the bound in \eqref{b_fe0} and Assumption \ref{smal_cub_ass}, one gets
   \begin{align*}
       \mathbb{E}\big[\max_{1\le n\le N_\varepsilon}\|\mathcal{I}_5^{n-1}\|^{2}_{H^1}\big]& = \varepsilon^4\mathbb{E}\big[\max_{1\le n\le N_\varepsilon}\sum_{\ell\in\mathcal{T}_{N_0}} (1+|\ell|^2)\big|
      \sum_{(\ell_1,\ell_2,\ell_3)\in\mathcal{I}_{\ell}^{N_0}} S_{n-1,\ell,\ell_1,\ell_2,\ell_3}r_{\ell,\ell_1,\ell_2,\ell_3}\\
      &\hspace{2cm}\times \theta_R(\eta^{\varepsilon,R}(t_{n-1}))(\widehat{\eta^{\varepsilon,R}}_{\ell_1}(t_{n-1}))^*\widehat{\eta^{\varepsilon,R}}_{\ell_2}(t_{n-1}) \widehat{\eta^{\varepsilon,R}}_{\ell_3}(t_{n-1})
       \big|^2\big]\\
       &\hspace{-1.5cm}\lesssim (\varepsilon^2 \tau)^{2}\mathbb{E}\big[\max_{1\le n\le N_\varepsilon}\sum_{\ell\in\mathcal{T}_{N_0}}\\
       &\times \big(\sum_{(\ell_1,\ell_2,\ell_3)\in\mathcal{I}_{\ell}^{N_0}}
       \big|\widehat{\eta^{\varepsilon,R}}_{\ell_1}(t_{n-1})\big| \big|\widehat{\eta^{\varepsilon,R}}_{\ell_2}(t_{n-1})\big| \big|\widehat{\eta^{\varepsilon,R}}_{\ell_3}(t_{n-1})\big| \prod_{j=1}^{3}(1+|\ell_j|)\big)^2\big]\\
       &\hspace{-1.5cm}\lesssim (\varepsilon^2 \tau)^{2}\mathbb{E}\big[\max_{1\le n\le N_\varepsilon}\|\eta^{\varepsilon,R}(t_{n-1})\|^{6}_{H^2}\big]\\
       &\hspace{-1.5cm}\lesssim (\varepsilon^2 \tau)^{2},
   \end{align*}
   where we used the auxiliary function $\xi^{\varepsilon,R}(x)=\sum_{\ell \in \mathbb{Z}} (1+|\ell|) |\widehat{\eta^{\varepsilon,R}}_{\ell}(t_n)| {\rm e}^{{\bf i}\ell x}$, with $\|\xi^{\varepsilon,R}\|_{H^1}\lesssim \|\eta^{\varepsilon,R}(t_n)\|_{H^2}$, noting that
   $$
   |\xi^{\varepsilon,R}(x)|^2\xi^{\varepsilon,R}(x)=\sum_{\ell\in\mathbb{Z}}\sum_{\ell=-\ell_1+\ell_2+\ell_3}\prod_{j=1}^{3} (1+|\ell_j|)|\widehat{\eta^{\varepsilon,R}}_{\ell_j}(t_n)|{\rm e}^{{\rm i}\ell x}.
   $$
   %The estimates for the other terms in \eqref{b_fe0} are more complicated since using a result of type \eqref{t_reg_cub} does not provide the desires order, i.e., it produces a term of size $\mathcal{O}(\varepsilon \tau^{\frac{1}{2}})$. Note that this introduces further difficulties respect to the deterministic setting addressed in \cite{fe_ka}. To avoid this issue, from \eqref{mild_cub_twis_trunc}, we note that for any $\ell\in\mathbb{Z}$, one has
   We now estimate the other terms in \eqref{start_dec}. From the definition of twisted variable, we note that for any $\ell\in\mathbb{Z}$, one has
   \begin{equation}
       \label{dec0}
   \widehat{\eta^{\varepsilon,R}}_{\ell}(t_k)-\widehat{\eta^{\varepsilon,R}}_{\ell}(t_{k+1})={\bf i}\varepsilon^2 \widehat{\mathcal{A}^{\varepsilon,R}}_{l}(t_k)-\varepsilon \int_{t_k}^{t_{k+1}}{\rm e}^{{\bf i}s \ell^2}q_{\ell}^{\frac{1}{2}}\,d\beta_\ell(s),
     \end{equation}
   where $\widehat{\mathcal{A}^{\varepsilon,R}}_{\ell}(t_k)$ is the $\ell$-th Fourier coefficient of the term
   \begin{equation}
       \label{dec0a}
   \mathcal{A}^{\varepsilon,R}(t_k):=\int_0^\tau {\rm e}^{-{\bf i}(t_k+s)\Delta}\theta_R(\eta^{\varepsilon,R}(t_k+s))\big| {\rm e}^{{\bf i}(t_k+s)\Delta} \eta^{\varepsilon,R}(t_k+s)\big|^2 {\rm e}^{{\bf i}(t_k+s)\Delta} \eta^{\varepsilon,R}(t_k+s)\,ds.
     \end{equation}
   Then, thanks to \eqref{dec0}, we get
   $$
   \mathbb{E}\big[\max_{1\le n\le N_\varepsilon}\|\mathcal{I}_1^{n-1}\|^{2}_{H^1}\big]\lesssim \mathbb{E}\big[\max_{1\le n\le N_\varepsilon}\|\mathcal{I}_{1,1}^{n-1}\|^{2}_{H^1}\big]+\mathbb{E}\big[\max_{1\le n\le N_\varepsilon}\|\mathcal{I}_{1,2}^{n-1}\|^{2}_{H^1}\big],
   $$
   where
   \begin{align*}
       & \mathcal{I}_{1,1}^{n-1}=-\varepsilon^4 \sum_{\ell\in\mathcal{T}_{N_0}}{\rm e}^{{\bf i}\ell x}\sum_{(\ell_1,\ell_2,\ell_3)\in\mathcal{I}_{\ell}^{N_0}}
    r_{\ell,\ell_1,\ell_2,\ell_3}\sum_{k=0}^{n-2}S_{k,\ell,\ell_1,\ell_2,\ell_3}\\
    &\hspace{4cm}\times \theta_R(\eta^{\varepsilon,R}(t_k)) (\widehat{\eta^{\varepsilon,R}}_{\ell_1}(t_k))^* \widehat{\mathcal{A}^{\varepsilon,R}}_{\ell_2} (t_k)\widehat{\eta^{\varepsilon,R}}_{\ell_3}(t_k),\\
    &\mathcal{I}_{1,2}^{n-1} = -{\bf i}\varepsilon^3  \sum_{\ell\in\mathcal{T}_{N_0}}{\rm e}^{{\bf i}\ell x}\sum_{(\ell_1,\ell_2,\ell_3)\in\mathcal{I}_{\ell}^{N_0}}
    r_{\ell,\ell_1,\ell_2,\ell_3}\sum_{k=0}^{n-2}S_{k,\ell,\ell_1,\ell_2,\ell_3}\\
    &\hspace{4cm}\times \theta_R(\eta^{\varepsilon,R}(t_k)) (\widehat{\eta^{\varepsilon,R}}_{\ell_1}(t_k))^* \widehat{\eta^{\varepsilon,R}}_{\ell_3}(t_k)
     \int_{t_k}^{t_{k+1}}{\rm e}^{{\bf i}s \ell_2^2}q_{\ell_2}^{\frac{1}{2}}\,d\beta_{\ell_2}(s).
   \end{align*}
   Thanks to the bound in \eqref{b_fe0}, similar computations as before, Hölder inequality and Assumption \ref{smal_cub_ass}, we get
   \begin{align*}
      \mathbb{E}\big[ \|\mathcal{I}_{1,1}^{n-1} \|^{2}_{H^1}\big]&=\varepsilon^8\mathbb{E}\big[\max_{1\le n\le N_\varepsilon}\|\sum_{k=0}^{n-2}\sum_{\ell\in\mathcal{T}_{N_0}}{\rm e}^{{\bf i}\ell x}\\
      &\hspace{0.4cm}\times\sum_{(\ell_1,\ell_2,\ell_3)\in\mathcal{I}_{\ell}^{N_0}}S_{k,\ell,\ell_1,\ell_2,\ell_3}
    r_{\ell,\ell_1,\ell_2,\ell_3}  \theta_R(\eta^{\varepsilon,R}(t_k)) \widehat{\eta^{\varepsilon,R}}_{\ell_1}(t_k))^* \widehat{\mathcal{A}^{\varepsilon,R}}_{\ell_2} (t_k)\widehat{\eta^{\varepsilon,R}}_{\ell_3}(t_k) \|_{H^1}^{2}\big]\\
    &\lesssim \varepsilon^8 (N_\varepsilon-1)\sum_{k=0}^{N_\varepsilon-2}\mathbb{E}\big[\sum_{\ell\in\mathcal{T}_{N_0}} (1+|\ell|^2)\\
    &\hspace{0.4cm}\big|\sum_{(\ell_1,\ell_2,\ell_3)\in\mathcal{I}_{\ell}^{N_0}}S_{k,\ell,\ell_1,\ell_2,\ell_3}
    r_{\ell,\ell_1,\ell_2,\ell_3}  \theta_R(\eta^{\varepsilon,R}(t_k)) (\widehat{\eta^{\varepsilon,R}}_{\ell_1}(t_k))^* \widehat{\mathcal{A}^{\varepsilon,R}}_{\ell_2} (t_k)\widehat{\eta^{\varepsilon,R}}_{\ell_3}(t_k)\big|^2\big]\\
    &\lesssim (\varepsilon^4 \tau)^{2} (N_\varepsilon-1)\sum_{k=0}^{N_\varepsilon-2}\mathbb{E}\big[\sum_{\ell\in\mathcal{T}_{N_0}}\\
    &\hspace{2cm}\times \big(\sum_{(\ell_1,\ell_2,\ell_3)\in\mathcal{I}_{\ell}^{N_0}}
    \prod_{j=1}^{3}(1+|\ell_j|) \big|\widehat{\eta^{\varepsilon,R}}_{\ell_1}(t_k)| |\widehat{\mathcal{A}^{\varepsilon,R}}_{\ell_2} (t_k)| |\widehat{\eta^{\varepsilon,R}}_{\ell_3}(t_k)|\big)^2\big]\\
    &\lesssim (\varepsilon^4 \tau)^{2} (N_\varepsilon-1)\sum_{k=0}^{N_\varepsilon-2}\mathbb{E}\big[ \|\eta^{\varepsilon,R}(t_k)\|^{4}_{H^2}\|\mathcal{A}^{\varepsilon,R}(t_k)\|^{2}_{H^2}\big]\\
    &\lesssim (\varepsilon^4 \tau)^{2} (N_\varepsilon-1)\sum_{k=0}^{N_\varepsilon-2}\big(\mathbb{E}\big[ \|\eta^{\varepsilon,R}(t_k)\|^{8}_{H^2}\big]\big)^{\frac{1}{2}}\big(\mathbb{E}\big\|\mathcal{A}^{\varepsilon,R}(t_k)\|^{4}_{H^2}\big]\big)^{\frac{1}{2}}\\
    &\lesssim  (\varepsilon^4 \tau)^{2} (N_\varepsilon-1)\sum_{k=0}^{N_\varepsilon-2} \tau^{2}\\
    &\lesssim (\varepsilon^2 \tau)^{2}.
   \end{align*}
   Moreover, we have
   \begin{align*}
        \mathbb{E}\big[ \max_{1\le n\le N_\varepsilon}\|\mathcal{I}_{1,2}^{n-1} \|^2_{H^1}\big] = \varepsilon^6 \mathbb{E}\big[\max_{1\le n\le N_\varepsilon}\|\sum_{\ell\in\mathcal{T}_{N_0}} {\rm e}^{{\bf i}\ell x} \sum_{\ell_2\in\mathcal{T}_{N_0}}\sum_{k=0}^{n-2}\int_{t_k}^{t_{k+1}} G_{\ell,\ell_2}(s)d\beta_{\ell_2}(s)\|^2_{H^1}\big],
   \end{align*}
   where
   $$
   G_{\ell,\ell_2}(s)=\sum_{\substack{\ell_1,\ell_3\in\mathcal{T}_{N_0}\\ \ell-\ell_2=\ell_3-\ell_1\\ \ell^2-\ell_2^2\ne \ell_3^2-\ell_1^2}} r_{\ell,\ell_1,\ell_2,\ell_3}
   S_{[\frac{s}{\tau}],\ell,\ell_1,\ell_2,\ell_3} \theta_R(\eta^{\varepsilon,R}([s/\tau]\tau)) (\widehat{\eta^{\varepsilon,R}}_{\ell_1}([s/\tau]\tau))^* \widehat{\eta^{\varepsilon,R}}_{\ell_3}([s/\tau]\tau)
  {\rm e}^{{\bf i}s \ell_2^2}q_{\ell_2}^{\frac{1}{2}}.
   $$
   Then, we get
   \begin{align*}
       \mathbb{E}\big[\max_{1\le n\le N_\varepsilon} \|\mathcal{I}_{1,2}^{n-1} \|^2_{H^1}\big] &= \varepsilon^6 \mathbb{E}\big[\max_{1\le n\le N_\varepsilon}\sum_{\ell\in\mathcal{T}_{N_0}} (1+|\ell|^2)\big|\sum_{\ell_2\in\mathcal{T}_{N_0}}\int_0^{t_{n-1}} G_{\ell,\ell_2}(s)\,d\beta_{\ell_2}(s) \big|^2\big]\\
       &\le \varepsilon^6 \sum_{\ell\in\mathcal{T}_{N_0}} (1+|\ell|^2)\mathbb{E}\big[\max_{1\le n\le N_\varepsilon}\big|\int_0^{t_{n-1}}  \sum_{\ell_2\in\mathcal{T}_{N_0}}G_{\ell,\ell_2}(s)\,d\beta_{\ell_2}(s) \big|^2\big].
     %  &\lesssim \varepsilon^6 \sum_{\ell\in\mathcal{T}_{N_0}}  \sum_{\ell_2\in\mathcal{T}_{N_0}}\mathbb{E}\big[\int_0^{t_n} |(1+|\ell|)G_{\ell,\ell_2}(s)|^2\,ds\big]\\
      % &\hspace{1cm}+\varepsilon^6 \sum_{\ell\in\mathcal{T}_{N_0}} (1+|\ell|^2) \sum_{\substack{\ell_2,m_2\in\mathcal{T}_{N_0}\\ \ell_2\ne m_2}}\mathbb{E}\big[\int_0^{t_n} G_{\ell,\ell_2}(s)\,d\beta_{\ell_2}(s)\int_0^{t_n} G_{\ell,m_2}(s)\,d\beta_{m_2}(s)\big].
   \end{align*}
   It is convenient to define the auxiliary linear operator
   $$
   G_{\ell}(s):L^2(\mathbb{T})\to\mathbb{C}, \quad G_{\ell}(s)w = \sum_{m\in\mathbb{Z}} G_{\ell,m}(s)\hat{w}_{m},
   $$
   with the notation $w(x)=\frac{1}{\sqrt{2\pi}}\sum_{m\in\mathbb{Z}} {\rm e}^{{\bf i}m x} \hat{w}_m$. 
   Then, one has
   $$
   \int_0^{t_{n-1}} \sum_{\ell_2\in\mathcal{T}_{N_0}}G_{\ell,\ell_2}(s)\,d\beta_{\ell_2}(s) = \int_0^{t_{n-1}} G_{\ell}(s)\mathcal{P}_{N_0} \,dW(s).
   $$
   We see that, by definition, $G_{\ell}(s)$ is independent on later Wiener increments. i.e., $G_{\ell}$ is non-anticipative and, also, $G_{\ell}(s){\rm e}^{{\bf i}\ell_2 x} = \sqrt{2\pi}G_{\ell,\ell_2}(s), \ell_2\in\mathbb{Z}$.
   Then, by BDG inequality, the definition of $G_{\ell,\ell_2}(s)$ and the bound in \eqref{b_fe0}, one has
   \begin{align*}
       \mathbb{E}\big[\max_{1\le n\le N_\varepsilon} \|\mathcal{I}_{1,2}^{n-1} \|^2_{H^1}\big]&\lesssim \varepsilon^6 \sum_{\ell\in\mathcal{T}_{N_0}} (1+|\ell|^2)\mathbb{E}\big[
       \int_0^{T_{\varepsilon}-\tau} \|G_{\ell}(s)\mathcal{P}_{N_0}\|^2_{\mathcal{L}_2(L^2;\mathbb{C})}\,ds\big]\\
       &\lesssim \varepsilon^6 \sum_{k=0}^{N_\varepsilon-2}\int_{t_k}^{t_{k+1}} \mathbb{E}\big[\sum_{\ell\in\mathcal{T}_{N_0}}\sum_{\ell_2\in\mathcal{T}_{N_0}}\big|(1+|\ell|)G_{\ell,\ell_2}(s)\big|^2 \big]\,ds\\
      & \lesssim \varepsilon^6 \tau^2 \sum_{k=0}^{N_\varepsilon-2}\int_{t_k}^{t_{k+1}} \mathbb{E}\big[ \sum_{\ell\in\mathcal{T}_{N_0}} 
      \sum_{\ell_2\in\mathcal{T}_{N_0}} (1+|\ell_2|^2) q_{\ell_2}\\
       &\hspace{0.5cm}\times \big(\sum_{\substack{\ell_1,\ell_3\in\mathcal{T}_{N_0}\\ \ell-\ell_2=\ell_3-\ell_1\\ \ell^2-\ell_2^2\ne \ell_3^2-\ell_1^2}} (1+|\ell_1|)(1+|\ell_3|) \big|\widehat{\eta^{\varepsilon,R}}_{\ell_1}(t_k)\big| |\widehat{\eta^{\varepsilon,R}}_{\ell_3}(t_k)| \big)^2\big]\,ds\\
       &\lesssim \varepsilon^6 \tau^2 \sum_{k=0}^{N_\varepsilon-2}\int_{t_k}^{t_{k+1}} \mathbb{E}\big[ 
      \sum_{\ell_2\in\mathcal{T}_{N_0}} (1+|\ell_2|^2) q_{\ell_2} \\
       &\hspace{0.5cm}\times \sum_{\ell\in\mathcal{T}_{N_0}}  \big(\sum_{\ell_3\in\mathcal{T}_{N_0}} (1+|\ell_3-\ell+\ell_2|)(1+|\ell_3|) \big|\widehat{\eta^{\varepsilon,R}}_{\ell_3-\ell+\ell_2}(t_k)\big| \big|\widehat{\eta^{\varepsilon,R}}_{\ell_3}(t_k)\big| \big)^2\big]\,ds.
   \end{align*}
   Introducing the auxiliary function $\tilde{\eta}^{\varepsilon,R}(x)=\sum_{\ell \in \mathbb{Z}} \widehat{\tilde{\eta}^{\varepsilon,R}}_\ell {\rm e}^{{\bf i} \ell x}$ with $\widehat{\tilde{\eta}^{\varepsilon,R}}_\ell=\widehat{\eta^{\varepsilon,R}}_{-\ell+\ell_2}$, using the Young's inequality and then the Hölder inequality, we get
   \begin{align*}
       &\sum_{\ell_2\in\mathcal{T}_{N_0}} (1+|\ell_2|^2) q_{\ell_2}\sum_{\ell\in\mathcal{T}_{N_0}}  \big(\sum_{\ell_3\in\mathcal{T}_{N_0}} (1+|\ell_3-\ell+\ell_2|)(1+|\ell_3|) |\widehat{\eta^{\varepsilon,R}}_{\ell_3-\ell+\ell_2}(t_k)| |\widehat{\eta^{\varepsilon,R}}_{\ell_3}(t_k)| \big)^2\\
       &\quad = \sum_{\ell_2\in\mathcal{T}_{N_0}} (1+|\ell_2|^2) q_{\ell_2}\sum_{\ell\in\mathcal{T}_{N_0}}  \big(\sum_{\ell_3\in\mathcal{T}_{N_0}} (1+|\ell_3-\ell+\ell_2|)(1+|\ell_3|) |\widehat{\tilde{\eta}^{\varepsilon,R}}_{\ell-\ell_3}(t_k)| |\widehat{\eta^{\varepsilon,R}}_{\ell_3}(t_k)| \big)^2\\
       &\quad \lesssim \sum_{\ell_2\in\mathcal{T}_{N_0}} (1+|\ell_2|^2) q_{\ell_2} \|(\tilde{\eta}^{\varepsilon,R})^{(1)}(t_k)\star (\eta^{\varepsilon,R})^{(1)}(t_k)\|^2_{\ell_2}\\
       &\hspace{5cm} +  \sum_{\ell_2\in\mathcal{T}_{N_0}} (1+|\ell_2|^4) q_{\ell_2} \|\tilde{\eta}^{\varepsilon,R}(t_k)\star (\eta^{\varepsilon,R})^{(1)}(t_k)\|^2_{\ell_2}\\
       &\quad \lesssim \|\eta^{\varepsilon,R}(t_k)\|^2_{H^1}\sum_{\ell_2\in\mathcal{T}_{N_0}} (1+|\ell_2|^2) q_{\ell_2} \|(\tilde{\eta}^{\varepsilon,R})^{(1)}(t_k)\|^2_{\ell_1}\\
       &\hspace{5cm} +  \sum_{\ell_2\in\mathcal{T}_{N_0}} (1+|\ell_2|^4) q_{\ell_2} \|(\tilde{\eta}^{\varepsilon,R})^{(1)}(t_k)\|^2_{\ell_2} \|(\eta^{\varepsilon,R})^{(1)}(t_k)\|^2_{\ell_1}\\
       &\quad \lesssim \|\eta^{\varepsilon,R}(t_k)\|^2_{H^1}\sum_{\ell_2\in\mathcal{T}_{N_0}} (1+|\ell_2|^2) q_{\ell_2} \big(\sum_{\ell} (1+|\ell|) |\widehat{\eta^{\varepsilon,R}}_{-\ell+\ell_2}|\big)^2\\
       &\hspace{5cm}+\|\eta^{\varepsilon,R}(t_k)\|^2_{H^{\frac{3}{2}+}} \sum_{\ell_2\in\mathcal{T}_{N_0}} (1+|\ell_2|^4) q_{\ell_2} \sum_{\ell} |\widehat{\eta^{\varepsilon,R}}_{-\ell+\ell_2}|^2\\
       &\quad \lesssim \|\eta^{\varepsilon,R}(t_k)\|^2_{H^1} \sum_{\ell_2\in\mathcal{T}_{N_0}} (1+|\ell_2|^2) q_{\ell_2} \big(\sum_{\ell} (1+|\ell-\ell_2|) |\widehat{\eta^{\varepsilon,R}}_{-\ell+\ell_2}|\big)^2\\
       &\hspace{5cm}+ \|\eta^{\varepsilon,R}(t_k)\|^2_{H^1} \sum_{\ell_2\in\mathcal{T}_{N_0}} (1+|\ell_2|^4) q_{\ell_2} \big(\sum_{\ell}  |\widehat{\eta^{\varepsilon,R}}_{-\ell+\ell_2}|\big)^2\\
       &\hspace{5cm}+\|\eta^{\varepsilon,R}(t_k)\|^2_{H^{\frac{3}{2}+}} \|\eta^{\varepsilon,R}(t_k)\|^2\|Q^{\frac{1}{2}}\|^2_{\mathcal{L}_2^2}\\
       &\quad \lesssim \|\eta^{\varepsilon,R}(t_k)\|^4_{H^{\frac{3}{2}+}}\|Q^{\frac{1}{2}}\|^2_{\mathcal{L}_2^2}.
   \end{align*}
   Then, we get
   \begin{align*}
       %\varepsilon^6 \sum_{\ell\in\mathcal{T}_{N_0}}  \sum_{\ell_2\in\mathcal{T}_{N_0}}\mathbb{E}\big[\int_0^{t_n} |(1+|\ell|)G_{\ell,\ell_2}(s)|^2\,ds\big]
       %&\lesssim \varepsilon^6\tau^2 \|Q^{\frac{1}{2}}\|^2_{\mathcal{L}_2^2} \sum_{k=0}^{n-1}\int_{t_k}^{t_{k+1}} \mathbb{E}\big[\|v^R(t_k)\|_{H^{\frac{3}{2}+}}^{4}\big]\,ds\\
       %&\lesssim\varepsilon^6 \tau^2 t_n\\
       %&\lesssim(\varepsilon^2 \tau)^2.
        \mathbb{E}\big[\max_{1\le n\le N_\varepsilon} \|\mathcal{I}_{1,2}^{n-1} \|^2_{H^1}\big]&\lesssim\varepsilon^6 \tau^2 \|Q^{\frac{1}{2}}\|^2_{\mathcal{L}_2^2}\sum_{k=0}^{N_\varepsilon-2} \int_{t_k}^{t_{k+1}}\mathbb{E}\big[\|\eta^{\varepsilon,R}(t_k)\|^4_{H^{\frac{3}{2}+}}\big]\,ds\\
        &\lesssim \varepsilon^6 \tau^2 N_\varepsilon \tau\\
        &\lesssim(\varepsilon^2\tau)^2,
   \end{align*}
   %Moreover, defining
  % $$
   %G_{\ell,\ell_2}^k=\sum_{\substack{\ell_1,\ell_3\in\mathcal{T}_{N_0}\\ \ell-\ell_2=\ell_3-\ell_1\\ \ell^2-\ell_2^2\ne \ell_3^2-\ell_1^2}} r_{\ell,\ell_1,\ell_2,\ell_3}
   %S_{k,\ell,\ell_1,\ell_2,\ell_3} \theta_R(v^R(t_k)) \overline{\hat{v}^R_{\ell_1}}(t_k) \hat{v}^R_{\ell_3}(t_k), 
   %$$
   %we have
   %\begin{align*}
    %  & \sum_{\substack{\ell_2,m_2\in\mathcal{T}_{N_0}\\ \ell_2\ne m_2}}\mathbb{E}\big[\int_0^{t_n} G_{\ell,\ell_2}(s)\,d\beta_{\ell_2}(s)\int_0^{t_n} G_{\ell,m_2}(s)\,d\beta_{m_2}(s)\big]\\
     %  &\quad = \sum_{\substack{\ell_2,m_2\in\mathcal{T}_{N_0}\\ \ell_2\ne m_2}}\mathbb{E}\big[\sum_{k,j=0}^{n-1} G_{\ell,\ell_2}^k G_{\ell,m_2}^j \int_{t_k}^{t_{k+1}} {\rm e}^{{\bf i}s\ell_2^2}q_{\ell_2}^{\frac{1}{2}}\,d\beta_{\ell_2}(s)   \int_{t_j}^{t_{j+1}} {\rm e}^{{\bf i}sm_2^2}q_{m_2}^{\frac{1}{2}}\,d\beta_{m_2}(s) \big]\\
     %  &\quad =0,
   %\end{align*}
   %since, for the case $k=j$, the random variable $G_{\ell,\ell_2}^k G_{\ell,m_2}^k$ depends on values of the exact solution up to time $t_k$ and the stochastic integrals depend on Brownian increments later than $t_k$, while for the case $k<j$, the random variable $G_{\ell,\ell_2}^k G_{\ell,m_2}^j \int_{t_k}^{t_{k+1}} {\rm e}^{{\bf i}s l_2^2} q_{\ell_2}^{\frac{1}{2}}\,d\beta_{\ell_2}(s)$ is independent from the stochastic integral $\int_{t_j}^{t_{j+1}} {\rm e}^{{\bf i}s m_2^2} q_{m_2}^{\frac{1}{2}}\,d\beta_{m_2}(s)$ by similar arguments than above. Similarly, the same argument can apply for the case $k>j$.
   %This finally shows
   %$$
   %\mathbb{E}\big[\|\mathcal{I}_{1,2}^n\|_{H^1}^2\big] \lesssim (\varepsilon^2 \tau)^2
   %$$
   that finally gives
   $$
    \mathbb{E}\big[\max_{1\le n\le N_\varepsilon}\|\mathcal{I}_{1}^{n-1}\|_{H^1}^2\big] \lesssim (\varepsilon^2 \tau)^2.
   $$
   With similar technique, we can bound the other terms $\mathcal{I}_{2}^{n-1}, \mathcal{I}_3^{n-1}$ and $\mathcal{I}_4^{n-1}$. We have
   $$
   \mathbb{E}\big[\max_{1\le n\le N_\varepsilon}\|I_2^{n-1}\|^2_{H^1}\big]\lesssim \mathbb{E}\big[\max_{1\le n\le N_\varepsilon}\|I_{2,1}^{n-1}\|^2_{H^1}\big]+\mathbb{E}\big[\max_{1\le n\le N_\varepsilon}\|I_{2,2}^{n-1}\|^2_{H^1}\big],
   $$
   where
   \begin{align*}
      &\mathcal{I}_{2,1}^{n-1}:={\bf i}\varepsilon^2 \sum_{\ell\in\mathcal{T}_{N_0}}{\rm e}^{{\bf i}\ell x}\sum_{(\ell_1,\ell_2,\ell_3)\in\mathcal{I}_{\ell}^{N_0}}
    r_{\ell,\ell_1,\ell_2,\ell_3}\sum_{k=0}^{n-2}S_{k,\ell,\ell_1,\ell_2,\ell_3}\\
    &\hspace{2cm}\times \theta_R(\eta^{\varepsilon,R}(t_k)) (\widehat{\eta^{\varepsilon,R}}_{\ell_1}(t_k)-\widehat{\eta^{\varepsilon,R}}_{\ell_1}(t_{k+1}))^*(\widehat{\eta^{\varepsilon,R}}_{\ell_2}(t_{k+1})-\widehat{\eta^{\varepsilon,R}}_{\ell_2}(t_{k}))\widehat{\eta^{\varepsilon,R}}_{\ell_3}(t_k),\\
    &\mathcal{I}_{2,2}^{n-1}:={\bf i}\varepsilon^2 \sum_{\ell\in\mathcal{T}_{N_0}}{\rm e}^{{\bf i}\ell x}\sum_{(\ell_1,\ell_2,\ell_3)\in\mathcal{I}_{\ell}^{N_0}}
    r_{\ell,\ell_1,\ell_2,\ell_3}\sum_{k=0}^{n-2}S_{k,\ell,\ell_1,\ell_2,\ell_3}\\
    &\hspace{2cm}\times \theta_R(\eta^{\varepsilon,R}(t_k)) (\widehat{\eta^{\varepsilon,R}}_{\ell_1}(t_k)-\widehat{\eta^{\varepsilon,R}}_{\ell_1}(t_{k+1}))^*\widehat{\eta^{\varepsilon,R}}_{\ell_2}(t_k)\widehat{\eta^{\varepsilon,R}}_{\ell_3}(t_k).
   \end{align*}
   The term $\mathbb{E}\big[\displaystyle\max_{1\le n\le N_\varepsilon}\|I_{2,2}^{n-1}\|^2_{H^1}\big]$ can be bounded in the same way of the term $\mathbb{E}\big[\displaystyle\max_{1\le n\le N_\varepsilon}\|I_{1}^{n-1}\|^2_{H^1}\big]$ and then we have
   $$
   \mathbb{E}\big[\max_{1\le n\le N_\varepsilon}\|I_{2,2}^{n-1}\|^2_{H^1}\big]\lesssim (\varepsilon^2 \tau)^2.
   $$
   For the term $\mathbb{E}\big[\displaystyle\max_{1\le n\le N_\varepsilon}\|I_{2,1}^{n-1}\|^2_{H^1}\big]$, as done for the term $\mathcal{I}_1$, using the bound in \eqref{b_fe0}, the algebra property of $H^2$, the Hölder inequality, the estimate in \ref{t_reg_cub} and Assumption \ref{smal_cub_ass}, we have
   \begin{align*}
       \mathbb{E}\big[\max_{1\le n\le N_\varepsilon}\|I_{2,1}^{n-1}\|^2_{H^1}\big]& \lesssim \varepsilon^4 (N_\varepsilon-1)\sum_{k=0}^{N_\varepsilon-2}\mathbb{E}\big[\sum_{\ell\in\mathcal{T}_{N_0}}(1+|\ell|^2)\big(\sum_{(\ell_1,\ell_2,\ell_3)\in\mathcal{I}_{\ell}^{N_0}} \big|r_{\ell,\ell_1,\ell_2,\ell_3} S_{k,\ell,\ell_1,\ell_2,\ell_3}\big|\\
       &\hspace{1cm}\times \big|\widehat{\eta^{\varepsilon,R}}_{\ell_1}(t_k)-\widehat{\eta^{\varepsilon,R}}_{\ell_1}(t_{k+1})\big| \big|\widehat{\eta^{\varepsilon,R}}_{\ell_2}(t_{k+1})-\widehat{\eta^{\varepsilon,R}}_{\ell_2}(t_{k})\big| \big|\widehat{\eta^{\varepsilon,R}}_{\ell_3}(t_k)\big|\big)^2\big]\\
       &\lesssim  \varepsilon^4 \tau^2 (N_\varepsilon-1)\sum_{k=0}^{N_\varepsilon-2}\mathbb{E}\big[\big(\sum_{(\ell_1,\ell_2,\ell_3)\in\mathcal{I}_{\ell}^{N_0}} \prod_{j=1}^{3}(1+|\ell_j|)\\
       &\hspace{1cm}\times \big|\widehat{\eta^{\varepsilon,R}}_{\ell_1}(t_k)-\widehat{\eta^{\varepsilon,R}}_{\ell_1}(t_{k+1})\big| \big|\widehat{\eta^{\varepsilon,R}}_{\ell_2}(t_{k+1})-\widehat{\eta^{\varepsilon,R}}_{\ell_2}(t_{k})\big| \big|\widehat{\eta^{\varepsilon,R}}_{\ell_3}(t_k)\big|\big)^2\big]\\
       &\lesssim \varepsilon^4 \tau^2 (N_\varepsilon-1)\sum_{k=0}^{N_\varepsilon-2}\mathbb{E}\big[\|\eta^{\varepsilon,R}(t_{k+1})-\eta^{\varepsilon,R}(t_k)\|^4_{H^2} \|\eta^{\varepsilon,R}(t_k)\|_{H^2}^2\big]\\
       &\lesssim \varepsilon^4 \tau^2 (N_\varepsilon-1)\sum_{k=0}^{N_\varepsilon-2}\big(\mathbb{E}\big[\|\eta^{\varepsilon,R}(t_{k+1})-\eta^{\varepsilon,R}(t_k)\|^8_{H^2}\big]\big)^{\frac{1}{2}}\big(\mathbb{E}\big[ \|\eta^{\varepsilon,R}(t_k)\|_{H^2}^4\big]\big)^{\frac{1}{2}}\\
       &\lesssim \varepsilon^4 \tau^2 N_\varepsilon^2 (\varepsilon \tau^{\frac{1}{2}})^4\\
       &\lesssim(\varepsilon^2 \tau)^2.
   \end{align*}
   Then, we finally get
   $$
   \mathbb{E}\big[\max_{1\le n\le N_\varepsilon}\|I_{2}^{n-1}\|^2_{H^1}\big]\lesssim (\varepsilon^2 \tau)^2.
   $$
   With a similar technique, we can get the same bound for the term $\mathbb{E}\big[\displaystyle\max_{1\le n\le N_\varepsilon}\|I_{3}^{n-1}\|^2_{H^1}\big]$. We leave the computation to the reader.

   It is remained to bound the term $\mathbb{E}\big[\displaystyle\max_{1\le n\le N_\varepsilon}\|I_{4}^{n-1}\|^2_{H^1}\big]$. This term may be decomposed as
   \begin{align*}
   \mathbb{E}\big[\max_{1\le n\le N_\varepsilon}\|I_{4}^{n-1}\|^2_{H^1}\big]&\lesssim \mathbb{E}\big[\max_{1\le n\le N_\varepsilon}\|I_{4,1}^{n-1}\|^2_{H^1}\big]+\mathbb{E}\big[\max_{1\le n\le N_\varepsilon}\|I_{4,2}^{n-1}\|^2_{H^1}\big]\\
   &\quad +\mathbb{E}\big[\max_{1\le n\le N_\varepsilon}\|I_{4,3}^{n-1}\|^2_{H^1}\big]+\mathbb{E}\big[\max_{1\le n \le N_\varepsilon}\|I_{4,4}^{n-1}\|^2_{H^1}\big],
   \end{align*}
   where 
   \begin{align*}
       &\mathcal{I}_{4,1}^{n-1}:={\bf i}\varepsilon^2 \sum_{\ell\in\mathcal{T}_{N_0}}{\rm e}^{{\bf i}\ell x}\sum_{(\ell_1,\ell_2,\ell_3)\in\mathcal{I}_{\ell}^{N_0}}
    r_{\ell,\ell_1,\ell_2,\ell_3}\sum_{k=0}^{n-2}S_{k,\ell,\ell_1,\ell_2,\ell_3}\\
    &\hspace{0.5cm}\times \big(\widehat{\eta^{\varepsilon,R}}_{\ell_1}(t_{k+1})-\widehat{\eta^{\varepsilon,R}}_{\ell_1}(t_{k})\big)^*\widehat{\eta^{\varepsilon,R}}_{\ell_2}(t_{k+1})\widehat{\eta^{\varepsilon,R}}_{\ell_3}(t_{k+1})(\theta_R(\eta^{\varepsilon,R}(t_k))-\theta_R(\eta^{\varepsilon,R}(t_{k+1}))),\\ 
    &\mathcal{I}_{4,2}^{n-1}:={\bf i}\varepsilon^2 \sum_{\ell\in\mathcal{T}_{N_0}}{\rm e}^{{\bf i}\ell x}\sum_{(\ell_1,\ell_2,\ell_3)\in\mathcal{I}_{\ell}^{N_0}}
    r_{\ell,\ell_1,\ell_2,\ell_3}\sum_{k=0}^{n-2}S_{k,\ell,\ell_1,\ell_2,\ell_3}\\
    &\hspace{0.5cm}\times \big(\widehat{\eta^{\varepsilon,R}}_{\ell_1}(t_{k})\big)^*\big(\widehat{\eta^{\varepsilon,R}}_{\ell_2}(t_{k+1})-\widehat{\eta^{\varepsilon,R}}_{\ell_2}(t_{k})\big)\widehat{\eta^{\varepsilon,R}}_{\ell_3}(t_{k+1})(\theta_R(\eta^{\varepsilon,R}(t_k))-\theta_R(\eta^{\varepsilon,R}(t_{k+1}))),\\ 
    &\mathcal{I}_{4,3}^{n-1}:={\bf i}\varepsilon^2 \sum_{\ell\in\mathcal{T}_{N_0}}{\rm e}^{{\bf i}\ell x}\sum_{(\ell_1,\ell_2,\ell_3)\in\mathcal{I}_{\ell}^{N_0}}
    r_{\ell,\ell_1,\ell_2,\ell_3}\sum_{k=0}^{n-2}S_{k,\ell,\ell_1,\ell_2,\ell_3}\\
    &\hspace{0.5cm}\times \big(\widehat{\eta^{\varepsilon,R}}_{\ell_1}(t_{k})\big)^*\widehat{\eta^{\varepsilon,R}}_{\ell_2}(t_{k})\big(\widehat{\eta^{\varepsilon,R}}_{\ell_3}(t_{k+1})-\widehat{\eta^{\varepsilon,R}}_{\ell_3}(t_{k})\big)(\theta_R(\eta^{ \varepsilon,R}(t_k))-\theta_R(\eta^{\varepsilon,R}(t_{k+1}))),\\ 
    &\mathcal{I}_{4,4}^{n-1}:={\bf i}\varepsilon^2 \sum_{\ell\in\mathcal{T}_{N_0}}{\rm e}^{{\bf i}\ell x}\sum_{(\ell_1,\ell_2,\ell_3)\in\mathcal{I}_{\ell}^{N_0}}
    r_{\ell,\ell_1,\ell_2,\ell_3}\sum_{k=0}^{n-2}S_{k,\ell,\ell_1,\ell_2,\ell_3}\\
    &\hspace{0.5cm}\times \big(\widehat{\eta^{\varepsilon,R}}_{\ell_1}(t_{k})\big)^*\widehat{\eta^{\varepsilon,R}}_{\ell_2}(t_{k})\widehat{\eta^{\varepsilon,R}}_{\ell_3}(t_{k})(\theta_R(\eta^{\varepsilon,R}(t_k))-\theta_R(\eta^{\varepsilon,R}(t_{k+1}))).
   \end{align*}
   The estimation of the first three above terms is similar. So, here we give only the estimation of the first and of the last of the above terms. With similar computations as above and taking into account the fact that the map $\theta_R: H^1\to \mathbb{R}^+$ is globally Lipschitz, we have
   \begin{align*}
       &\mathbb{E}\big[\max_{1\le n\le N_\varepsilon}\|I_{4,1}^{n-1}\|^2_{H^1}\big]\lesssim\varepsilon^4 (N_\varepsilon-1) \sum_{k=0}^{N_\varepsilon-2}\mathbb{E}\big[\big|\theta_R(\eta^{\varepsilon,R}(t_k))-\theta_R(\eta^{\varepsilon,R}(t_{k+1}))\big|^2\sum_{\ell\in\mathcal{T}_{N_0}} (1+|\ell|^2\\
       &\hspace{1cm}\times \big|\sum_{(\ell_1,\ell_2,\ell_3)\in\mathcal{I}_{\ell}^{N_0}}r_{\ell,\ell_1,\ell_2,\ell_3}S_{k,\ell,\ell_1,\ell_2,\ell_3}
       \big(\widehat{\eta^{\varepsilon,R}}_{\ell_1}(t_{k})-\widehat{\eta^{\varepsilon,R}}_{\ell_1}(t_{k+1})\big)^*\widehat{\eta^{\varepsilon,R}}_{\ell_2}(t_{k+1})\widehat{\eta^{\varepsilon,R}}_{\ell_3}(t_{k+1})\big|^2\big]\\
       &\lesssim \varepsilon^4 \tau^2 (N_\varepsilon-1)\sum_{k=0}^{N_\varepsilon-2}\mathbb{E}\big[\big|\theta_R(\eta^{\varepsilon,R}(t_k))-\theta_R(\eta^{\varepsilon,R}(t_{k+1}))\big|^2\sum_{\mathcal{T}_{N_0}}\big(\sum_{(\ell_1,\ell_2,\ell_3)\in\mathcal{I}_{\ell}^{N_0}}\prod_{j=1}^{3}(1+|\ell_j|)\\
       &\hspace{1cm}\times \big|\widehat{\eta^{\varepsilon,R}}_{\ell_1}(t_{k})-\widehat{\eta^{\varepsilon,R}}_{\ell_1}(t_{k+1})\big| \big|\widehat{\eta^{\varepsilon,R}}_{\ell_2}(t_{k+1})\big| \big|\widehat{\eta^{\varepsilon,R}}_{\ell_3}(t_{k+1})\big| \big)^2\big]\\
       &\lesssim  \varepsilon^4 \tau^2 (N_\varepsilon-1)\sum_{k=0}^{N_\varepsilon-2}\mathbb{E}\big[\big|\theta_R(\eta^{\varepsilon,R}(t_k))-\theta_R(\eta^{\varepsilon,R}(t_{k+1}))\big|^2 \|\eta^{\varepsilon,R}(t_{k+1})-\eta^{\varepsilon,R}(t_{k})\|_{H^2}^2 \|\eta^{\varepsilon,R}(t_{k+1})\|^4_{H^2}\big]\\
       &\lesssim  \varepsilon^4 \tau^2 (N_\varepsilon-1)\sum_{k=0}^{N_\varepsilon-2}\mathbb{E}\big[\|\eta^{\varepsilon,R}(t_{k+1})-\eta^{\varepsilon,R}(t_{k})\|_{H^2}^4 \|\eta^{\varepsilon,R}(t_{k+1})\|^4_{H^2}\big]\\
       &\lesssim \varepsilon^4 \tau^2 (N_\varepsilon-1)\sum_{k=0}^{N_\varepsilon-2}\big(\mathbb{E}\big[\|\eta^{\varepsilon,R}(t_{k+1})-\eta^{\varepsilon,R}(t_{k})\|_{H^2}^8\big)^{\frac{1}{2}} \big]\\
       &\lesssim (\varepsilon^2 \tau)^2,
   \end{align*}
   and similarly for the term $\mathcal{I}_{4,2}^{n-1}$ and $\mathcal{I}^{n-1}_{4,3}$. It now remains to estimate the term $\mathbb{E}\big[\displaystyle\max_{1\le n \le N_\varepsilon}\|I_{4,4}^{n-1}\|_{H^1}^2\big]$. Similarly to the proof of Lemma \ref{strong_lem_trunc}, we can split this term into
   \begin{align*}
   \mathbb{E}\big[\max_{1\le n \le N_\varepsilon}\|\mathcal{I}^{n-1}_{4,4}\|^2_{H^1}\big]&\lesssim  \mathbb{E}\big[\max_{1\le n \le N_\varepsilon}\|\mathcal{I}^{n-1}_{4,4,a}\|^2_{H^1}\big]+ \mathbb{E}\big[\max_{1\le n \le N_\varepsilon}\|\mathcal{I}^{n-1}_{4,4,b}\|^2_{H^1}\big]\\
   &\quad + \mathbb{E}\big[\max_{1\le n \le N_\varepsilon}\|\mathcal{I}^{n-1}_{4,4,c}\|^2_{H^1}\big],
   \end{align*}
   where 
   \begin{align*}
       &\mathbb{E}\big[\max_{1\le n\le N_\varepsilon}\|\mathcal{I}^{n-1}_{4,4,a}\|^2_{H^1}\big] = \varepsilon^8\mathbb{E}\big[\max_{1\le n \le N_\varepsilon}\| \sum_{\ell\in\mathcal{T}_{N_0}}{\rm e}^{{\bf i}\ell x}\sum_{(\ell_1,\ell_2,\ell_3)\in\mathcal{I}_{\ell}^{N_0}}
    r_{\ell,\ell_1,\ell_2,\ell_3}\sum_{k=0}^{n-2}S_{k,\ell,\ell_1,\ell_2,\ell_3}\\
    &\hspace{1cm}\times \big(\widehat{\eta^{\varepsilon,R}}_{\ell_1}(t_{k})\big)^*\widehat{\eta^{\varepsilon,R}}_{\ell_2}(t_{k})\widehat{\eta^{\varepsilon,R}}_{\ell_3}(t_{k})D\theta_R(\eta^{\varepsilon,R}(t_k))\mathcal{A}^{\varepsilon,R}(t_k)\|^2_{H^1}\big],\\
     &\mathbb{E}\big[\max_{1\le n \le N_\varepsilon}\|\mathcal{I}^{n-1}_{4,4,b}\|^2_{H^1}\big] = \varepsilon^6 \mathbb{E}\big[\max_{1\le n\le N_\varepsilon}\|\sum_{\ell\in\mathcal{T}_{N_0}}{\rm e}^{{\bf i}\ell x}\sum_{(\ell_1,\ell_2,\ell_3)\in\mathcal{I}_{\ell}^{N_0}}
    r_{\ell,\ell_1,\ell_2,\ell_3}\sum_{k=0}^{n-2}S_{k,\ell,\ell_1,\ell_2,\ell_3}\\
    &\hspace{1cm}\times \big(\widehat{\eta^{\varepsilon,R}}_{\ell_1}(t_{k})\big)^*\widehat{\eta^{\varepsilon,R}}_{\ell_2}(t_{k})\widehat{\eta^{\varepsilon,R}}_{\ell_3}(t_{k})D\theta_R(\eta^{\varepsilon,R}(t_k))\int_{t_k}^{t_{k+1}}{\rm e}^{-{\bf i}s\Delta}Q^{\frac{1}{2}}\,dW(s)\|^2_{H^1}\big],\\
    &\mathbb{E}\big[\max_{1\le n\le N_\varepsilon}\|\mathcal{I}^{n-1}_{4,4,c}\|^2_{H^1}\big]\\
    &\quad = \varepsilon^4 \mathbb{E}\big[\max_{1\le n\le N_\varepsilon}\|\sum_{\ell\in\mathcal{T}_{N_0}}{\rm e}^{{\bf i}\ell x}\sum_{(\ell_1,\ell_2,\ell_3)\in\mathcal{I}_{\ell}^{N_0}}
    r_{\ell,\ell_1,\ell_2,\ell_3}\sum_{k=0}^{n-2}S_{k,\ell,\ell_1,\ell_2,\ell_3} \big(\widehat{\eta^{\varepsilon,R}}_{\ell_1}(t_{k})\big)^*\widehat{\eta^{\varepsilon,R}}_{\ell_2}(t_{k})\widehat{\eta^{\varepsilon,R}}_{\ell_3}(t_{k})\\
    &\hspace{0.1cm}\times \int_0^1 D^2\theta_R(\mu \eta^{\varepsilon,R}(t_{k+1})+(1-\mu)\eta^{\varepsilon,R}(t_k))
    \big(\eta^{\varepsilon,R}(t_{k+1})-\eta^{\varepsilon,R}(t_k),\eta^{\varepsilon,R}(t_{k+1})-\eta^{\varepsilon,R}(t_k)\big)\,d\mu\|^2_{H^1}\big],
   \end{align*}
   with $\mathcal{A}^{\varepsilon,R}(t_k)$ defined in \eqref{dec0a}. By similar computation as in the proof of Lemma \ref{strong_lem_trunc}, due to the fact that $\theta_R(v)=0$ for $\|v\|_{H^1}\ge 2R$, for any $k$ one has 
   \begin{equation}
       \label{a1}
       \|\mathcal{A}^{\varepsilon,R}(t_k)\|_{H^1}\lesssim \tau, \qquad \mathbb{P}-{\rm a.s.}
   \end{equation}
   Then, using similar bounds shown in previous estimates, the uniform boundedness of $D\theta_R$, the Hölder inequality, Assumption \ref{smal_cub_ass} and the estimate in \eqref{a1} , one gets
   \begin{align*}
       \mathbb{E}\big[\max_{1\le n\le N_\varepsilon}\|\mathcal{I}^{n-1}_{4,4,a}\|^2_{H^1}\big]&\lesssim \varepsilon^8 \tau^2 (N_\varepsilon-1) \sum_{k=0}^{N_\varepsilon-2} \mathbb{E}\big[|D\theta_R(\eta^{\varepsilon,R}(t_k)) \mathcal{A}^{\varepsilon,R}(t_k)|^2 \|\eta^{\varepsilon,R}(t_k)\|^6_{H^2}\big]\\
       &\lesssim \varepsilon^8 \tau^2 (N_\varepsilon-1)  \sum_{k=0}^{N_\varepsilon-2} \mathbb{E}\big[\|\mathcal{A}^{\varepsilon,R}(t_k)\|^2_{H^1} \|\eta^{\varepsilon,R}(t_k)\|^6_{H^2} \big]\\
       &\lesssim \varepsilon^8 \tau^2 (N_\varepsilon-1)  \sum_{k=0}^{N_\varepsilon-2} \big(\mathbb{E}\big[\|\mathcal{A}^{\varepsilon,R}(t_k)\|^4_{H^1}\big]\big)^{\frac{1}{2}} \\
       &\lesssim \varepsilon^8 \tau^4 (N_\varepsilon-1)^2\\
       &\lesssim (\varepsilon^2\tau)^2.
   \end{align*}
   similarly, using the uniform boundedness of $D^2\theta_R$ (see, e.g.,\cite{deb0}) and also \eqref{t_reg_cub}, one get
   \begin{align*}
       \mathbb{E}\big[\max_{1\le n\le N_\varepsilon}\|\mathcal{I}^{n-1}_{4,4,c}\|^2_{H^1}&\lesssim \varepsilon^8 \tau^2 (N_\varepsilon-1) \sum_{k=0}^{N_\varepsilon-2} \mathbb{E}\big[\|\eta^{\varepsilon,R}(t_{k+1})-\eta^{\varepsilon,R}(t_k)\|_{H^1}^4 \|\eta^{\varepsilon,R}(t_k)\|^6_{H^2}\big]\\
       &\lesssim \varepsilon^8 \tau^2 (N_\varepsilon-1) \sum_{k=0}^{N_\varepsilon-2} \big(\mathbb{E}\big[\|\eta^{\varepsilon,R}(t_{k+1})-\eta^{\varepsilon,R}(t_k)\|_{H^1}^8\big]\big)^{\frac{1}{2}}\\
       &\lesssim   \varepsilon^8 \tau^2 (N_\varepsilon-1)^2 (\varepsilon \tau^{\frac{1}{2}})^4\\
       &\lesssim (\varepsilon^2 \tau)^2.
   \end{align*}
   The last term to estimate is $\mathbb{E}\big[\displaystyle\max_{1\le n\le N_\varepsilon}\|\mathcal{I}^{n-1}_{4,4,b}\|^2_{H^1}\big]$. The BDG inequality here yields
   \begin{align*}
      & \mathbb{E}\big[\max_{1\le n\le N_\varepsilon}\|\mathcal{I}^{n-1}_{4,4,b}\|^2_{H^1}\big]=\varepsilon^6 \mathbb{E}\big[\max_{1\le n\le N_\varepsilon}\|\sum_{\ell\in \mathcal{T}_{N_0}}{\rm e}^{{\bf i}\ell x} \sum_{k=0}^{n-2} \int_{t_k}^{t_{k+1}} 
       \sum_{(\ell_1,\ell_2,\ell_3)\in\mathcal{I}_{\ell}^{N_0}} r_{\ell,\ell_1,\ell_2,\ell_3}S_{[\frac{s}{\tau}],\ell,\ell_1,\ell_2,\ell_3}\\
       &\hspace{1cm}  \times \big(\widehat{\eta^{\varepsilon,R}}_{\ell_1}([s/\tau]\tau)\big)*\widehat{\eta^{\varepsilon,R}}_{\ell_2}([s/\tau]\tau)\widehat{\eta^{\varepsilon,R}}_{\ell_3}([s/\tau]\tau) D\theta_R(\eta^{\varepsilon,R}([s/\tau]\tau)){{\rm e}}^{-{\bf i}s\Delta}Q^{\frac{1}{2}}\,dW(s)\|^2_{H^1}\big]\\
       &= \varepsilon^6 \sum_{\ell\in\mathcal{T}_{N_0}}(1+|\ell|^2)\mathbb{E}\big[\max_{1\le n\le N_\varepsilon}\big|\int_0^{t_{n-1}} \sum_{(\ell_1,\ell_2,\ell_3)\in\mathcal{I}_{\ell}^{N_0}}r_{\ell,\ell_1,\ell_2,\ell_3}S_{[\frac{s}{\tau}],\ell,\ell_1,\ell_2,\ell_3}\\
       &\hspace{1cm}\times \big(\widehat{\eta^{\varepsilon,R}}_{\ell_1}([s/\tau]\tau)\big)*\widehat{\eta^{\varepsilon,R}}_{\ell_2}([s/\tau]\tau)\widehat{\eta^{\varepsilon,R}}_{\ell_3}([s/\tau]\tau) D\theta_R(\eta^{\varepsilon,R}([s/\tau]\tau)){{\rm e}}^{-{\bf i}s\Delta}Q^{\frac{1}{2}}\,dW(s)\big|^2\big]\\
       &\lesssim  \varepsilon^6 \mathbb{E}\big[\sum_{\ell\in\mathcal{T}_{N_0}}(1+|\ell|^2) \sum_{k=0}^{N_\varepsilon-2}\int_{t_k}^{t_{k+1}}
       \|\sum_{(\ell_1,\ell_2,\ell_3)\in\mathcal{I}_{\ell}^{N_0}} r_{\ell,\ell_1,\ell_2,\ell_3}S_{k,\ell,\ell_1,\ell_2,\ell_3}\\
       &\hspace{1cm}\times \big(\widehat{\eta^{\varepsilon,R}}_{\ell_1}(t_k)\big)*\widehat{\eta^{\varepsilon,R}}_{\ell_2}(t_k)\widehat{\eta^{\varepsilon,R}}_{\ell_3}(t_k) D\theta_R(\eta^{\varepsilon,R}(t_k)){{\rm e}}^{-{\bf i}s\Delta}Q^{\frac{1}{2}} \|^2_{\mathcal{L}_2(L^2;\mathbb{C})}\,ds\big]\\
       &\lesssim  \varepsilon^6 \mathbb{E}\big[\sum_{\ell\in\mathcal{T}_{N_0}} \sum_{k=0}^{N_\varepsilon-2}\int_{t_k}^{t_{k+1}}
       \big(\sum_{(\ell_1,\ell_2,\ell_3)\in\mathcal{I}_{\ell}^{N_0}}\| (1+|\ell|) r_{\ell,\ell_1,\ell_2,\ell_3}S_{k,\ell,\ell_1,\ell_2,\ell_3}\\
       &\hspace{1cm}\times \big(\widehat{\eta^{\varepsilon,R}}_{\ell_1}(t_k)\big)*\widehat{\eta^{\varepsilon,R}}_{\ell_2}(t_k)\widehat{\eta^{\varepsilon,R}}_{\ell_3}(t_k) D\theta_R(\eta^{\varepsilon,R}(t_k)){{\rm e}}^{-{\bf i}s\Delta}Q^{\frac{1}{2}}\|_{\mathcal{L}_2(L^2;\mathbb{C})} \big)^2\,ds\big].
   \end{align*}
   Now, using the bound in \eqref{b_fe0} and considering an orthonormal basis of $L^2(\mathbb{T})$ $\{e_m\}_{m\in\mathbb{Z}}$, we compute
   \begin{align*}
       &\| (1+|\ell|) r_{\ell,\ell_1,\ell_2,\ell_3}S_{k,\ell,\ell_1,\ell_2,\ell_3}
        \big(\widehat{\eta^{\varepsilon,R}}_{\ell_1}(t_k)\big)^*\widehat{\eta^{\varepsilon,R}}_{\ell_2}(t_k)\widehat{\eta^{\varepsilon,R}}_{\ell_3}(t_k) D\theta_R(\eta^{\varepsilon,R}(t_k)){{\rm e}}^{-{\bf i}s\Delta}Q^{\frac{1}{2}}\|_{\mathcal{L}_2(L^2;\mathbb{C})}^2\\
        &\quad = \sum_{m\in\mathbb{Z}} \big|(1+|\ell|) r_{\ell,\ell_1,\ell_2,\ell_3}S_{k,\ell,\ell_1,\ell_2,\ell_3}
       \big(\widehat{\eta^{\varepsilon,R}}_{\ell_1}(t_k)\big)^*\widehat{\eta^{\varepsilon,R}}_{\ell_2}(t_k)\widehat{\eta^{\varepsilon,R}}_{\ell_3}(t_k) D\theta_R(\eta^{\varepsilon,R}(t_k)){{\rm e}}^{-{\bf i}s\Delta}Q^{\frac{1}{2}}e_m\big|^2\\
        &\quad\lesssim \tau^2 \big(\prod_{j=1}^{3} (1+|\ell_j|)\big)^2 \big|\widehat{\eta^{\varepsilon,R}}_{\ell_1}(t_k)\big|^2 \big|\widehat{\eta^{\varepsilon,R}}_{\ell_2}(t_k)\big|^2 \big|\widehat{\eta^{\varepsilon,R}}_{\ell_3}(t_k)\big|^2 
        \sum_{m\in\mathbb{Z}} \big|D\theta_R(\eta^{\varepsilon,R}(t_k)){{\rm e}}^{-{\bf i}s\Delta}Q^{\frac{1}{2}}e_m\big|^2.
   \end{align*}
   Using the fact that ${\rm e}^{-{\bf i}s \Delta}Q^{\frac{1}{2}}e_m \in H^1$, the isometry property of the operator ${\rm e}^{-{\bf i}s\Delta}$ in $H^1$ and the uniform boundedness of $D\theta_R(u) :H^1\to \mathbb{R}, \ u\in H^1$, (see, e.g., \cite{deb0}), we finally get
   \begin{align*}
       &\| (1+|\ell|) r_{\ell,\ell_1,\ell_2,\ell_3}S_{k,\ell,\ell_1,\ell_2,\ell_3}
        \big(\widehat{\eta^{\varepsilon,R}}_{\ell_1}(t_k)\big)^*\widehat{\eta^{\varepsilon,R}}_{\ell_2}(t_k)\widehat{\eta^{\varepsilon,R}}_{\ell_3}(t_k) D\theta_R(\eta^{\varepsilon,R}(t_k)){{\rm e}}^{-{\bf i}s\Delta}Q^{\frac{1}{2}}\|_{\mathcal{L}_2(L^2;\mathbb{C})}\\
        &\quad \lesssim \tau \prod_{j=1}^{3} (1+|\ell_j|) \big|\widehat{\eta^{\varepsilon,R}}_{\ell_1}(t_k)\big| \big|\widehat{\eta^{\varepsilon,R}}_{\ell_2}(t_k)\big| \big|\widehat{\eta^{\varepsilon,R}}_{\ell_3}(t_k)\big| \big(\sum_{m\in\mathbb{Z}} \|Q^{\frac{1}{2}}e_m\|^2_{H^1}\big)^{\frac{1}{2}}\\
        &\quad = \tau \prod_{j=1}^{3} (1+|\ell_j|) \big|\widehat{\eta^{\varepsilon,R}}_{\ell_1}(t_k)\big| \big|\widehat{\eta^{\varepsilon,R}}_{\ell_2}(t_k)\big| \big|\widehat{\eta^{\varepsilon,R}}_{\ell_3}(t_k) \|Q^{\frac{1}{2}}\|_{\mathcal{L}_2^1}.
   \end{align*}
   Then, we end up with
   \begin{align*}
        \mathbb{E}\big[\max_{1\le n\le N_\varepsilon}\|\mathcal{I}^{n-1}_{4,4,b}\|^2_{H^1}\big]&\lesssim\varepsilon^6 \tau^3 \sum_{k=0}^{N_\varepsilon-2}\mathbb{E}\big[\sum_{l\in\mathcal{T}_{N_0}}\\
        &\hspace{1cm}\times \big(\sum_{(\ell_1,\ell_2,\ell_3)\in\mathcal{I}_{\ell}^{N_0}} \prod_{j=1}^{3} (1+|\ell_j|)\big|\widehat{\eta^{\varepsilon,R}}_{\ell_1}(t_k)\big| \big|\widehat{\eta^{\varepsilon,R}}_{\ell_2}(t_k)\big| \big|\widehat{\eta^{\varepsilon,R}}_{\ell_3}(t_k)\big| \big)^2\big]\\
        &\lesssim \varepsilon^6 \tau^3 \sum_{k=0}^{N_\varepsilon-2}\mathbb{E}\big[\|\eta^{\varepsilon,R}(t_k)\|^6_{H^2}\big]\\
        &\lesssim (\varepsilon^2 \tau)^2.
   \end{align*}
   Putting all together, we achieve 
   $$
\mathbb{E}\big[\max_{1\le n\le N_\varepsilon}\|\mathcal{I}^{n-1}_{4,4}\|^2_{H^1}\big]\lesssim(\varepsilon^2\tau)^2,
   $$
   that finally gives
   \begin{equation}
       \label{fin_es0}
   \mathbb{E}\big[\max_{1\le n\le N_\varepsilon}\|\sum_{k=0}^{n-1}\mathcal{P}_{N_0}\mathcal{R}^{\varepsilon,R}(\eta^{\varepsilon,R}(t_k))\|^2_{H^1}\big]\lesssim(\varepsilon^2 \tau)^2.
    \end{equation}
    Inserting \eqref{fin_es0} into \eqref{in_es0}, we get
    $$
    \mathbb{E}\big[\max_{1\le n\le N_\varepsilon}\|e^{\varepsilon,R}_{n}\|_{H^1}^2\big]\lesssim \tau_0^{2(\sigma-1)}+ (\varepsilon^2\tau)^2+\varepsilon^2 \tau\sum_{k=0}^{N_\varepsilon-1}\mathbb{E}\big[\max_{0\le l \le k}\|e^{\varepsilon,R}_{l}\|_{H^1}^2\big].
    $$
    An application of the discrete Gronwall lemma then gives the estimate in \eqref{long_term_mean_trunc_cub}. 

    In order to obtain \eqref{long_term_path_trunc_cub}, note that, when $\tau_0^{\sigma-1}$ is neglectable, for $\tau>0$ and $\delta\in(0,1)$, the Markov inequality and \eqref{long_term_mean_trunc_cub} yield
    \begin{align*}
    \mathbb{P}\big(\max_{1\le n\le N_\varepsilon}\|\psi^{\varepsilon,R}(t_n)-\psi_n^{\varepsilon,R}\|_{H^1}\ge \varepsilon^2 \tau^\delta\big)
    &\le \frac{1}{\varepsilon^4 \tau^{2\delta}}\mathbb{E}\big[\max_{1\le n\le N_\varepsilon}\|\psi^{\varepsilon,R}(t_n)-\psi_n^{\varepsilon,R}\|^2_{H^1}\big]\\
    &\lesssim \tau^{2(1-\delta)}.
   \end{align*}
   Choosing $\tau_k=q^{-k}$ for some $q>1$, one has
   $$
   \sum_{k=0}^\infty \mathbb{P}\big(\max_{1\le n\le N^k_\varepsilon}\|\psi^{\varepsilon,R}(t_n)-\psi_n^{\varepsilon,R}\|_{H^1}\ge \varepsilon^2 \tau_k^\delta\big)
   \lesssim \sum_{k=0}^\infty q^{-2k(1-\delta)}<\infty
   $$
   independently on $\varepsilon$. The result \eqref{long_term_path_trunc_cub} then follows by an application of the Borel-Cantelli lemma.
\end{proof}

\begin{remark}
    According to Lemma \ref{long_term_cub_trunc}, in case of smooth solutions, i.e., $\sigma$ large enough, one then gets that the method (SNRLR2), when applied to the truncated cubic equation \eqref{trunc_cub2}, shows an error of size $\mathcal{O}(\varepsilon^2\tau)$ in $L^2(\Omega; H^1)$ up to large time scales of length $\mathcal{O}(\varepsilon^{-2})$. Moreover, almost surely, one has that the $H^1$ error remains bounded and of size $\mathcal{O}(\varepsilon^2\tau^{\delta})$ on the same time scale, with any $\delta<1$, if $\tau$ is sufficiently small but independent on $\varepsilon$.
\end{remark}

%We also need the following lemma.
%\begin{lemma}
 %   \label{path_stab_cub}
  %  Let us assume $\|\psi^\varepsilon_n\|_{H^1}\le R_1(\omega)$, $\|\mathcal{W}_n\|_{H^1}\le R_2(\omega)$. Then, there exists a positive constant $C_{GN}$ such that, for all $\tau\le % C_{GN}^{-1}$, one has $\|\psi^{\varepsilon}_{n+1}\|_{H^1}\le \tilde{R}(\omega)$, where $\tilde{R}(\omega):=R_1(\omega)(1+R_1(\omega)^2)+R_2(\omega)$.
% \end{lemma}
With Lemma \ref{long_term_cub_trunc} at hand, we are now in a position to state and prove a pathwise long-term error bound for the numerical integrator (SNRLR2) applied to the cubic stochastic Schrödinger equation \eqref{cub_small}.
\begin{theorem}
\label{long_term_path_orig_cub}
    Let us assume the setting of Lemma \ref{long_term_cub_trunc}, under the condition that the term $\tau_0^{\sigma-1}$ is neglectable. Then, the numerical integrator (SNRLR2), applied to \eqref{cub_small} satisfies 
    $$
    \max_{1\le n \le N_\varepsilon}\| \psi^\varepsilon(t_n)-\psi^\varepsilon_n\|_{H^1}<\varepsilon^2\tau^\delta, \qquad \mathbb{P}-{\rm a.s.}, \quad N_\varepsilon=\frac{T}{\varepsilon^2 \tau},
    $$
    for any $\delta<1$, asymptotically in $\tau$.
\end{theorem}
\begin{proof}
    Similarly to the proof of Proposition \ref{conv_path1}, let us denote
$$
R_{1,\varepsilon}(\omega)\ge \sup_{t\in[0,T_\varepsilon]} \|u(t)\|_{H^1}+\max_{1\le n\le N_\varepsilon} \|\mathcal{W}_n\|_{H^1},
$$
where we have explicit the fact that, in principle, one cannot assume that the random variable $R_{1,\varepsilon}$ is independent from $\varepsilon$.

     Now, let $\tau$ sufficiently small (independently on $\varepsilon$) and assume there exists $1\le \tilde{n}\le N_\varepsilon$ such that one has
     $$
         \|\psi^\varepsilon (t_n)-\psi^\varepsilon_{n}\|_{H^1}<\varepsilon^2 \tau^\delta, \ n=1,\dots,\tilde{n}-1, \quad \|\psi^\varepsilon (t_{\tilde{n}})-\psi^\varepsilon_{\tilde{n}}\|_{H^1}\ge \varepsilon^2 \tau^\delta.
     $$
     Then, for $1\le n\le \tilde{n}-1$ and considering $\varepsilon^2 \tau^\delta<1$, one has
     \begin{align*}
         \|\psi^\varepsilon_n\|_{H^1}&\le \|\psi^\varepsilon(t_n)\|_{H^1}+\|\psi^\varepsilon(t_n)-\psi^\varepsilon_n\|_{H^1}\\
         &\le \|\psi^\varepsilon(t_n)\|_{H^1}+\varepsilon^2 \tau^\delta\\
         &\le R_{1,\varepsilon}(w)+1.
     \end{align*}
     Then, using this, the triangular inequality, the unitary property of the operator ${\rm e}^{{\bf i}t\Delta}$ and the algebra property of $H^\sigma$, one has
     \begin{align*}
         \|\psi^\varepsilon_{\tilde{n}}\|_{H^1}\le R_{1,\varepsilon}(w)+1&+\varepsilon^2\tau C_{GN} \|\psi^\varepsilon_{\tilde{n}-1}\|^3_{H^1}
         +\varepsilon^2 \tau \big\| (\widehat{g(\psi^\varepsilon_{\tilde{n}-1})}\big)_0 {\rm e}^{{\bf i} \tau \Delta} \psi^\varepsilon_{\tilde{n}-1}\big\|_{H^1}\\
         &+\varepsilon^2 \tau \|{\rm e}^{{\bf i}\tau \Delta} h(\psi^\varepsilon_{\tilde{n}-1})\|_{H^1}+\varepsilon R_{1,\varepsilon}(w),
     \end{align*}
     for some $C_{GN}>0$.
% \begin{align*}
 % \| \psi_{n+1}\|_{H^\sigma} \le 2 R_0&+ \tau \|(\psi^n)^2(\varphi_1(-2 {\bf i} \tau \Delta)\overline{\psi_n})\|_{H^\sigma}\\
% &\le 2 R_0+C \tau \|\psi_n\|_{H^\sigma}^3\\
% &\le 2 R_0+C \tau R_0^2\|\psi_n\|_{H^\sigma}.
%+2 \tau \| (\widehat{g(u^n)}\big)_0 {\rm e}^{{\bf i}\tau \Delta } u^n\|_{H^\sigma}
%+\tau\| h(u^n)\|_{H^\sigma}.
%\end{align*}

By definition of $g$ and using the Holder's inequality, we have
\begin{align*}
|\big(\widehat{g(\psi^\varepsilon_{\tilde{n}-1})}\big)_0|&\le \displaystyle\frac{1}{2\pi}\big[\left\|\psi_{\tilde{n}-1}^\varepsilon\right\|^2+\int_{-\pi}^{\pi}\left|\psi^\varepsilon_{\tilde{n}-1}(x)\right|\left|\varphi(-2{\bf i}\tau \Delta)\overline{\psi_{\tilde{n}-1}^\varepsilon}(x)\right| \,{\rm d}x\big]\\
&\le \displaystyle\frac{1}{2\pi}\big[\left\|\psi_{\tilde{n}-1}^\varepsilon\right\|^2+\big(\int_{-\pi}^{\pi}\left| \psi_{\tilde{n}-1}^\varepsilon(x)\right|^2\,{\rm d}x\big)^{\frac{1}{2}}
\big(\int_{-\pi}^{\pi}\left|\varphi(-2{\bf i}\tau \Delta)\overline{\psi_{\tilde{n}-1}^\varepsilon}(x)\right|^2\,{\rm d}x\big)^{\frac{1}{2}}\big]\\
&=\displaystyle\frac{1}{2\pi}\left[\left\|\psi^\varepsilon_{\tilde{n}-1}\right\|^2_{H^1}+\left\|\psi^\varepsilon_{\tilde{n}-1}\right\|_{H^1}\left\|\varphi(-2{\bf i}\tau \Delta) \overline{\psi_{\tilde{n}-1}^\varepsilon}\right\|_{H^1}\right]\\
&\lesssim  \left\| \psi_{\tilde{n}-1}^\varepsilon\right\|_{H^1}^2.
\end{align*}
This then yields
$$
\big\| (\widehat{g(\psi_{\tilde{n}-1}^\varepsilon)}\big)_0 {\rm e}^{{\bf i} \tau \Delta} \psi_{\tilde{n}-1}^\varepsilon\big\|_{H1}\le |  (\widehat{g(\psi_{\tilde{n}-1}^\varepsilon)}\big)_0|
\| {\rm e}^{{\bf i} \tau \Delta} \psi^\varepsilon_{\tilde{n}-1}\|_{H^1}\lesssim  \|\psi^\varepsilon_{\tilde{n}-1} \|^3_{H^1}.
$$
It remains to bound the term $\left\| h(\psi^\varepsilon_{\tilde{n}-1})\right\|_{H^1}$. First note that for $w \in H^1$, we have
$$
\left| \varphi_1\left(2 {\bf i} \tau \ell^2\right) \hat{w}_\ell \right|\le C \left|\hat{w}_\ell \right|, \quad \ell \in \mathbb{Z}.
$$
Then, we get
\begin{align*}
\| h(\psi^\varepsilon_{\tilde{n}-1})\|_{H^1}^2 &= \displaystyle\sum_{\ell \in\mathbb{Z}}\left(1+\left|\ell\right|\right)^{2}| (\widehat{h(\psi^\varepsilon_{\tilde{n}-1})})_\ell |^2\\
&\le  \displaystyle\sum_{\ell \in\mathbb{Z}}\left(1+\left|\ell \right|\right)^{2} \big|\big( \widehat{\psi^\varepsilon}_{\tilde{n}-1,\ell} \big)^*\widehat{\psi^\varepsilon}_{\tilde{n}-1,\ell}\widehat{\psi^\varepsilon}_{\tilde{n}-1,\ell} \big|^2\\
&\hspace{2cm}+\displaystyle\sum_{\ell \in\mathbb{Z}}\left(1+\left|\ell \right|\right)^{2} \big| \varphi\left(2 {\bf i} \tau \ell ^2\right)  \big( \widehat{\psi^\varepsilon}_{\tilde{n}-1,\ell} \big)^*\widehat{\psi^\varepsilon}_{\tilde{n}-1,\ell}\widehat{\psi^\varepsilon}_{\tilde{n}-1,\ell}\big|^2\\
&\lesssim \displaystyle\sum_{\ell\in\mathbb{Z}}\left(1+\left|\ell\right|\right)^{2} \big| \big( \widehat{\psi^\varepsilon}_{\tilde{n}-1,\ell} \big)^*\widehat{\psi^\varepsilon}_{\tilde{n}-1,\ell}\widehat{\psi^\varepsilon}_{\tilde{n}-1,\ell}\big|^2\\
&\lesssim \|\psi^\varepsilon_{\tilde{n}-1}\|^3_{H^1}.
\end{align*}
Then, we end up with
\begin{align*}
    \|\psi^\varepsilon_{\tilde{n}}\|_{H^1}&\le (1+\varepsilon)R_{1,\varepsilon}(w)+1+\varepsilon^2\tau C_{GN}  \|\psi^\varepsilon_{\tilde{n}-1}\|^3_{H^1}\\
    &\le  (1+\varepsilon)R_{1,\varepsilon}(w)+1+\varepsilon^2\tau C_{GN} (R_{1,\varepsilon}(w)+1)^3.
\end{align*}
This shows that for sufficiently small $\tau$ (independently on $\varepsilon$), there exists a random variable $R_{2,\varepsilon}(\omega)\ge R_{1,\varepsilon}(\omega)+1$ such that
$$
\|\psi^\varepsilon_{\tilde{n}}\|_{H^1}\le R_{2,\varepsilon}(\omega).
$$
This shows that $\psi^\varepsilon_n=\psi_n^{\varepsilon,R_{2,\varepsilon}(\omega)}$, for $n=1,\dots,\tilde{n}$ and 
$$
\max_{1\le n \le N_\varepsilon} \|\psi^{\varepsilon, R_{2,\varepsilon}(\omega)}(t_n)-\psi^{\varepsilon, R_{2,\varepsilon}(\omega)}_n\|\ge \varepsilon^2 \tau^\delta,
$$
provided $\tau$ sufficiently small. In particular, this means that one can always find a sufficiently large $k$ such that 
$$
\max_{1\le n \le N^k_\varepsilon} \|\psi^{\varepsilon, R_{2,\varepsilon}(\omega)}(t_n)-\psi^{\varepsilon, R_{2,\varepsilon}(\omega)}_n\|\ge \varepsilon^2 \tau_k^\delta,
$$
being $\{\tau_k\}$ the subsequence defined in Lemma \ref{long_term_cub_trunc}. By \eqref{long_term_path_trunc_cub}, this is impossible and then a contradiction holds. The proof is concluded.
\end{proof}
\section{Numerical experiments}
In this section, we give several numerical experiments to confirm all our theoretical results. 

Based on Assumption \ref{ass_q}, for $\sigma\ge 0$, we will consider the following
$$
q_0=1,\qquad q_{\ell}=|\ell|^{-(2\sigma+1)-\epsilon}, \ \ell\in\mathbb{Z}_{\ne 0}.
$$
Then, $Q^{\frac{1}{2}}\in\mathcal{L}_2^\sigma$. 

Moreover, in all the experiments, we will adopt a pseudo-spectral approximation in space with $\kappa=2^{11}$ Fourier modes, i.e., $\ell= -\frac{\kappa}{2}+1,\dots,\frac{\kappa}{2}$. Following \cite{os_ka,fe_ka}, w we will consider the following initial data. First, we take 
$$
U_0={\rm rand} +{\bf i} \ {\rm rand}, \quad U_\ell=|\ell|^{-(\frac{1}{2}+\epsilon)}\big({\rm rand}+{\bf i} \ {\rm rand}\big), \qquad \ell=-k/2+1,\dots,k/2, \ \ell\ne 0,
$$
for an arbitrary small $\epsilon$. Then, for some $\theta\ge 0$, we take $u_0$ given by the following Fourier coefficients 
$$
\hat{u_0}_{0}=U_0, \quad \hat{u_0}_{\ell} = |\ell|^{-\theta} U_\ell, \quad \ell=-\frac{\kappa}{2}+1,\dots,\frac{\kappa}{2}, \ \ell\ne 0.
$$
Then, $u_0 \in H^\theta$ almost surely (see, e.g., \cite{fe_ka}). Note that we fix $u_0$ independently from the realizations, so that it is deterministic w.r.t. the $\sigma$-algebra generated by $\{W(t)\}_{t\ge 0}$.

%\subsection{Linear SSE with potential}
We start by testing strong and pathwise convergence rates for method (SLR1) and method (SNRLR1), applied to the case of the linear SSE with rough potential, with $\varepsilon=1$ and, hence, $\Phi(x)=V(x)$, i.e., Equation \eqref{lin_eq}.
\begin{comment}
We test the integrator (SLR1) in \eqref{stand_lr_met} that reads
\begin{equation}
\label{lr_lin_met}
u_{n+1}={\rm e}^{{\bf i}\tau\Delta}\big[u_{n}-{\bf i}\tau u_n \big(\varphi_1(-{\bf i}\tau \Delta)V(x)\big)\big]+\mathcal{W}_n,
\end{equation}
and the integrator (SNRLR1) in \eqref{non_res1}
\begin{align}
u^\varepsilon_{n+1}={\rm e}^{{\bf i}\tau\Delta} u_n&-{\bf i} \tau {\rm e}^{{\bf i}\tau \Delta} \big(u_n\varphi_1(-{\bf i}\tau\Delta)V(x)\big)\notag \\
&-{\bf i} \tau (1-\varphi_1(-4{\bf i}\tau\Delta)){\rm e}^{{\bf i}\tau\Delta}g_2(V(x),u_n)+\mathcal{W}_n \label{lin_met_nr}.
\end{align}
\end{comment}
For $\mu\ge 0$, we take the potential $V$ defined through its Fourier coefficients as follows
\begin{align*}
  &  V_0={\rm rand}+{\bf i}\ {\rm rand}, \quad V_{\ell}=|\ell|^{-(\frac{1}{2}+\epsilon)}\big({\rm rand}+{\bf i}\ {\rm rand}\big), \ \ell=1,\dots,\frac{\kappa}{2},\\
  &  \hat{V}_{0}=V_{0}, \quad \hat{V}_{\ell}=|\ell|^{-\mu} V_{\ell}, \ \ell=1,\dots,\frac{\kappa}{2},
  \quad \hat{V}_{\ell}=\overline{\hat{V}_{-\ell}}, \ \ell=-\frac{\kappa}{2},\dots,-1.
\end{align*}
Then, one has $V\in H^\mu$ almost surely.

\begin{figure}
\centering
\subfigure{\includegraphics[width=0.48\textwidth]{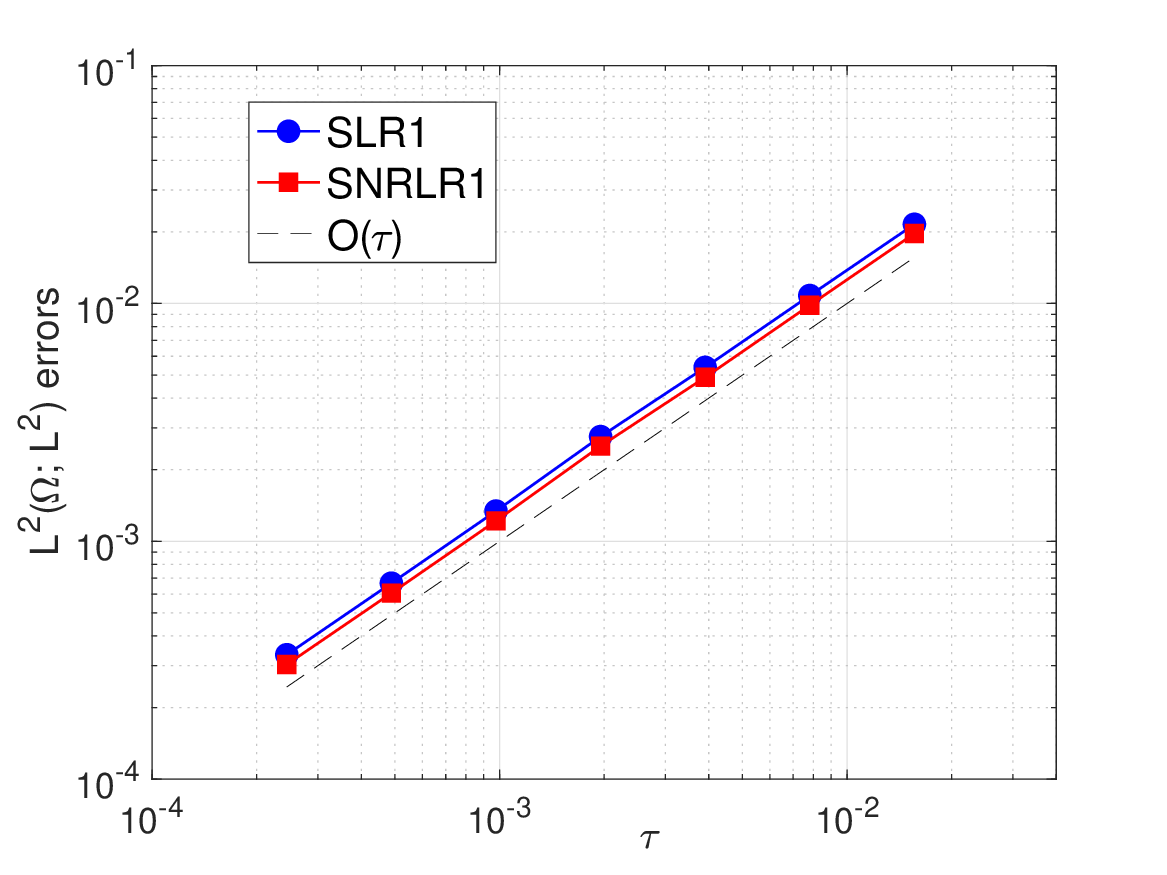}}\quad \subfigure{\includegraphics[width=0.48\textwidth]{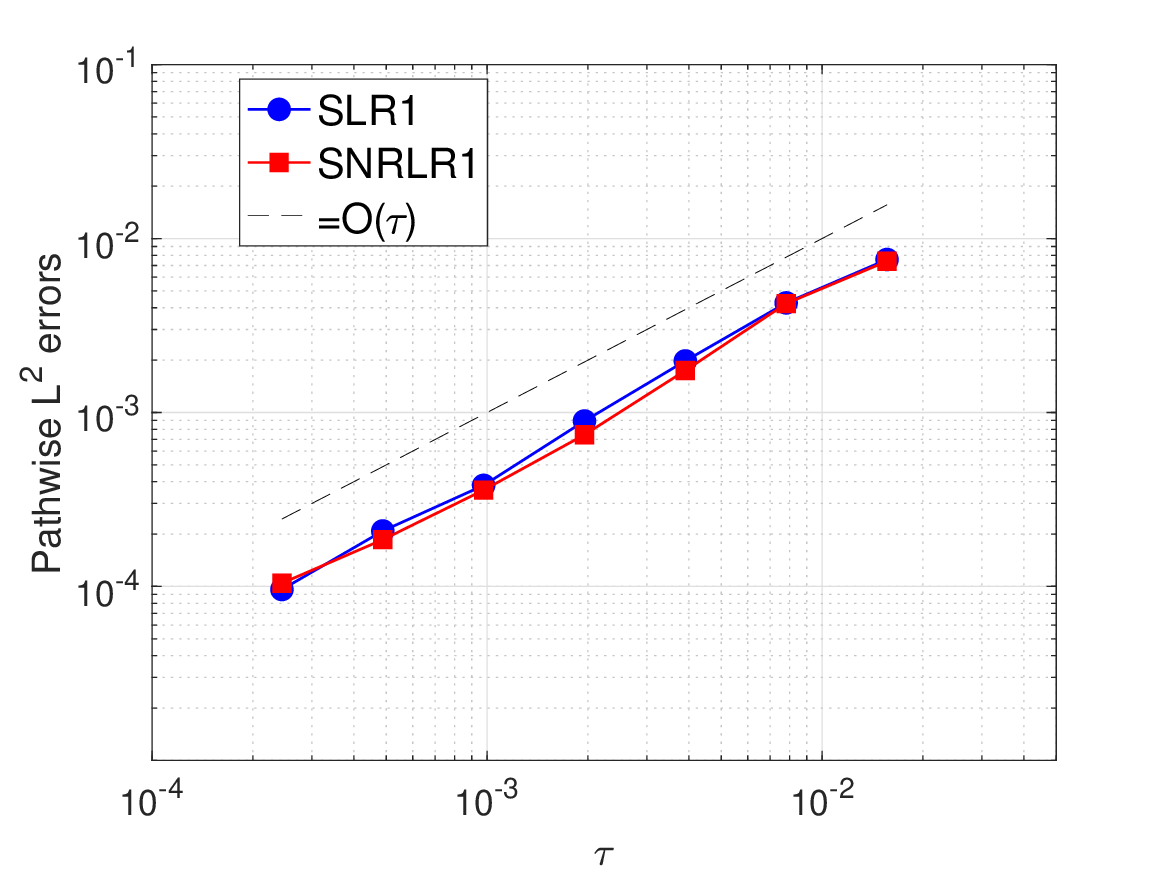}}
          \caption{$L^2(\Omega; L^2)$ errors (on the left) and pathwise $L^2$ errors (on the right) for methods (SLR1) and (SNRLR1) applied to Equation \eqref{lin_eq} in $T=1$, under $H^1$ regularity.}
           \label{fig1}
\end{figure}

The $L^2(\Omega; L^2)$ (on the left) and the pathwise $L^2$ (on the right) errors at $T=1$ for methods (SLR1) and (SNRLR1), applied to such equation, under $H^1$ regularity, are depicted in Figure \ref{fig1}. We take $\theta=\sigma=1$ and $\mu>\frac{3}{2}$.

The reference solution has been computed with a modification of the classical exponential method (EXP) defined in \cite{an_co}
\begin{equation}
    \label{mod_exp}
u_{n+1}={\rm e}^{{\bf i}\tau\Delta}[u_n-{\bf i}\tau V(x)u_n]+\mathcal{W}_n
\end{equation}
with step size $\tau_{\rm ref}=2^{-19}$. The errors have been computed with step sizes $\tau=2^j \tau_{\rm ref}, \ j=7,\dots,13$. The expectations have been approximated using $M=100$ sample paths.
 Strong order one has been confirmed for both methods.

\begin{figure}
\centering
\subfigure{\includegraphics[width=0.5\textwidth]{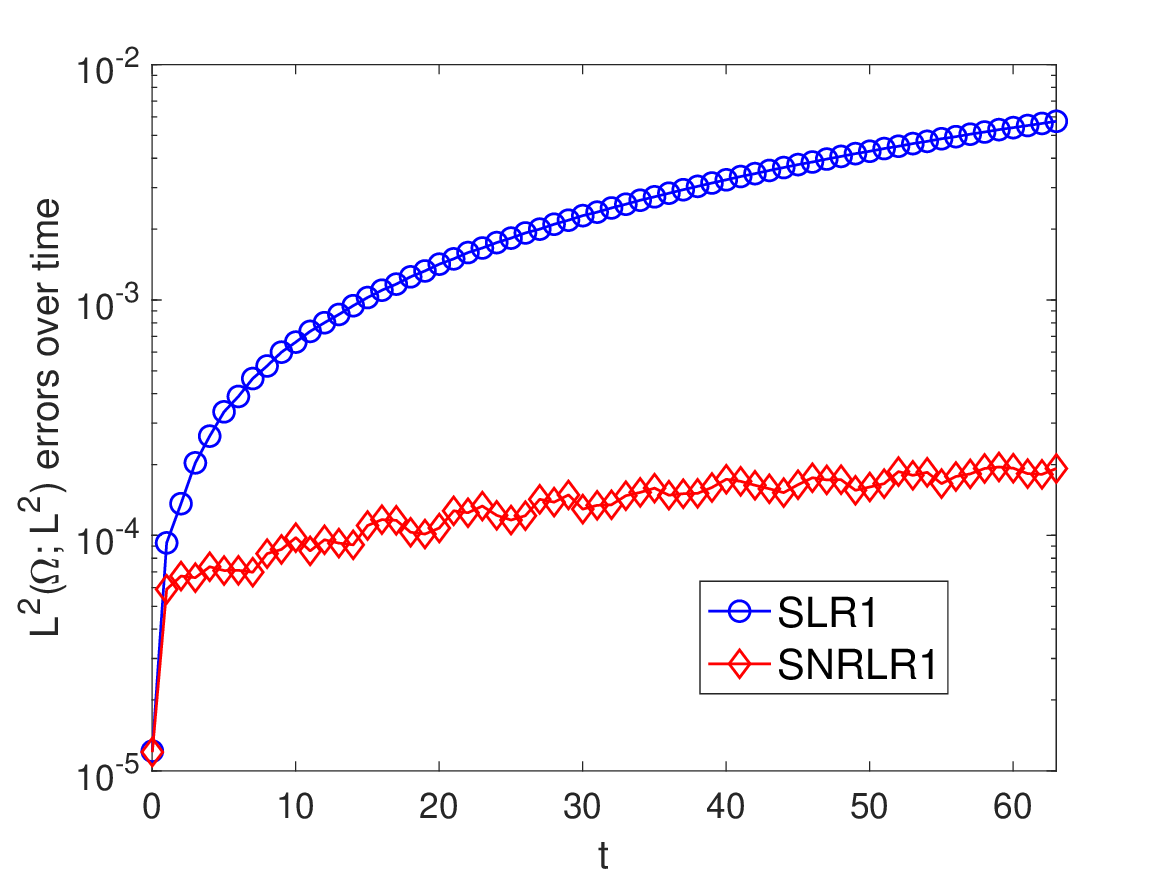}}
          \caption{Long-term $L^2(\Omega; L^2)$ errors under $H^1$ regularity, for methods (SLR1) and (SNRLR1), applied to \eqref{lin_small} with $\varepsilon=2^{-3}$.}
           \label{fig_mean1}
\end{figure}

In Figure \ref{fig_mean1}, we display the $L^2(\Omega; L^2)$ errors over time  for (SLR1) method and (SNRLR1) method, applied to the stochastic linear Schrödinger equations with small potential and noise (see \eqref{lin_small}), with $\theta=\sigma=1, \mu>\frac{3}{2}$, $T=1$, $\tau=2^{-6}$, for $\varepsilon=2^{-3}$. The reference solution has been computed with (SNRLR1) method, using a step size $\tau_{\rm ref}= 2^{-16}$. The expectations have been approximatex with $M=30$ sample paths. 

\begin{figure}
\centering
\subfigure{\includegraphics[width=0.48\textwidth]{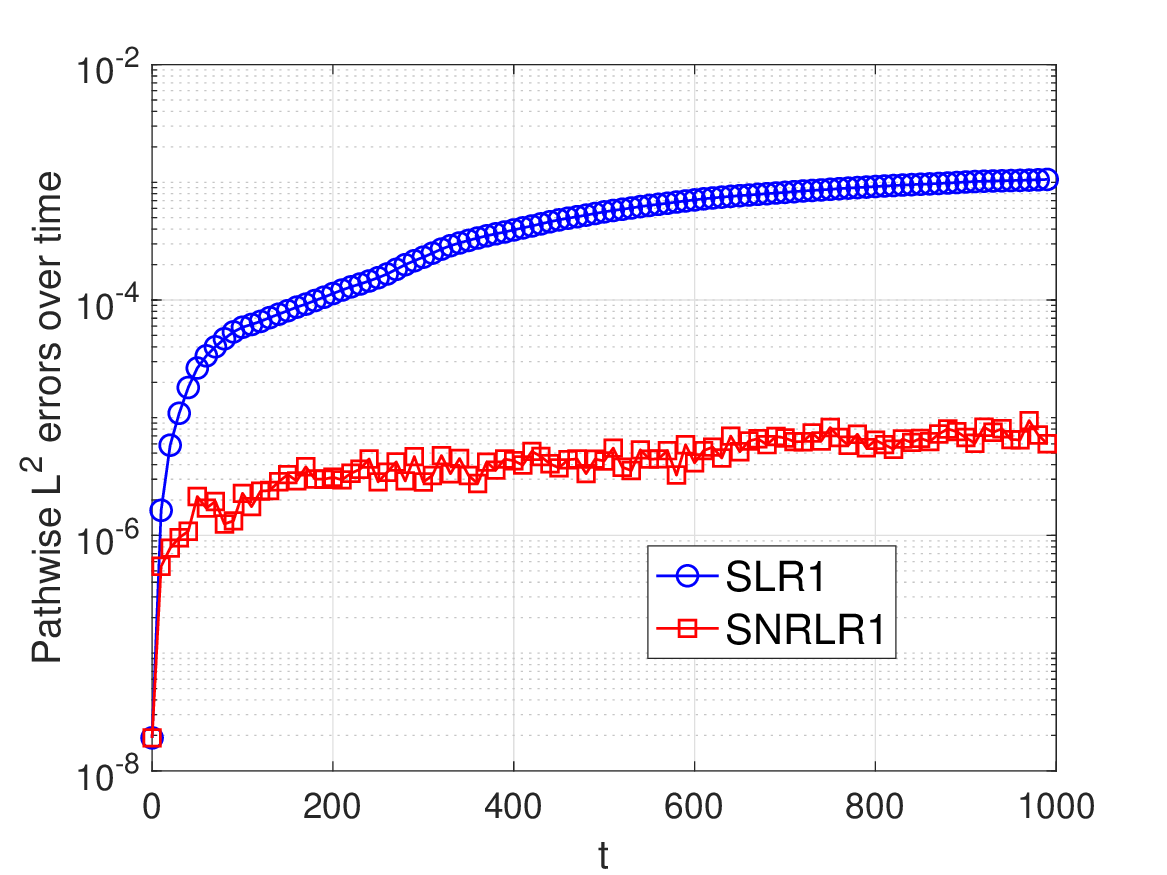}}\quad \subfigure{\includegraphics[width=0.48\textwidth]{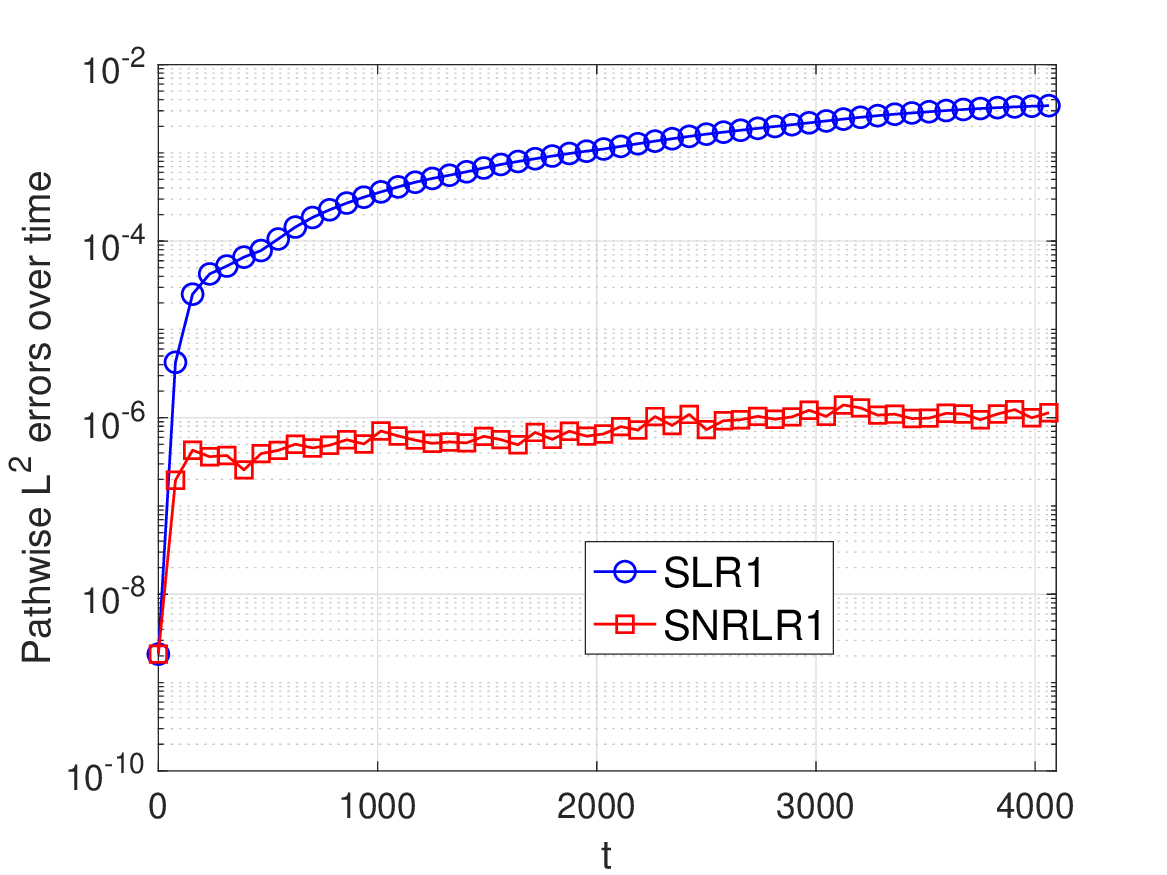}}
          \caption{Long-term pathwise $L^2$ errors under $H^1$ regularity, for methods (SLR1) and (SNRLR1), applied to \eqref{lin_small} with $\varepsilon=10^{-\frac{3}{2}}$ (on the left) and $\varepsilon=2^{-6}$ (on the right).}
           \label{fig3}
\end{figure}

Smaller values of $\varepsilon$, i.e., larger time windows, are considered in Figure \ref{fig3}, where we report the comparison between pathwise $L^2$ errors over long time for (SLR1) method and (SNRLR1) method, applied to the stochastic linear Schrödinger equations with small potential and noise (see \eqref{lin_small}), with $\theta=\sigma=1, \mu>\frac{3}{2}$, $T=1$, $\tau=0.01$ (for the left plot) and $\tau=2^{-6}$ (for the right plot), for $\varepsilon=10^{-\frac{3}{2}}$ (on the left) and $\varepsilon=2^{-6}$ (on the right). The reference solution has been computed with (SNRLR1) method, using a step size $\tau_{\rm ref}=10^{-5}$ (for the left plot) and $\tau_{\rm ref}=2^{-16}$ (for the right plot).

\begin{figure}
\centering
\subfigure{\includegraphics[width=0.5\textwidth]{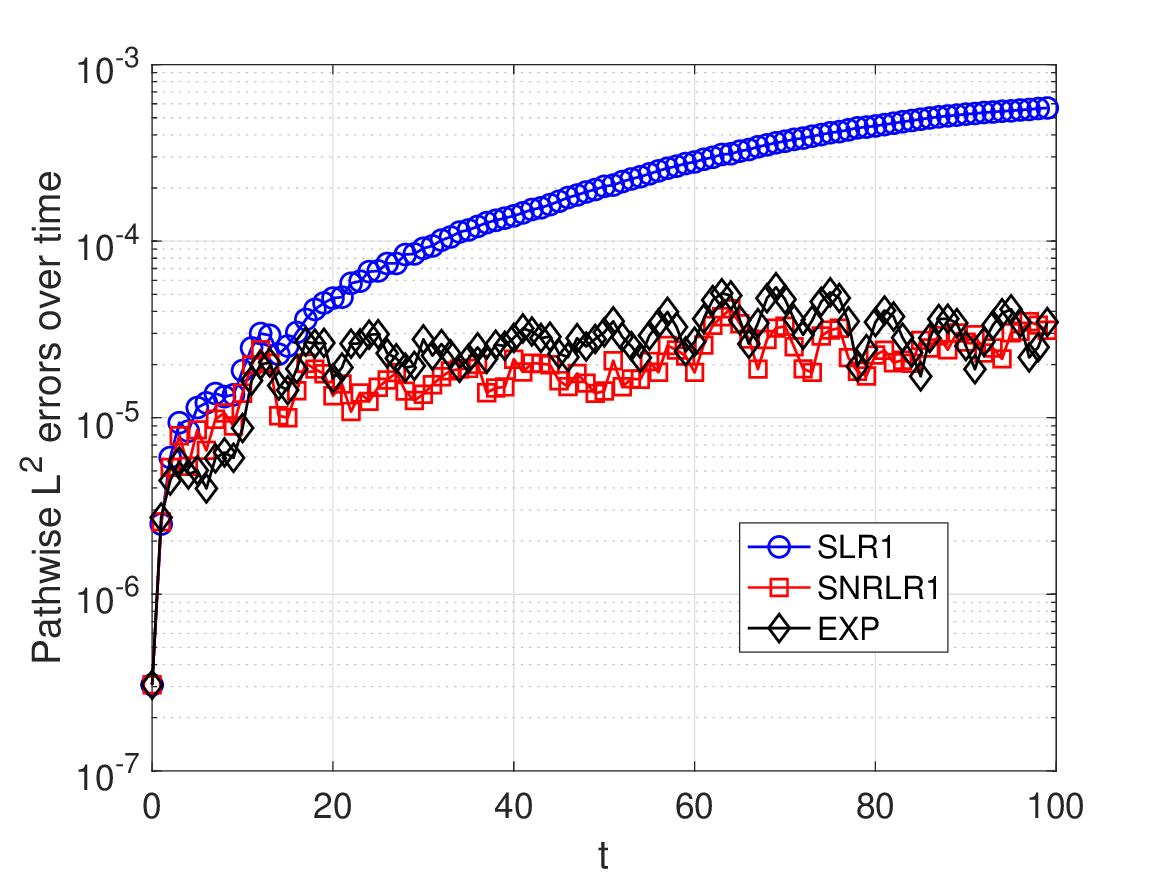}}%\quad \subfigure{\includegraphics[width=0.44\textwidth]{linear_H1_comp.eps}}
          \caption{Long-term pathwise $L^2$ errors under $H^2$ regularity, for methods (SLR1) and (SNRLR1), applied to \eqref{lin_small}, with $\varepsilon=0.01$.}
           \label{fig2}
\end{figure}

A long-term simulation under $H^2$ regularity has been done and depicted in Figure \ref{fig2}, where we report the comparison between pathwise $L^2$ errors over long time for (SLR1) method, (SNRLR1) method and the classical (EXP) method with exact simulation of the stochastic convoluiont (see \eqref{mod_exp}), applied to the same equation \eqref{lin_small}. We take $\theta=\sigma=\mu=2$, $\varepsilon=0.1$, $T=1$, i.e., $T_\varepsilon=100$ and $\tau=0.005$. The reference solution has been computed with (SNRLR1) method, using a step size $\tau_{\rm ref}=10^{-6}$.

In all the aforementioned plots, we can observe the boundedness of the errors for (SNRLR1) method, up time of length $\mathcal{O}(\varepsilon^{-2})$, predicted by Theorem \ref{lobg_term_res1}, and, hence, the improvement achieved with respect the method (SLR1). Also the classical (EXP) method \eqref{mod_exp} shows good behaviour in $H^2$ regularity when the stochastic convolution is simulated exactly.

\begin{figure}
\centering
\subfigure{\includegraphics[width=0.5\textwidth]{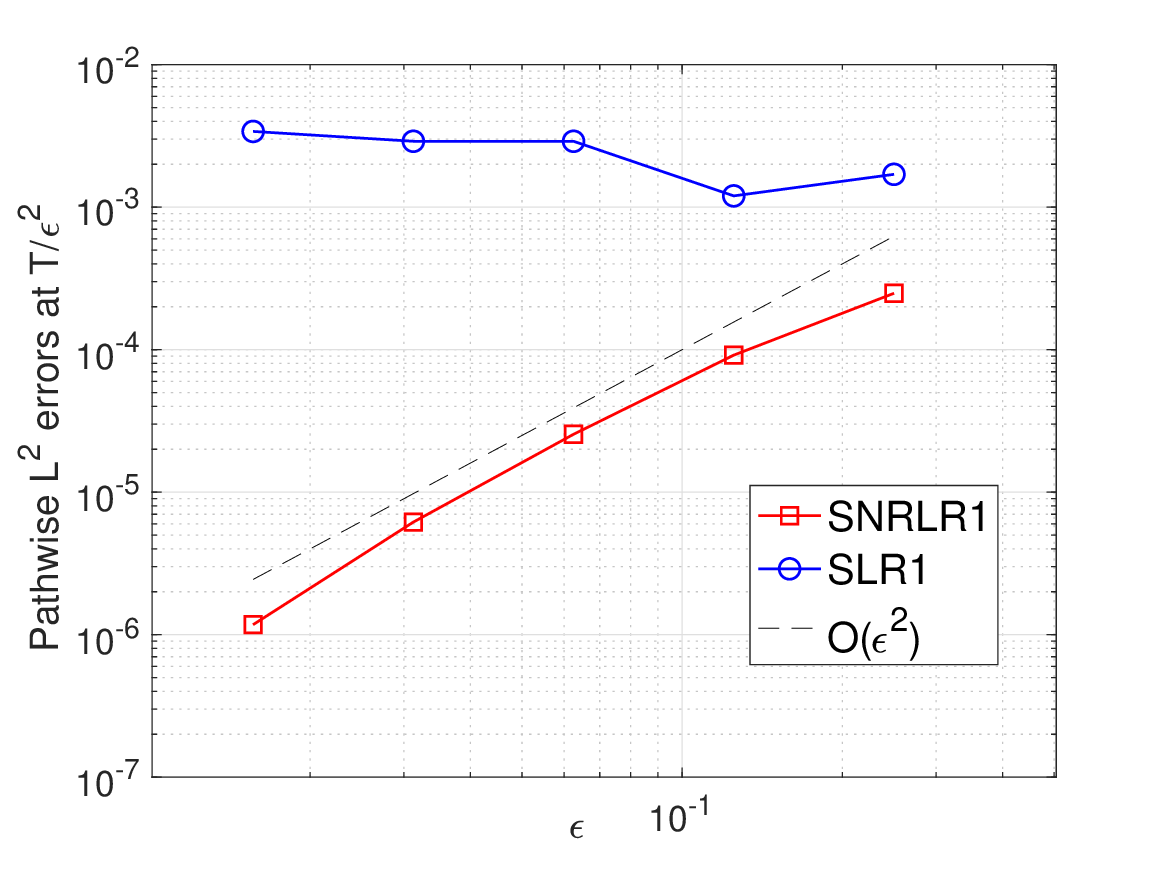}}
          \caption{Pathwise $L^2$ errors at $\frac{T}{\varepsilon^2}$, under $H^1$ regularity, over $\varepsilon$, for methods (SLR1) and (SNRLR1), applied to \eqref{lin_small} for $\varepsilon=2^{-j}, \ j=2,\dots,6$.}
           \label{fig4}
\end{figure}

In Figure \ref{fig4}, we display the pathwise convergence in $L^2$ with respect the parameter $\varepsilon$,  for methods (SLR1) and (SNRLR1), under $H^1$ regularity, i.e., we plot the $L^2$ errors at time $\frac{T}{\varepsilon^2}$, for $\varepsilon=2^{-j}$, for $j=2,\dots,6$, computed with $\tau=2^{-6}$. The reference solutions have been computed with method (SNRLR1) with $\tau_{\rm ref}=2^{-16}$. We can observe an $\mathcal{O}(\varepsilon^2)$ behaviour for our non-resonant integrator (SLR1), accordingly to the rexsult of Theorem \ref{lobg_term_res1}.

\begin{figure}
\centering
\subfigure{\includegraphics[width=0.5\textwidth]{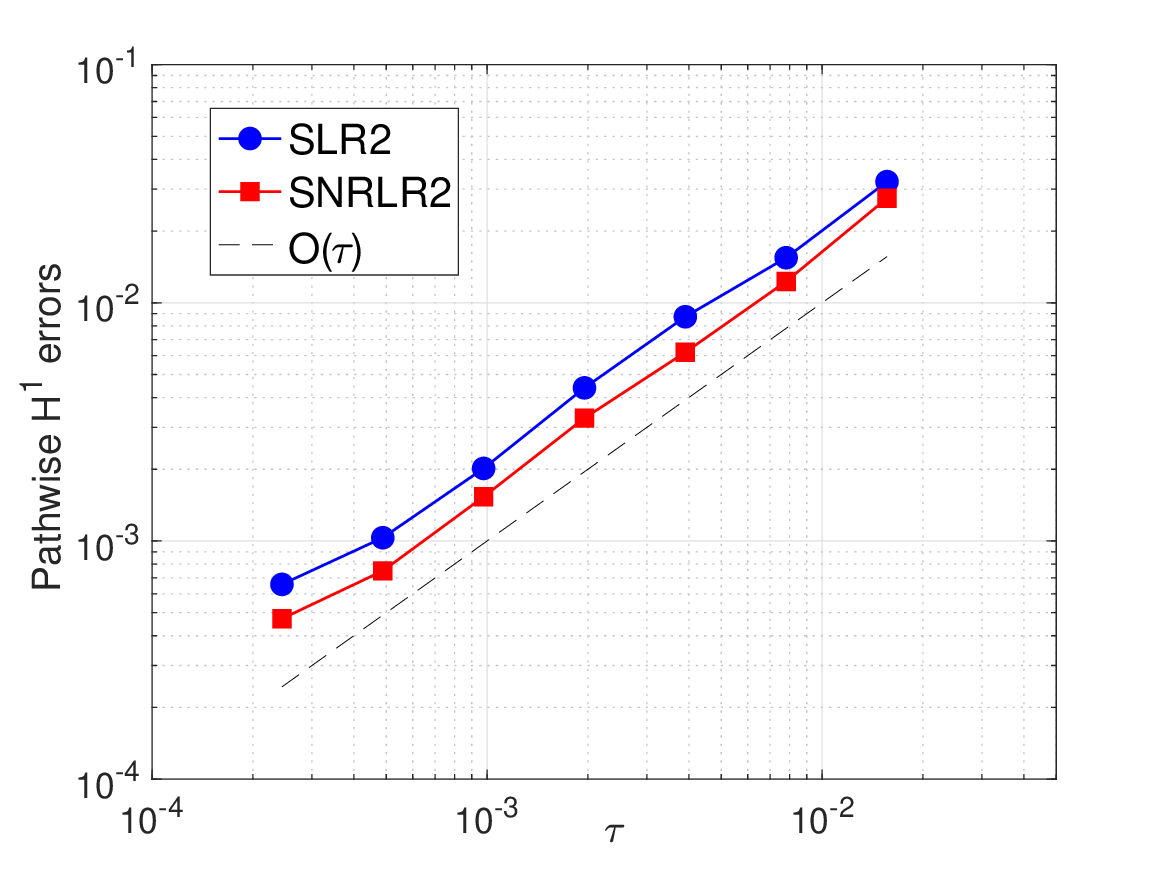}}
          \caption{Pathwise $H^1$ errors at $T=1$, under $H^2$ regularity, for methods (SLR2) and (SNRLR2), applied to the cubic equation \eqref{cub_eq}.}
           \label{fig_cub1}
\end{figure}

%\subsection{Cubic SSE}
Finally, we here consider the numerical approximation of the cubic SSE \eqref{cub_eq}-\eqref{cub_small} by methods (SLR2) and (SNRLR2). We test the pathwise convergence in $H^1$, under $H^2$ regularity for the exact solution, in $T=1$. As reference solution, we take the method (SNRLR2) itself with a very small step size $\tau_{\rm ref}=2^{-19}$ and we test $\tau_j=2^j \tau_{\rm ref}, \ j=7,\dots,13$. The results are depicted in Figure \ref{fig_cub1}, where order one of convergence is confirmed for both methods.

\end{document}